\documentclass{amsart}

\usepackage{amssymb}
\usepackage{amsfonts}
\usepackage{amscd}
\usepackage{lscape}
\usepackage{verbatim}
\usepackage{epsfig}
\usepackage{xr}

\input xy
\xyoption{all}

\DeclareFontFamily{OT1}{rsfs}{}
\DeclareFontShape{OT1}{rsfs}{n}{it}{<-> rsfs10}{}
\DeclareMathAlphabet{\mathscr}{OT1}{rsfs}{n}{it}

\newtheorem{thm}{Theorem}[subsection]
\newtheorem{lem}[thm]{Lemma}
\newtheorem{prop}[thm]{Proposition}

\newtheorem{cor}[thm]{Corollary}
\newtheorem{rems}[thm]{Remarks}

\newtheorem{sth}{Theorem}[thm]

\newtheorem{sprop}[sth]{Proposition}

\theoremstyle{definition}

  \newtheorem{defi}[thm]{Definition}
  \newtheorem{rem}[thm]{Remark}

  \newtheorem{ack}{Acknowledgements}  

\numberwithin{equation}{thm}

\newcommand{\Cref}[1]{Corollary~\textup{\ref{#1}}}
\newcommand{\Dref}[1]{Definition~\textup{\ref{#1}}}

\newcommand{\Lref}[1]{Lemma~\textup{\ref{#1}}}
\newcommand{\Pref}[1]{Proposition~\textup{\ref{#1}}}
\newcommand{\Rref}[1]{Remark~\textup{\ref{#1}}}
\newcommand{\Sref}[1]{Section~\textup{\ref{#1}}}
\newcommand{\Ssref}[1]{Subsection~\textup{\ref{#1}}}
\newcommand{\Tref}[1]{Theorem~\textup{\ref{#1}}}

\def\bilap#1{\hbox to 0pt{\hss#1\hss}}
  \def\Rarrow#1{\bilap{\hbox to#1{\rightarrowfill}}}
  \def\Larrow#1{\bilap{\hbox to#1{\leftarrowfill}}}
  \def\Equals#1{\bilap
                   {\raise 4pt\hbox
                     {\vrule width#1 height.5pt}%
                    \kern-#1\raise 1pt\hbox
                     {\vrule width#1 height.5pt}%
                   }}

\newcommand{\EQAL}[1]%
{\,\begin{picture}(#1,0)%
\put(0,3){\line(1,0){#1}}%
\put(0,1){\line(1,0){#1}}%
\end{picture}\,}%

\newcommand{\vlto}[1]%
{\,\begin{picture}(#1,3)%
\put(0,2){\vector(1,0){#1}}%
\end{picture}\,}%

\newcommand{\vllarrow}[1]%
{\,\begin{picture}(#1,3)%
\put(#1,2){\vector(-1,0){#1}}%
\end{picture}\,}%

\newcommand{\dirlm}[1]%
   {
      {\lim\hskip-1.58em\lower.65ex
        \hbox{$
                 {}_{\stackrel{\lower1ex\hbox
                                         {$\scriptstyle -\!\!\!\longrightarrow$}
                                       }{\vbox to0pt{\vss\vskip.6ex
                                             \hbox{$\scriptstyle{}^{#1}$}\vss}}
                    }
             $}
      }
\:}

\newcommand{\subdirlm}[1]%
   {
      {\lim\hskip-1.5em\lower.6ex
        \hbox{$
                    {}_{\stackrel{\lower1ex\hbox
                                            {$\scriptstyle\longrightarrow$}
                                 }{ ^{#1} }
                       }
              $}
      }
\:}

\newcommand{\inlm}[1]%
    {
       {\lim\hskip-1.58em\lower.65ex
         \hbox{$
                  {}_{\stackrel{\lower1ex\hbox
                                         {$\scriptstyle \longleftarrow\!\!\<-$}
                               }{\vbox to0pt{\vss\vskip.6ex
                                             \hbox{$\scriptstyle{}^{#1}$}\vss}}
                     }
              $}
       }
\:}

\def\hz#1{{\hbox to 0pt{#1}}}
\def\Iso{\vbox to 0pt{\vss\hbox{$\widetilde{\phantom{nn}}$}\vskip-7pt}}

\def\>{\mspace {1mu}}
\def\<{\mspace{-1mu}}
\def\({{\textup(}}
\def\){{\textup)}}
\def\bigl#1{{\textup{\begin{large}#1\end{large}}}}
\def\bigr#1{{\textup{\begin{large}#1\end{large}}}}

\newcommand{\btrg}{\blacktriangle}
\newcommand{\fm}{{\mathfrak{m}}}

\newcommand{\X}{{\mathscr X}}

\newcommand{\Y}{{\mathscr Y}}
\newcommand{\Z}{{\mathscr Z}}
\newcommand{\V}{{\mathscr V}}
\newcommand{\U}{{\mathscr U}}
\newcommand{\W}{{\mathscr W}}
\newcommand{\I}{{\mathscr I}}
\newcommand{\J}{{\mathscr J}}

\newcommand{\eE}{{\mathscr E}}
\newcommand{\eH}{{\mathscr H}}
\newcommand{\eP}{{\mathscr P}}

\newcommand{\eN}{{\mathscr N}}
\newcommand{\co}{{\mathscr O}}
\newcommand{\eF}{{\mathscr F}}
\newcommand{\eG}{{\mathscr G}}

\newcommand{\Spec}{{\mathrm {Spec}}}
\newcommand{\Spf}{{\mathrm {Spf}}}

\newcommand{\ares}[1]{\ush{\boldsymbol {\mathrm{res}}}_{\<\<\<\<\<\<\<\<\<{}_#1}}
\newcommand{\sres}[1]{{\boldsymbol {\mathrm{res}}}^{\<{}_S}_{\<{}_#1}}
\newcommand{\res}[1]{{\boldsymbol {\mathrm{res}}}_{#1}}

\newcommand{\Res}[4]{{\mathrm{Res}}_{#1}\begin{bmatrix}{#2}\\ {#3}, \dots, {#4}\end{bmatrix}}

\newcommand{\A}{{\mathcal A}}

\newcommand{\De}{{\Delta}}

\newcommand{\F}{{\mathcal F}}

\newcommand{\Hr}{{\mathrm H}}

\newcommand{\Rr}{{\mathrm R}}

\newcommand{\eK}{{\mathscr K}}

\newcommand{\eL}{{\mathscr L}}
\newcommand{\eQ}{{\mathscr Q}}
\newcommand{\D}{{\mathbf D}}
\newcommand{\K}{{\mathbf K}}

\newcommand{\bbP}{{\mathbb P}}

\newcommand{\bbG}{{\mathbb G}}

\newcommand{\vc}{{\vec{\mathrm{c}}}}

\newcommand{\Dqc}{\D_{\mkern-1.5mu\mathrm {qc}}}
\newcommand{\wDqc}{ \widetilde
          {\vbox to6.5pt{\vss\hbox{$\mathbf D$}}}
    _{\mkern-1.5mu\mathrm {qc}} }

\newcommand{\wDqcp}{\wDqc^{\lower.5ex\hbox{$\scriptstyle+$}}}

\newcommand{\Dvc}{\D_{\<\vc}}
\newcommand{\Dqct}{\D_{\mkern-1.5mu\mathrm{qct}}}

\newcommand{\Dc}{\D_{\mkern-1.5mu\mathrm c}}

\newcommand{\qc}{{\mathrm{qc}}}

\newcommand{\R}{{\mathbf R}}
\newcommand{\Rp}[1]{{{\Rr}'}_{\<\<\<{{#1}}}}
\newcommand{\Rfs}{{\mathbf R f_{\!*}}}
\newcommand{\bL}{{\mathbf L}}
\newcommand{\Hom}{{\mathrm {Hom}}}

\newcommand{\Homb}{{\mathrm {Hom}}^{\bullet}}
\newcommand{\Ac}{\A_{\mathrm c}}
\newcommand{\Aqc}{\A_{\qc}}
\newcommand{\Avc}{\A_{\vec {\mathrm c}}}
\newcommand{\Aqct}{\A_{\mathrm {qct}}\<}

\newcommand{\fs}{f^!}

\newcommand{\ush}[1]{{#1^{\textup{\texttt\#}}}}
\newcommand{\trc}[1]{{\mathrm{Trc}}_{{#1}}}
\newcommand{\Tr}[1]{{\mathrm {Tr}}_{{#1}}}

\newcommand{\ttr}[1]{{\boldsymbol{\ush{\tau}_{\!\!\!\!{}_{#1}}}}}
\newcommand{\vttr}[1]{{\boldsymbol{\tau_{{}_{#1}}}}}
\newcommand{\vin}[1]{{\mathrm{tr}_{{#1}}}}
\newcommand{\tin}[1]{{\ush{\mathrm{tr}}_{{{#1}}}}}
\newcommand{\omgs}[1]{{\ush{\omega_{\<\<{#1}}}}}
\newcommand{\oreg}[1]{{\omega^{\mathrm{reg}}_{#1}}}
\newcommand{\rin}[1]{{\boldsymbol{\int_{#1}^{\mathrm{reg}}}}}
\newcommand{\wnor}[1]{{{\eN}}^r_{{#1}}}
\newcommand{\wI}[2]{({\wedge_{#1}^r{#2}/{#2}^2})^*}
\newcommand{\wdd}[2]{{\mathrm d}{#1}\wedge\dots\wedge{\mathrm d}{#2}}
\newcommand{\Ext}{\operatorname{\eE{\mathit{xt}}}}

\newcommand{\sHom}{\eH{om}}

\newcommand{\iGp}[1]{{\varGamma_{\<\!#1}'}}
\newcommand{\iG}[1]{{\varGamma_{\<\!#1}^{\phantom\prime}}}

\newcommand{\wid}[1]{\widehat{#1}}
\newcommand{\wit}[1]{\widetilde{#1}}
\newcommand{\set}{\!:=}
\newcommand{\lra}{\longrightarrow}
\newcommand{\iso}%
{{\mkern8mu\longrightarrow \mkern-25.5mu{}^\sim\mkern17mu}}
\newcommand{\osi}%
{{\mkern8mu\longleftarrow \mkern-24.5mu{}^\sim\mkern16mu}}
\newcommand{\Otimes}{\underset
   {\vbox to 0pt {\vskip-1ex\hbox{$\scriptscriptstyle=$}\vss}}
     {\otimes}\vadjust{\kern.4pt}}

\newcommand{\BL}{{\boldsymbol\Lambda}}

\newcommand{\smcirc}%
   {{\raise.15ex\hbox to.7em{$\hss \scriptstyle\circ\hss$}}}

\title[Transitivity II]%
{Grothendieck Duality and Transitivity II: Traces and residues via
Verdier's isomorphism}

\author[S.\,Nayak]{Suresh Nayak}%
{\smash{}}
\address{Indian Statistical Institute \\
8th Mile, Mysore Road, Bangalore \\
Karnataka-560059, INDIA}
\email{snayak@isibang.ac.in}

\author[P.\,Sastry]{Pramathanath Sastry}%
{\smash{}}
\address{Chennai Mathematical Institute\\
                 Sipcot IT park, Siruseri\\
                 Kanchipuram Dist TN, 603103, INDIA}
\email{pramath@cmi.ac.in}

\makeindex

\date{\today}

\begin{document}

\begin{abstract}{For a smooth map between noetherian schemes, Verdier relates 
the top relative differentials of the map with the twisted inverse image functor
``upper shriek" \cite{verdier}. We show that the associated traces for smooth
proper maps can be rendered concrete by showing that the resulting theory of residues
satisfy the residue formulas (R1)--(R10) in Hartshorne's \emph{Residues and Duality} \cite{RD}.
We show that the resulting abstract transitivity map relating the twisted image functors
for the composite of two smooth maps satisfies an explicit formula involving differential
forms.
We also give explicit formulas for traces of differential forms for finite flat maps
 (arising from Verdier's isomorphism) between
schemes which are smooth over a common base, and use this  to relate Verdier's isomorphism 
to Kunz and Waldi's regular differentials. These results also give concrete
realisations of traces and residues for Lipman's fundamental class map via the results
of Lipman and Neeman \cite{fund-loc} relating the fundamental class to Verdier's isomorphism.}
\end{abstract}

\maketitle



All schemes formal or ordinary are assumed to be noetherian. The category of $\co_\X$-modules
for a formal scheme $\X$ is denoted $\A(\X)$ and its derived category $\D(\X)$ as in
\cite{fub-abs}. 
In general we use the notations of \textit{ibid}. Thus the torsion functor $\iGp{\X}$ on $\co_\X$-
modules is defined by the formula
\[\iGp{\X} \set \dirlm{n}\sHom_{\co_\X}(\co_\X/\I^n, -)\] 
where $\I$ is any ideal of definition of the formal scheme $\X$. A {\emph{torsion module}}
$\eF$ is an object in $\A(\X)$ such that $\iGp{\X}\eF=\eF$.
The reader is advised to look at \cite[\S\,2]{fub-abs} for further definitions
and notations, especially the definitions of $\Ac(\X)$, $\Avc(\X)$, $\Aqc(\X)$, $\Aqct(\X)$,
and the triangulated full subcategories of $\D(X)$, $\Dc(\X)$, $\Dvc(\X)$, $\Dqc(\X)$, $\Dqct(\X)$,
and their various bounded versions (e.g., $\Dc^+(\X)$ \dots).

\section{\bf Introduction}

The principal aim of this paper is to describe
explicitly the residues---and the trace, when the map in question is proper---associated with
Verdier's isomorphism 
\[f^!\cong f^*(\boldsymbol{-})\otimes\Omega^n_{X/Y}[n]\] 
(see \cite[p.\,397, Thm.\,3]{verdier}) for a smooth map 
$f\colon X\to Y$ of relative dimension $n$. 
The foundations of Grothendieck duality (GD) that we use are the ones initiated by 
Deligne in \cite{del}. 

For this introduction, unless otherwise stated, schemes are ordinary noetherian schemes. In the main
body of the paper, we use formal schemes as way around compactifications of separated
finite type maps, so that complications involving
 compatibilities between different compactifications 
do not need to be addressed. 
\pagebreak
The principal input, when we use formal schemes, is \cite{fub-abs}
which should be regarded as a companion paper, written mainly with this manuscript in mind.

There are two quite different ways that GD is constructed. 
The foundations for GD used in \cite{RD} and \cite{conrad} are based on residual complexes.
In this approach, the functor $f^!$ (for a suitable finite type map $f$), 
as well as its attendant trace $\Tr{f}\colon \Rfs f^!\to {\bf 1}$ (when $f$ is proper) have
a certain concreteness built into their construction. One then has to work out a large array
of compatibilites between the various concrete representations of $f^!$ and $\Tr{f}$ and there
often are different concrete representations of these for the same map, e.g., a finite map which
also factors as a closed immersion followed by a smooth map. 
In a different direction, in his appendix to \cite{RD}, Deligne initiated an approach to GD which is
conceptually attractive \cite{del}. From this point of view, $f^!$ for a proper map $f$ is the right adjoint
to $\Rfs$, and exists for very general category theoretic reasons. These foundations have
been worked on, extended, and new techniques introduced over the years by Lipman, Neeman, 
and their collaborators. Residual complexes and dualizing complexes are not needed to build GD in 
this approach. We mention \cite{del}, \cite{del-sga}, \cite{del-sga'}, \cite{verdier} for literature
on this approach before the 1980s, and  recent work found in \cite{bous}, \cite{nee-bc},
which do much to extend (via a conceptually different approach to finding right adjoints)
the work initiated by Deligne and Verdier to more general situations, often bypassing the old
annoying hypotheses on boundedness for the existence of $f^!$ or for its base change. The stable
version of these can be found in Lipman's elegant and carefully written book \cite{notes}.
We also mention Neeman's recent manuscript \cite{nee-simp} which gives a coherent account
of the difficulties and the recent simplifications of many matters. Since we rely on formal schemes
as a way to our results on maps between ordinary schemes, we have been influenced
enormously by the work of Lipman, Alonso Tarr\'{\i}o, and Jerem\'{\i}as L\'opez, in
\cite{loch}, \cite{dfs}, and \cite{gm}. We rely on our results on abstract transitivity on formal
schemes and related matter in \cite{fub-abs}. 

Given the highly abstract methods of construction, and definitions based on
universal properties, the question arises: 
\begin{quotation}
\emph{To what extent can we render concrete
realisations of the various
constructions occurring in this version of GD?}
\end{quotation}

The issue of concrete representations of $f^!$ (for our preferred version of GD)
was addressed partially, soon after 
\cite{del} appeared, by Verdier when $f$ is smooth \cite{verdier}. The answer is
$f^!\cong f^*(\boldsymbol{-})\otimes_{\co_X}\Omega^n_{X/Y}[n]$ where $n$ is the relative dimension
of $f\colon X\to Y$. This isomorphism in turn depends on the concrete representation 
$i^!\cong \bL i^*(\boldsymbol{-})\otimes \wedge^d_{\co_U}\eN[-d]$ (via the 
fundamental local isomorphism) for a regular immersion $i\colon U\hookrightarrow V$
of codimension $d$, with $\eN$ the normal bundle
of $U$ in $V$. Verdier's answer for smooth maps is only a partial answer because
the associated trace map (when $f$ is proper)
\[\vin{f}\colon \Rr^nf_*(\Omega^n_{X/Y})\lra \co_Y,\]
denoted $\int_f$ in \cite{verdier},\footnote{The map $\vin{f}$ is 
$\Hr^0(\boldsymbol{-})$ applied to the composite $\Rfs\Omega^n_{X/Y}[n] \iso \Rfs f^!\co_Y
\xrightarrow{\Tr{f}(\co_Y)} \co_Y$, where the first arrow is Verdier's isomorphism.} 
is seemingly intractable via this approach. In fact $\vin{f}$ has not been worked out in the literature
 even when $A$ is
the spectrum of a field $k$ (e.g., $k=\mathbb{C}$), $X=\bbP^n_k$, and $f$ the structure
map. However, from the abstract properties of $f^!$ and the fact that $f^!\co_Y$ is concentrated
in degree $-n$ (for example by Verdier's isomorphism), the pair $(\Omega^n_{X/Y}, \vin{f})$
is easily seen to represent the functor $\Hom_Y(\Rr^nf_*(\boldsymbol{-}), \co_Y)$ on quasi-coherent
sheaves on $X$ when $f$ is proper (the only situation where $\vin{f}$ is defined), and in fact 
$\vin{f}$ and the composite
$\Rfs\Omega^n_{X/Y}[n] \iso \Rfs f^!\co_Y \xrightarrow{\Tr{f}(\co_Y)} \co_Y$ determine
each other.

When $Z$ is a closed subscheme of $X$, finite over $Y$, defined locally
by an $\co_X$-sequence, and $\Ext^i_f(\co_Z, \boldsymbol{-})$ the $i^{\mathrm{th}}$ right
derived functor of $f_*\sHom_X(\co_Z, \boldsymbol{-})$, Verdier asserts (see top of p.\,400 of
\cite{verdier}) that the composite 
\stepcounter{thespecial}
\begin{equation*}
\label{intro:ext}\tag{\thethespecial}
\Ext^n_f(\co_Z, \Omega^n_{X/Y}) \lra \Rr^nf_*(\Omega^n_{X/Y}) \xrightarrow{\vin{f}}\co_Y
\end{equation*}
is governed by the residue symbol of \cite[Chap.\,III, \S\,9]{RD}.\footnote{In \cite{RD}, the proofs
of the assertions about the residue symbol are not given. They are provided later by Conrad 
\cite{conrad}, the construction and definition of various traces being those developed in \cite{RD}. 
They do not apply to our situation since we use a different foundation for GD.}  It is certainly
true that if this is so, following (essentially) the argument given in \cite[bottom of p.\,399]{verdier},
the trace map
$\vin{f}$ can be realised in an explicit way. However, the proof that \eqref{intro:ext} (denoted
${\mathrm{Res}}_Z$ in \cite{verdier}) is governed by the residue symbol is not there in the literature.
In the over 50 years that have passed since Verdier's assertion, it has been recognised
by experts that this is in fact a non-trivial problem (see our quote of Conrad below).
One difficulty is the assertion (R4) in \cite[p.\,400]{verdier}, namely that
\eqref{intro:ext} commutes with arbitrary base change. This needs, at the very
least, for one to show that the isomorphism $f^!\co_Y\cong \Omega^n_{X/Y}[n]$ of Verdier commutes
with arbitrary base change, in a sense we will make more precise in a moment. This
compatibility with base change was only established in 2004 by the second auhor \cite{cm}. In
slightly greater detail, here is what the compatibility entails.
Suppose $u\colon Y'\to Y$ is a map and $g\colon X\times_YY'\to Y'$ and $v\colon X\times_YY'\to X$
are the two projections. To show the compatibility of Verdier's isomorphism 
 $f^!\co_Y\cong \Omega^n_{X/Y}[n]$ with
arbitrary base change, first one needs
to show that there is an isomorphism 
\[\theta_u^f\colon v^*f^!\co_Y\iso g^!\co_{Y'}\]  
for our smooth $f$ (even when it is not proper). This
is a delicate point, especially if one demands that in the proper case $\Tr{f}(\co_Y)$ should
be compatible with this base change isomorphism (remember, $\Tr{f}$ in Deligne's approach,
is defined as a co-adjoint unit and is not explicit), and that the isomorphism
is also compatible with open immersions into $X$. After this is established, one has to check
that this base change isomorphism $\theta_u^f$ when grafted on to Verdier's isomorphisms for
$f$ and for $g$ give the canonical isomorphism of differential forms. It is easier to carry out
the first part in the slightly more general situation of $f$ being Cohen-Macaulay, and this
is one of the main results of \cite{cm}. In \cite{conrad},
the base change isomorphism $\theta_u^f$ is proven using the foundations of GD in \cite{RD}.
However, the isomorphism between $f^!\co_Y$ and $\Omega^n_{X/Y}[n]$ in \cite{RD}
and \cite{conrad} is by fiat, and it is not clear that it is the same as Verdier's isomorphism.
In other words, it is not clear that the trace $\Rr^nf_*(\Omega^n_{X/Y})\to\co_Y$
built using the foundations of GD in \cite{RD} is the same as the one that arises
when using the foundations initiated in \cite{del}.
In fact, we are back to the frustrating detail that 
 we do not know the $\vin{f}$ explicitly when we work with the foundations 
initiated in \cite{del}. 
Even with the compatibility of Verdier's isomorphism with arbitrary base
change in hand, showing that \eqref{intro:ext} is governed by the residue symbol of
\cite[Chap.\,III, \S\,9]{RD} is not trivial. In fact it takes all this paper.
We can do no better than quote 
Conrad from his introduction to his book \cite{conrad} (using however our labelling of the citations 
given there):

\begin{quotation} ``\dots The methods in \cite{verdier} take place in derived
categories with ``bounded below" conditions. This leads to technical problems
for a base change such as $p\colon \Spec{(A/\fm)}\hookrightarrow \Spec{(A)}$
with $(A,\,\fm)$ a non-regular local ring, in which case the right exact $p^*$ does not
have finite homological dimension (so $\bL p^*$ does not make sense as
a functor between ``bounded below" derived categories). 
Moreover, Deligne's construction of the trace map in \cite[Appendix]{RD},
upon which \cite{verdier} is based, is so abstract that it is a non-trivial task to relate Deligne's
construction to the sheaf $\Rr^nf_*\(\Omega^n_{X/Y})$. However, a direct relation between the
duality theorem and differential forms is essential for many important calculations (e.g., 
\cite[\S6, \S14(p.121)]{maz})."
\end{quotation}

In other words the task of finding a concrete expression for $\vin{f}$
is not simple. In this paper we take up this task, and believe we give a satisfactory answer
to the problem. Briefly, any theory of traces comes with an associated theory of residues,
and we show that residues associated  with $\vin{f}$ satisfy the formulas \cite[III, \S\,9]{RD}, which are
stated without proof in \textit{loc.cit.}\footnote{though all the formulas stated in \cite{RD} (labelled (R1)--(R10) there) have been proved with great care by Conrad in \cite[Appendix A]{conrad}.} 

The prime object of study in this paper is Verdier's isomorphism \cite[p.\,397, Thm.\,3]{verdier}
\[ \Omega^n_{X/Y}[n] \iso f^!\co_Y\]
for a smooth separated morphism $f\colon X\to Y$ of ordinary schemes of relative dimension $n$.
Strictly speaking, the isomorphism in \textit{loc.cit.}~is from $f^!\co_Y$ to $\Omega^n_{X/Y}[n]$,
and thus, we are talking about the inverse of the map in \textit{loc.cit}. In view of recent results of Lipman and Neeman, this is the {\emph{fundamental class map $c_f$ associated with $f$ 
\cite[p.\,152, (4.4.1)]{fund-loc}}, 
but we use the description given in \cite{verdier} and hence call it the Verdier isomorphism. In
\cite{mexico}, Lipman outlines a programme for a global residue theorem via the fundamental
class map (see [\textit{ibid.}, \S\,5.5 and \S\,5.6]). This paper is intimately related to that
programme via the just mentioned results of Lipman and Neeman. However, we do not
use the results on the fundamental class map of [\textit{ibid.}]. Since the isomorphism
we use (between $\Omega^n_{X/Y}[n]$ and $f^!\co_Y$) is that described by Verdier, we call
it the Verdier isomorphism rather than the fundamental class.

We also recommend Beauville's expository paper \cite{beauville} for an overview (without proofs)
of residues, especially for the concrete expressions for them. Our attention was drawn to it recently
by Joe Lipman. 

We now give a more more detailed description of the contents of the paper. We are
concerned, mainly with three (intertwined) aspects:
\begin{enumerate}
\item[{\bf 1.}] Understanding the abstract traces
\[\Tr{f}(\co_Y)\colon \Rfs f^!\co_Y \to \co_Y\]
and
\[\Tr{f, Z}(\co_Y) \colon \R_Z f_* f^!\co_Y \to \co_Y\]
in concrete terms (using differential forms via Verdier's isomorphism)
when $f\colon X\to Y$ is smooth and separated. The first map is meaningful when $f$ is proper, as
the co-adjont unit for the adjoint pair $(\Rfs, f^!)$ 
\cite[(1.1.2)]{fub-abs}. The second
is meaningful when $Z$ is a closed subscheme of $X$ proper over $Y$ 
\cite[(3.3.1)]{fub-abs}. If $Z=X$, $\Tr{f,Z}=\Tr{f}$. 
In fact, we will concentrate on the case when $Z$ is finite over $Y$,
in which case we are talking about abstract residues. The aim to is realise these abstractions
concretely when we substitute $\Omega^n_{X/Y}[n]$ for $f^!\co_Y$ via
Verdier's isomorphism ($n$ being the relative dimension of $f$). Understanding $\Tr{f,Z}$ for
such $Z$ is tantamount to understanding $\Tr{f}$ for $f$ proper via the so-called Residue Theorem.

\item[{\bf 2.}]\label{item:foo} Making concrete the abstract transitivity map
\[\chi_{[g,f]}\colon \bL f^*g^!\co_Z\overset{\bL}{\otimes}_{\co_X}f^!\co_Y \lra (gf)^!\co_Z\]
of \cite[\S\,4.9]{notes} and \cite[(7.2.16)]{fub-abs} 
concrete in terms of differential forms (again using Verdier's isomorphism)
when $f\colon X\to Y$ and $g\colon Y\to Z$ are separated finite-type maps in certain situations.
Our main interest is in the following two situations:
\begin{enumerate}
\item[(i)] The maps $f$ and $g$ are smooth, say of relative dimensions $m$ and $n$ respectively,
and we use Verdier's isomorphisms to identify $g^!\co_Z$, $f^!\co_Y$, and $(gf)^!\co_Z$ with
$\Omega^n_{Y/Z}[n]$, $\Omega^m_{X/Y}[m]$, and 
$\Omega^{m+n}_{X/Z}[m+n]$ respectively. This is closely related
to the results in \cite{jag}.
\item[(ii)] The map $f$ is a closed immersion say of codimension $d$, 
and the maps $g$ and $gf$ are smooth, say of relative dimensions $n+d$ and $n$ respectively.
\end{enumerate}
In fact these two cases are essentially enough to develop a theory of residues which give the
formulas (R1) to (R10) in \cite[Chap.\,III, \S\,9]{RD}.

\item[{\bf 3.}] Finding a concrete expression for the abstract trace map
\[h_*f^!\co_Z\cong h_*h^!g^!\co_Z \xrightarrow{\Tr{h}}g^!\co_Z\]
where $f\colon X\to Z$ and $g\colon Y\to Z$ are smooth separated
maps and $h\colon X\to Y$ is a finite flat map. This
concrete expression is in terms of differential forms (via our now familiar way of identifying
$f^!\co_Z$ and $g^!\co_Z$ with differential forms). In fact we show that
it is the trace of Lipman and Kunz,
defined in \cite[\S\,16]{kd}. One consequence is that if $f\colon X\to Y$ is an equidimensional
map of relative dimension $n$ such that $X$ and $Y$ are excellent with
no embedded points and the smooth locus of $f$ is dense in $X$, then $\Hr^{-n}(f^!\co_Y)$
can be identified, via the Verdier isomorphism on the smooth locus of $f$, with a coherent 
subsheaf of the sheaf of meromorphic differentials $\wedge^n_{k(X)}\Omega^1_{k(X)/k(Y)}$,
namely the sheaf of regular differentials of Kunz and Waldi \cite[\S\,3, \S\,4]{kw}. 
One therefore recovers the main results in 
\cite{hk1}, \cite{hk2}, \cite{ajm}, and \cite{jag} via our approach.
\end{enumerate}

We elaborate on these points in the rest of this introduction.

\subsection{The twisted image pseudofunctor $\boldsymbol{-^!}$} GD is
concerned with constructing a variance theory, i.e., a pseudofunctor, ``upper shriek", which we denote 
$\boldsymbol{-^!}$, on a suitable subcategory of schemes and finite type maps.\footnote{One can have essentially finite type maps.} For a
fixed scheme $Z$, $Z^!$ is a suitable full subcategory of $\D(Z)$ containing $\Dc(Z)$. We will
say more about these subcategories later. For now we wish to paint with broad strokes. 
Whichever way one approaches the foundations of GD, the resulting pseudofunctor $\boldsymbol{-^!}$ 
should be {\em local} (more on that in a moment), stable under, at least, flat base change,
and such that when $f$ is {\em proper}, $f^!$ is right adjoint to $\Rfs$. 
By {\em local}, this is what we mean:
If $U\rightarrow Y$ is an open $Y$-subscheme
of $g\colon V\to Y$ as well as of $h\colon W\to Y$ ($g$, $h$ of finite type), then $g^!\vert_U$
and $h^!\vert_U$ are canonically isomorphic --- canonical enough that if we have a third
finite type $Y$-scheme $f\colon X\to Y$ which contains $U$ as an open $Y$-subscheme, then
the isomorphisms between $f^!\vert_U$, $g^!\vert_U$, and $h^!\vert_U$ are compatible.
All of this (and much more) can be found in \cite{notes} for the theory of $\boldsymbol{-^!}$
initiated in \cite{del}. For schemes with finite Krull dimension, the local nature of upper-shriek
was proved by Deligne in \cite{del}, and using his flat base change result, by Verdier in 
\cite{verdier}.

Additonally, one wants a theory which
specializes to the familiar Serre duality for smooth complete varieties, with the top differential forms
playing a critical dualizing role. For a slightly more general situation, this means that from the
theory of upper shriek one should recover the duality isomorphisms 
\[\Ext^i_f(\V,\,\Omega^d_{X/Y})\cong \sHom_Y(\Rr^{d-i}f_*\V,\,\co_Y) \qquad (0\le i \le d) \]
 when $f\colon X\to Y$ is {\em smooth and proper} of
relative dimension $d$, $\V$ is a finite rank vector bundle on $X$. This amounts to showing that 
$f^!\co_Y\cong \Omega^d_{X/Y}[d]$ for such a smooth map $f$.

\subsection{Traces and residues} 
Given a theory of upper shriek there is an associated theory of traces and 
residues related to it in a very close manner. Briefly, for $f\colon X\to Y$ proper, since $f^!$ is a right adjoint to $\Rfs$, there is a co-adjoint unit
\[\Tr{f}\colon \Rfs f^!\to {\bf 1}\] 
called the {\em trace map},
namely the image of the the identity transformation $f^!\to f^!$ under $\Hom_{X^!}(f^!, f^!)\iso
\Hom_{Y^!}(\Rfs f^!,\,{\bf 1})$ where ${\bf 1}$ = the identity functor on $Y^!$. More generally,
if $f$ is separated of finite type and $Z$ is closed subscheme of $X$ which is {\em proper} over $Y$, 
then one has a map (the {\emph{trace along $Z$}})
\[\Tr{f,Z}\colon \R_Zf_*f^! \to {\bf 1}\] 
which for this introduction can be defined as follows: If $F:P\to Y$ is a compactification\footnote{i.e., $F$ is proper and contains $X$ as an open $Y$-subscheme---such an $F$ can always be arranged
\cite{nagata}.} of $f$, 
and $i\colon X\to P$ the open $Y$-immersion associated to this compactification, then using the
isomorphism $f^!\iso i^*F^!$, we have the composite
\[\R_Zf_*f^! \iso \R_Zf_*i^*(F^!) \iso \R_{i(Z)} F_*F^! \lra \R F_*F^! \xrightarrow{\Tr{F}} {\bf 1}.\]
The above composite is independent of the compactification data $(i, F)$, and we define $\Tr{f,Z}$
to be this composite. 
We give an equivalent but more useful definition in  
\cite[\S\S\,3.3]{fub-abs}.
The map $\Tr{f,Z}$ depends only on the formal completion of $X$ along $Z$ and not on the exact
subscheme structure of $Z$. Note that when $Z=X$ (this implies $f$ is proper), then $\Tr{f,Z}=\Tr{f}$.

If $f$ is smooth of relative dimension $n$, $Z$ as above is cut out by 
${\bf t}=(t_1, \dots, t_n)\in\Gamma(X,\,\co_X)$ with $Z\to Y$ finite, 
and $Z_m$ is the thickening of $Z$ defined by
by ${\bf t}^m=(t_1^m, \dots, t_n^m)$, and say $Y=\Spec{\,A}$, then \eqref{intro:ext} gives us
maps (one for each $m$)
\[\mathrm{Ext}^n_X(\co_{Z_m},\, \Omega^n_{X/Y}) \lra A.\]
As Verdier argued in \cite[bottom of p.399]{verdier}, by passing to the completion of a localisation of $A$ (via flat base
change), and making \'etale base changes, to know the above map (for any $m$) is to know
$\Tr{f}$ when $f$ is proper. If one passes to the direct limit as $m\to\infty$, then we get
a map
\[
\Hr^n_Z(X,\,\Omega^n_{X/Y}) \lra A.
\]
The above map is easily seen to be $\Hr^0(\boldsymbol{-})$ applied to the composite
\stepcounter{thm}
\begin{equation*}\label{intro:TrZ}\tag{\thethm}
\R\Gamma_Z (X, \Omega^n_{X/Y}[n]) \iso \R\Gamma_Z (X,\, f^!\co_Y) \xrightarrow{\Tr{f, Z}} A.
\end{equation*}
This map, which we denote $\res{Z}$, also determines $\Tr{f}$ if $f$ is proper. We prefer to
work with cohomology with supports (rather than with
$\mathrm{Ext}^n_X(\co_Z,\,\Omega^n_{X/Y})$), following the general
philosophy underlying Lipman's body of work, especially \cite{ast117}. Berthelot in \cite{berthelot}
also makes the connection between the map on $\mathrm{Ext}^n_X(\co_Z, \Omega^n_{X/Y})$
and the map on cohomology with supports. However Berthelot
 uses the foundations of GD based on residual complexes.

The relationship between upper shriek and the associated traces is intimate. 
\emph{To assert that one has a concrete understanding of upper shriek 
in a particular situation is tantamount to asserting that
one understands $\Tr{f,Z}$ for a certain class of closed subschemes $Z$ which are
proper over the base scheme $Y$}. For example if $f$ is a 
Cohen-Macaulay map (i.e., a flat finite type map with Cohen-Macaulay fibres), then to ``know" 
$\Tr{f,Z}(\co_Y)$ for $Z$ which are finite and flat over $Y$ is to ``know" duality for $f$.

Returning to the case we are discussing ($f$ smooth of relative dimension $n$), suppose
$Z\hookrightarrow X$ is a closed immersion cut out by a sequence of global
sections ${\bf t}= (t_1, \dots, t_n)$ of $\co_X$, and
$Y=\Spec{\,A}$. Assume  $Z\to Y$ is an isomorphism and contained
in an affine open subscheme $U=\Spec{\,B}$ of $X$, something that can be achieved by
shrinking $Y$, since $Z\to Y$ is an isomorphism. Let
\[\res{Z}\colon \Hr^n_Z(X,\,\Omega^n_{X/Y})=\Hr^n_Z(U,\,\Omega^n_{U/Y})  
\xrightarrow{\phantom{XXX}} A\]
be $\Hr^0(\eqref{intro:TrZ})$.
It is well known that elements of $\Hr^n_{{\bf t}B}(\Omega^n_{B/A})$ are finite $A$-linear
combinations of elements of the form 
$\Bigr[\begin{smallmatrix} dt_1\wedge\dots\wedge t_n\\ t_1^{\beta_1},\dots, t_n^{\beta_n}
\end{smallmatrix}\Bigl]$
with $\beta_i$ positive integers. Ideally one would like 
\stepcounter{thm}
\begin{equation*}\label{intro:std-res}\tag{\thethm}
\res{Z}\begin{bmatrix} dt_1\wedge\dots\wedge t_n\\ t_1^{\beta_1},\dots, t_n^{\beta_n}
\end{bmatrix}=
\begin{cases}
1 & \text{when $\beta_i=1$ for all $i=1,\dots, n$}\\
0 & \text{otherwise}.
\end{cases}
\end{equation*}
The exact answer depends on the isomorphism $f^!\co_Y\iso \Omega^n_{X/Y}[n]$ chosen. 
This is at the 
heart of this paper, since our choice is the isomorphism Verdier gives in \cite[p.\,397, Thm.\,3]{verdier}.
In fact we show that Verdier's isomorphism does give the above formula in the case being
considered, i.e., when $Z$ is a section of $f$. This is the critical case, and we deduce
other residue formulas from this one by either making \'etale base changes, or base changing
$f$ by itself and using the diagonal section $X\hookrightarrow X\times_YX$ of the first projection
(which is to be thought of as the base change of $f$).

We could obtain the above explicit description of $\res{Z}$ when $Z$ is a section of $f$ because
of the results in \cite{cm}. The main results there state that if $f\colon X\to Y$ is Cohen-Macaulay
of relative dimension $d$, then for any base change $u\colon Y'\to Y$, there is a natural
isomorphism $\theta_u^f\colon v^*\omgs{f}\iso \omgs{g}$, where $v\colon X\times_YY'\to X$
and $g\colon X\times_YY'\to Y'$ are the respective projections. When $f$ is proper, this 
isomorphism is compatible with traces. If $f$ is smooth (proper or not), then this isomorphism
when transferred to $v^*\Omega^n_{X/Y}[n]$ and $\Omega^n_{X\times_YY'/Y'}[n]$ is the 
identity map under the standard identification of differential forms. These are very similar
to the main results in \cite{conrad}. The difference is that in \cite{cm} the foundations of GD
are based on the one initiated by Deligne in \cite{del}, whereas in \cite{conrad} it is the based
on residual complexes. In \cite{conrad} the identification of differential forms is built into the definition
of the base change isomorphism between $v^*\omgs{f}$ and $\omgs{g}$, since the strategy
is to embedd $X$  into schemes smooth over $Y$. 
The challenge in \cite{conrad} is to show that the result
is compatible with traces when $f$ is proper.

In our approach to finding explicit formulas for $\res{Z}$, 
the role played by $\theta_u^f$, \emph{when $u$ is non-flat}, is crucial. Roughly speaking,
Verdier's isomorphism can be regarded as the residue formula 
$\res{\De}\Bigr[\begin{smallmatrix}\wdd{s_1}{s_n}\\ s_1, \dots, s_n\end{smallmatrix}\Bigr]=1$
for the diagonal section $\De$ in $X\times_YX$ where the diagonal is cut out by ${\bf s}$ in
$X\times_YX$. If $Z\hookrightarrow X$ is a section of $f$ , 
cut out by $t_1, \dots, t_n\in\Gamma(X,\,\co_X)$
then pulling back the diagonal via the base change
$Z\to X$, we get 
\stepcounter{thm}
\begin{equation*}\label{intro:std-res=1}\tag{\thethm}
\res{Z}\begin{bmatrix}\wdd{t_1}{t_n}\\ t_1, \dots, t_n\end{bmatrix}=1.
\end{equation*}
We can do this because Verdier's isomorphism is compatible with arbitrary base change --
the result in \cite[p.740, Thm.\,2.3.5\,(b)]{cm}  that we alluded to above. The
formula \eqref{intro:std-res=1} says that \eqref{intro:std-res} is true when all the $\beta_i$ are $1$.
If $X=\bbP^n_Y$, $f$ the standard projection $\bbP^n_Y\to Y$, and $Z=\cap_{i=1}^n\{T_i\neq 0\}$, 
where $T_i$, $i=0, \dots, n$ are homogeneous co-ordinates on $\bbP^n_Y$ (and $t_i=T_i/T_0$,
for $i=1, \dots, n$), then one can show easily that  \eqref{intro:std-res=1} implies \eqref{intro:std-res}.
The crucial ingredient needed is the simple and 
elegant computation of Lipman in \cite[pp.79--80,\,Lemma\,(8.6)]{ast117}. The proof is
essentially carried out in the proof of \Pref{prop:TrS=int} (ii). Since $\res{Z}$ depends only
on the formal completion of $X$ along $Z$, therefore if $Z$ and ${\bf t}$ satisfy the hypotheses
given when stating \eqref{intro:std-res}, then formula \eqref{intro:std-res} holds. 
This is the first, and a very important step in our proofs in \Sref{s:res-sym}
of the residue formulas (R1)--(R10) of \cite[Chap. III, \S\,9]{RD}. 

If the closed subscheme $Z$ of $X$ cut out by ${\bf t}=(t_1, \dots, t_n)$ is finite over $Y$ 
(and hence necessarily flat over $Y$) of constant rank 
(not necessarily an isomorphism), and of constant rank then it turns out that the right side of 
\eqref{intro:std-res=1} needs to be replaced by $\mathrm{rank}{(Z/Y)}$.

\begin{rem}\label{intro:dirac} One way to think about formula \eqref{intro:std-res=1} is to
regard 
\[\varphi \mapsto \res{Z}\begin{bmatrix}\varphi\>\wdd{t_1}{t_n}\\ t_1,\dots, t_n\end{bmatrix},\]
(for $\varphi$ a section of $\co_X$ in an open, or even formal, neighbourhood of $Z$)
as the Dirac distribution along $Z$. Indeed, (with $Y=\Spec{(A)}$), since $Z$ is a section of $f$, the completion of  $X$ along $Z$ is the power series ring $A[|t_1, \dots, t_n|]$, and
according to \eqref{intro:std-res=1}, the right side is $\varphi(0,\dots, 0)$, after developing $\varphi$
as a power-series in ${\bf t}$.
If $A=\mathbb{C}$, and
field of complex numbers (so that $Z=\{p\}$, a point), 
this can be interpreted as the fact that the Dolbeault representative
of the Cauchy kernel $\wdd{t_1}{t_n}/t_1\dots t_n$ at the point $p$ is the Dirac distribution at $p$
(see \cite{tong}). 
\end{rem}

\subsection{Transitivity} Finding concrete expressions
(when we have  two finite type separable maps $f\colon X\to Y$ and $g\colon Y\to Z$)
 for the abstract transitivity map
\[\chi_{[g,f]}\colon \bL f^*g^!\co_Z\overset{\bL}{\otimes}_{\co_X}f^!\co_Y \lra (gf)^!\co_Z\]
of \cite[\S\,4.9]{notes} and \cite[(7.2.16)]{fub-abs} 
is perhaps the most important technical task undertaken in this paper. To establish this,
we rely heavily on the abstract transitivity results on formal schemes in \cite{fub-abs}. 
As we pointed out earlier (see item ({\bf 2.})~on p.\pageref{item:foo}) 
there are two key situations where concrete manifestations of $\chi_{[-,-]}$ are important. 
The first situation of importance is when we have two smooth separated maps , 
$f\colon X\to Y$ and $g\colon Y\to Z$, say of relative dimensions $m$ and $n$ respectively.
If we use Verdier's isomorphisms to identify $g^!\co_Z$, $f^!\co_Y$, and $(gf)^!\co_Z$ with
$\Omega^n_{Y/Z}[n]$, $\Omega^m_{X/Y}[m]$, and 
$\Omega^{m+n}_{X/Z}[m+n]$ respectively, then $\chi_{[f,g]}$ transforms to
the map $f^*\mu\otimes \nu \mapsto \nu\wedge f^*\mu$  (see \Tref{thm:fubini}).
In fact we show this at the level of formal schemes, and formal schemes enter in an essential
way in our proof (via transitivity for residues 
\cite[(8.3.2)]{fub-abs} and \Tref{thm:res-res})
even when $X$ and $Y$ are ordinary schemes. The proof is carried out in \Sref{s:transitivity}.

 The second situation of importance occurs when the smooth map
 $f\colon X\to Y$ factors as $f=\pi\smcirc i$, where $i\colon X\hookrightarrow P$ 
 is a closed immersion say of codimension $d$, 
 and $\pi\colon P\to Y$ is smooth of relative dimension $n+d$. The concrete expression
 for $\chi_{i, \pi}$ then is governed by
 \[i^*(\eta\wedge\wdd{t_1}{t_d})\otimes {\bf 1/t}\mapsto i^*\eta\] 
where $\eta$ is a section of $\Omega^{n+d}_{P/Y}$, $t_i\in\Gamma(P,\,\co_P)$, $i=1, \dots, d$, 
are sections which cut out $X$ and ${\bf 1/t}$ is a well-defined generating section, 
depending upon ${\bf t}=(t_1, \dots, t_d)$, of the top exterior
product $\wedge^d\eN$ of the normal bundle  $\eN$ of $X$ in $P$, which exterior product, by the 
\emph{fundamental local isomorphism} is identified with $f^!\co_P[d]$. The proof of this
concrete representation of $\chi_{[i,\pi]}$is carried
out in \Ssref{ss:restriction}.

In both situations, we need the residue formula \eqref{intro:std-res} for residues along
sections of smooth maps. We turn this around later, and use the concrete expressions for
$\chi_{[-,-]}$ to arrive at formulas for $\res{Z}$ for smooth maps 
$f\colon X\to Y$ when $Z\to Y$ is not an isomorphism (but is finite).

There is one interesting way in which \eqref{intro:std-res=1} brings in concrete
answers. Let $A$ be a ring, and $C=A[T_1, \dots, T_n]/(f_1, \dots, f_n)$ be a finite
flate algebra over $A$. Let $Z=\Spec{\,C}$, $X=\mathbb{A}_A^n$ and $Y=\Spec{A}$.
Let $I={\bf f}A[{\bf T}]$, so that $I$ is the ideal of $Z$ in the polynomial ring
$A[{\bf T}]$. 
By the general calculus of generalised fractions, if $p({\bf T})\in A[{\bf T}]$ then
the element 
$\Bigl[\begin{smallmatrix}p({\bf T})\wdd{T_1}{T_n}\\f_1, \dots, f_n\end{smallmatrix}\Bigr]\in \Hr^n_I(\Omega^n_{A[{\bf T}]/A}$ depends only on the image of $p({\bf T})$ in $C$.
We show that the map 
\[c\mapsto \res{Z}\begin{bmatrix} p({\bf T})\wdd{T_1}{T_n}\\f_1, \dots f_n\end{bmatrix}\]
with $p({\bf T})$ a pre-image of $c$, is the Tate trace described in \cite[Appendix]{tate}. 
We prove this in \Tref{thm:tate2}, and \eqref{intro:std-res=1} plays an important role. 
The point is, knowing the residue in a very special situation allows us to deduce 
formulas for residues in many other situations.

Perhaps the most important way that that \eqref{intro:std-res} comes into play is that
it characterises the Verdier isomorphism (or more accurately the fundamental class).
Continuing with the situation where $f\colon X\to Y$ is smooth of relative dimension $n$,
suppose we have some isomorphism $\psi\colon \Omega^n_{X/Y}[n]\iso f^!\co_Y$. When
is $\psi$ Verdier's isomorphism? The answer is, if and only if, for every \'etale base change
$u^*Y'\to Y$ and every section $Z$ of the base change map
$f'\colon X\times_YY'\to Y'$, the composite
\[\Rr^n_Z f'_*\Omega^n_{X\times_YY'/Y'} \xrightarrow[\text{via $\psi$}]{\Iso}
\Hr^0(\R_Z f'_*{f'}^!\co_{Y'}) \xrightarrow{\Tr{f', Z}} \co_{Y'}\]
is given by \eqref{intro:std-res}. Note that the first isomorphism involves flat
base change for $f^!$. The precise statement is given in \Tref{thm:main}. This
characterisation of Verdier's isomorphism allows us to relate the fundamental
class with the regular differentials of Kunz. We work this out in \Sref{s:reg-diffs}.
We give a different proof later of the relationship between the fundamental class and regular
differentials.

\subsection{Trace for finite flat maps} Suppose the smooth map $f\colon X\to Z$ of
relative dimension $n$ can be factored
as $f=g\smcirc h$, where $h\colon X\to Y$ is finite and $g\colon Y\to Z$ is smooth of relative
dimension $n$ (so that $h$ is in fact flat). Then the composite 
$h_*f^! \cong h_*h^!g^! \xrightarrow{\Tr{h}} g^!$, gives, via Verdier's isomorphism for $f$ and
$g$ a map 
\[\vin{h}\colon h_*\Omega^n_{X/Z}\lra \Omega^n_{Y/Z}.\]
In \cite{kd}, Kunz, based on a suggestion by Lipman (who in turn was influenced by residue
formulas stated without proof in \cite[Chap.\,III, \S\,9]{RD}) defined an explicit trace
$\sigma_h\colon h_*\Omega^n_{X/Z}\lra \Omega^n_{Y/Z}$. We show that $\vin{h}=\sigma_h$.
In fact, we use the two concrete versions of transitivity that we mention above. Assuming $h$
factors as a closed immersion $i\colon X\to P$ followed by a smooth map $\pi\colon P\to Y$,
where $P$ is an open subscheme of ${\mathbb A}^{n+d}_Y$ and $\pi$ the structure map,
(a situation we can achieve, retaining finiteness of $h$, if we pass to completions of local
rings of points on $Y$) then the assertion $\vin{h}=\sigma_h$ amounts to the compatibilities
between the abstract transitivity maps $\chi_{[h,g]}$, $\chi_{[i,\pi]}$, $\chi_{[\pi, g]}$,
and $\chi_{[i, g\pi]}$ given in 
\cite[Prop.-Def.\,7.2.4\,(ii)]{fub-abs} or in
\cite[p.\,238]{notes}. The map $\vin{h}$ occurs in formula (R10) for residues, and it is satisfying
that there is a more explicit description of it in terms of the Kunz-Lipman trace $\sigma_h$.

\subsection{Regular Differential Forms} The regular differentials of Kunz and Waldi 
developed in \cite{kw} is a vast generalisation of Rosenlicht's differentials for singular
curves \cite{rosen}.
We have already alluded to the connection
between the regular differential forms and Verdier's isomorphism.
Regular differntial forms
are defined when $f\colon X\to Y$ is a generically smooth equidimensional
map between excellent schemes having no embedded points. In such a case,
if $X^{\mathrm{sm}}$ is the smooth locus of $f$,
and $f^{\mathrm{sm}}\colon X^{\mathrm{sm}} \to Y$ the restriction of $f$, there is
an isomorphism $\Omega^n_{X^{\mathrm{sm}}/Y}[n] \to (f^{\mathrm{sm}})^!$.
The isomorphism is based on the construction of regular differentials in
\cite{kw} and the principal results of \cite{ajm}. What we show in this paper is
that this isomorphism is Verdier's isomorphism. We give two proofs. The first
using the characterisation of Verdier's isomorphism via \eqref{intro:std-res} that
we alluded to before. The other, more satisfying, proof relies on the equality of
traces $\vin{h}=\sigma_h$ for finite flat maps $h$ between schemes smooth over a base
that we discussed above. (See \Ssref{ss:reg-diff-again}.)

\section{\bf Preliminaries}\label{s:prelims}

\subsection{Flat base change}
\label{subsec:base-ch}
Suppose we have a cartesian square $\mathfrak s$ of noetherian formal schemes
\[
\xymatrix{
\V \ar[d]_{g}\ar[r]^v \ar@{}[dr]|\square& \X \ar[d]^{f} \\
\W \ar[r]_u & \Y
}
\]
with $f$ in $\bbG$ and $u$ flat. 
The \textit{flat-base-change theorem} 
for $\ush{\boldsymbol{(-)}}$ (see  \cite[\S\S\,3.2]{fub-abs}) states that is
if $\F \in \Dc^+(\Y)$, or if $u$ is open or if $\V$ is an ordinary scheme,
we have an isomorphism (see \cite[(3.2.2)]{fub-abs}):
\stepcounter{thm}
\begin{equation*}\label{iso:bc-sharp}\tag{\thethm}
v^*\ush{f}\F\iso \ush{g}u^*\F.
\end{equation*}
If $f$ is pseudo-proper, \eqref{iso:bc-sharp} is proved in
\cite[Theorem 8.1, Corollary 8.3.3]{dfs}.

There is a somewhat more explicit version \eqref{iso:bc-sharp} when $f$ is a 
{\em regular immersion} i.e., $\X$ is a closed subscheme of $\Y$ given by the
vanishing of a regular sequence. In such a case let $\I$ be the ideal sheaf of $\X$ in $\Y$, 
and suppose the co-dimension
of the immersion is $r$, i.e., $\I$ is locally generated by regular sequences of length $r$.
As in \cite[(C.2.10)]{fub-abs}, for $\eF\in\Dqc(\X)$, we write
\stepcounter{thm}
\begin{equation*}\label{def:btrg}\tag{\thethm}
f^\btrg(\eF)\set \bL f^*\eF\overset{\bL}{\otimes}_{\co_\X}(\wedge^r_{\co_\X}(\I/\I^2)^*[-r]).
\end{equation*}
According \cite[(C.2.13)]{fub-abs}, for $\eF\in\Dc^+(\Y)$ we have a functorial isomorphism
\stepcounter{thm}
\begin{equation*}\label{iso:eta'-i}\tag{\thethm}
\eta'_f(\eF)\colon f^\btrg(\eF) \iso \ush{f}(\eF).
\end{equation*}
Now suppose we have a cartesian diagram $\mathfrak s$ of formal schemes
\stepcounter{thm}
\[
\begin{aligned}\label{diag:bc0}
\xymatrix{
\W' \ar@{}[dr]|{\square} \ar[d]_{\kappa} \ar[r]^j & \W \ar[d]^{\kappa_{{}_0}}\\
\X' \ar[r]_i & \X
}
\end{aligned}\tag{\thethm}
\]
such that $i$ is a regular immersion 
and $\kappa_{{}_0}$ is the completion of $\X$ with respect to a closed subscheme given by
a coherent ideal. Then according to \cite[(C.4.2)]{fub-abs} we have
an isomorphism
\stepcounter{thm}
\begin{equation*}\label{iso:BCf}\tag{\thethm}
\kappa^*i^\btrg \iso j^\btrg\kappa_{{}_0}^*.
\end{equation*}

The compatibility of \eqref{iso:BCf} with \eqref{iso:bc-sharp} (with $f=i$, $g=j$, $u=\kappa_{{}_0}$,
and $v=\kappa$) given the isomorphisms
 $\eta'_i$ and $\eta'_j$ (see \eqref{iso:eta'-i}), is proven in \cite[Prop.\,C.4.3]{fub-abs}.

Further properties of the base-change map are explored \cite[\S\S\,A.1]{fub-abs}.

\subsection{} 
If $f\colon\X\to \Y$ is  a map in $\bbG$ which is {\emph{formally \'etale}}
and $\F \in \Dc^+(\Y)$, then we have
an isomorphism
\stepcounter{thm}
\begin{equation*}\label{eq:gm}\tag{\thethm}
f^*\F \iso \ush{f}\F
\end{equation*}
which is pseudofunctorial over the category of formally \'etale maps 
\cite[(3.1.3)]{fub-abs}.

\subsection{Completions and direct image with support} Let $Z$ be a closed subscheme
of an ordinary scheme $X$ and $\kappa\colon\X\to X$ the {\emph{completion of $X$ along $Z$}}.
The isomorphism $\kappa_*\R\iGp{\X}\kappa^*\iso \R\iG{Z}$ gives rise to 
isomorphisms (one for every $j$)
\stepcounter{thm}
\begin{equation*}\label{iso:loc-coh}\tag{\thethm}
\Rr^j_Z f_*\eF \iso \Hr^j({\R\wid{f}}_*\R\iGp{\X}\kappa^*\eF)=\Rp{\X}^j{\wid{f}}_*\kappa^*\eF
\end{equation*}
which are functorial in $\eF$ varying over quasi-coherent $\co_X$-modules. 
In affine terms, if $X=\Spec{\,R}$, $M$ an $R$-module, and $Z$ is given
by the ideal $I$, then writing $\wid{R}$ for the $I$-adic completion of $R$, and $J=I\wid{R}$,
the above isomorphism is the well-known one
\[\Hr^j_I(M) \iso\Hr^j_J(M\otimes_R\wid{R}).\]

In the above situation, suppose $f\colon X\to Y$ is a separated finite-type map of ordinary
schemes such that the induced map $Z\to Y$ is {\emph{proper}}. 
Let $\wid{f}=f\<\smcirc\kappa$. Then according to \cite[(A.3.1)]{fub-abs}
we have a functorial isomorphism
\stepcounter{thm}
\begin{equation*}\label{iso:iGZ-iGpX}\tag{\thethm}
 \Rfs \R\iG{Z}\ush{f} \iso \R{\wid{f}}_*\R\iGp{\X}\ush{\wid{f}}.
\end{equation*}

\subsection{Abstract traces and residues} 
Suppose $f\colon \X\to\Y$ is a pseudo-proper map. As in
\cite[(3.1.2)]{fub-abs}
 we define the {\emph{trace map}} associated to $f$,
\stepcounter{thm}
\begin{equation*}\label{map:Tr-f}\tag{\thethm}
\Tr{f} \colon \R f_*\R\iGp{\X}\ush{f} \to {\bf 1}.
\end{equation*}
to be the co-adjoint unit associated to the right adjoint $\ush{f}$ of $\R f_*\R\iGp{\X}$.

In \cite[Def.\,4.1.2]{fub-abs} we defined a
map of formal schemes $f\colon \X \to \Y$ to be {\em Cohen-Macaulay \textup{(CM)} of
relative dimension $r$} if it is flat, locally in $\bbG$ with $\Hr^i(\ush{f}\co_\Y)=0$ for
$i\neq -r$, and $\omgs{f}\set \Hr^{-r}(\ush{f}\co_\Y)$ is flat over~$\Y$. The coherent
$\co_\X$-module $\omgs{f}$ is called the \emph{relative dualizing sheaf} for the CM map
$f$. If such a map~$f$ is already in~$\bbG$, we shall make the identification 
$\ush{f}\co_\Y=\omgs{f}[r]$. 

For a CM map $f\colon\X\to \Y$ of relative dimension $r$ which is {\emph{pseudo-proper}},
\stepcounter{thm}
\begin{equation*}\label{map:tin}\tag{\thethm}
\tin{f}\colon \Rp{\X}^r f_*\omgs{f} \to \co_\Y
\end{equation*}
will denote the {\emph{abstract trace map on $\Rp{\X}^rf_*\omgs{f}$}} and is defined by
the formula
\[\tin{f}=\Hr^0(\Tr{f}(\co_\Y))\]
as in \cite[(5.1.2)]{fub-abs}.

Next, let $f\colon X\to Y$ be a separated map of finite-type between ordinary schemes, and $Z$  a 
closed subscheme of $X$ which is proper over $Y$.  We recall the
notion of the trace of $f$ along $Z$ from \cite[\S\S\,3.3]{fub-abs}.
Now, the completion map $\kappa\colon \X\to X$ of $X$ along $Z$,
is formally \'etale and affine and the composition $\wid{f} \set f\kappa$ is pseudo-proper. 
We define the {\emph{trace map for $f$ along $Z$}},
\stepcounter{thm}
\begin{equation*}\label{map:trZ}\tag{\thethm}
\Tr{f,Z}\colon \R_Zf_*\ush{f} \to {\bf 1}
\end{equation*}
to be the composite
\[\R f_*\R\iG{Z}\ush{f}  \xrightarrow[\eqref{iso:iGZ-iGpX}]{\Iso} \R{\wid{f}}_*\R\iGp{\X}\ush{\wid{f}}
\xrightarrow{\Tr{\wid{f}}} {\bf 1}. \]

Suppose $f$ above is {\emph{Cohen-Macaulay of relative dimension $r$}}. 
As in \cite[(5.2.2)]{fub-abs}, 
the {\em abstract residue along $Z$} 
\stepcounter{thm}
\begin{equation*}\label{def:ares}\tag{\thethm}
\ares{Z}\colon \Rr^r_Zf_*\omgs{f} \to \co_Y
\end{equation*}
is defined as the composite 
\[\Rr^r_Zf_*\omgs{f} \xrightarrow[\Hr^0\eqref{iso:iGZ-iGpX}]{\Iso} \Rp\X^r{\wid{f}}_*\omgs{\wid{f}} 
\xrightarrow{\tin{\wid{f}}} \co_Y.\]
It is clear from the definitions that
\[\ares{Z}=\Hr^0(\Tr{f,Z}(\co_Y)).\]

\section{\bf Verdier's isomorphism} 
\subsection{} Let $f\colon \X\to \Y$ be a smooth map of relative dimension $r$ between
(formal) schemes. Assume $f$ is a composite of compactifiable maps. 
Set $\X''\set \X\times_\Y\X$ and let $\De\colon \X\to \X''$ the diagonal immersion, which
is closed by our hypotheses. Denote by $p_{{}_1}$ and $p_{{}_2}$ the two projections
from $\X''$ on to $\X$, and by $\eN_{\<\<\<\<{}_\De}$ the locally free $\co_\X$-module corresponding
to the ``normal bundle" of the regular immersion~$\De$. In other words, if $\I_\De$ is the
ideal sheaf of $\X$ in $\X''$, then $\eN_{\<\<\<\<{}_\De}=(\I_\De/\I_\De^2)^*$, the dual of
$\I_\De/\I_\De^2$. As in  \cite[(C.2.8)]{fub-abs} and \eqref{def:btrg},
set 
\[\wnor{\De}=\wedge^r_{\co_\X}\eN_\De\] 
and
\[\De^\btrg=\bL\De^*({\boldsymbol{-}})\overset{\bL}\otimes_{\co_\X}\wnor{\De}[-r]\] 
We then have an isomorphism
\stepcounter{thm}
\begin{equation*}\label{iso:pre-verdier}\tag{\thethm}
\ush{f}\co_\Y\overset{\bL}\otimes_{\co_\X}\wnor{\De}[-r] \iso \co_\X
\end{equation*}
defined by the commutativity of the following diagram:
\stepcounter{thm}
\begin{equation*}\label{diag:pre-verdier}\tag{\thethm}
\xymatrix{
\ush{f}\co_\Y\overset{\bL}\otimes_{\co_\X}\wnor{\De}[-r] 
\ar[d]^{\rotatebox{90}{\makebox[0.1cm]{\Iso}}}_{\eqref{iso:pre-verdier}}
& \De^\btrg p_{{}_2}^*\ush{f}\co_\Y \ar[l]_-\Iso \ar[r]^{\Iso}_{\eqref{iso:bc-sharp}} 
& \De^\btrg\ush{p_{{}_1}}\co_\X 
\ar[d]_{\,\rotatebox{-90}{\makebox[-0.1cm]{\Iso}}}^{\eta'_\De} \\
\co_\X & &\ush{\De}\ush{p_{{}_1}}\co_\X \ar[ll]_\Iso 
}
\end{equation*}
The map $\eta'_\De$ is as in \eqref{iso:eta'-i}.
The unlabelled arrow on the top row is the one arising from $\bL\De^*p_{{}_2}^* \iso {\bf 1}$ and
the one on the bottom row from the functorial isomorphism 
$\ush{\De}\ush{p_{\<{}_1}}\iso \ush{{\bf 1}_\X}$. 


Writing $\eL^*$ for the dual of an invertible $\co_\X$-module $\eL$ we
see that $\wnor{\De}=\omega_f^*$.

Using this and \eqref{iso:pre-verdier} one deduces, as Verdier did in 
\cite[p.\,397, Theorem\,3]{verdier}, that
$\ush{f}\co_\Y$ and $\omega_f[r]$ are isomorphic. However, there is some ambiguity
about the exact isomorphism (see the discussion around (7.2) in
 p.\,758 of \cite{cm}). But at the very least we note that $f$ is 
 Cohen-Macaulay.
 We give the isomorphism we will
 work with in \Dref{def:verdier} after some necessary preliminaries.
 
 As usual, let $\omgs{f}=\Hr^{-r}(\ush{f}\co_\Y)$ and make the identification
 \[\ush{f}\co_\Y=\omgs{f}[r].\]
  Applying $\Hr^0$ to \eqref{iso:pre-verdier} we get
 an isomorphism 
 \stepcounter{thm}
 \begin{equation*}\label{iso:v-pairing}\tag{\thethm}
 \omgs{f}\otimes_{\co_\X}\omega_f^* \iso \co_{\X}
 \end{equation*}
 Let 
\stepcounter{thm}
 \begin{equation*}\label{def:verd-bar}\tag{\thethm}
{\bar{\bf v}}_{\<\<{}_f}(=\bar{\bf v})\colon \omega_f \iso \omgs{f}
 \end{equation*}
be the canonical isomorphism induced by \eqref{iso:v-pairing}.

\begin{defi}
\label{def:verdier} 
The {\em Verdier isomorphism} for the smooth map $f$ is the isomorphism
\stepcounter{sth}
\[
{\bf v}_{\<\<{}_f}(={\bf v})\colon\omega_f[r]\iso \omgs{f}[r]=\ush{f}\co_\Y.
\]
given by ${\bf v}_{\<\<{}_f}={\bar{\bf v}}_{\<\<{}_f}[r]$.
\end{defi}
We will often refer to $\bar{\bf v}_{\<\<{}_f}$ also as the Verdier isomorphism. Indeed 
${\bf v}_{\<\<{}_f}$ and ${\bar{\bf v}}_{\<\<{}_f}$ determine each other.


\begin{rem}
\label{rem:theta}
 The isomorphism $p_{{}_2}^*\ush{f}\co_\Y \iso \ush{p_{{}_1}}\co_\X$ of
\eqref{iso:bc-sharp}
induces (on applying $\Hr^0$) an isomorphism
\[\theta\colon p_{{}_2}^*\omgs{f} \iso \omgs{p_{{}_1}}.\]
Note that the original
isomorphism $p_{{}_2}^*\ush{f}\co_\Y \iso \ush{p_{{}_1}}\co_\X$ is $\theta[r]$ under
the identifications we have agreed to make throughout, namely, $\ush{f}\co_\Y= \omgs{f}[r]$
and $\ush{p_{{}_1}}\co_\X = \omgs{p_{{}_1}}[r]$. Applying the functor $\Hr^0$ to the commutative
diagram \eqref{diag:pre-verdier} we get the following commutative diagram, showing the
relationship between the pairing \eqref{iso:v-pairing} and maps of the form $\tau_h$ defined
in \cite[(5.3.2)]{fub-abs} (below, $h$ is the identity map).
\stepcounter{sth}
\begin{equation*}\label{diag:pairing-tau}\tag{\thesth}
\xymatrix{
\omgs{f}\otimes_{\co_\X}\omega_f^*
\ar[d]^{\,\rotatebox{-90}{\makebox[-0.1cm]{\Iso}}}_{\eqref{iso:v-pairing}} \ar@{=}[r]
& \De^*(p_{{}_2}^*\omgs{f}\otimes_{\co_{\X''}}\De_*\omega_f^*)
\ar[d]_{\,\rotatebox{-90}{\makebox[-0.1cm]{\Iso}}}^{{\text{via $\theta$}}} \\
\co_\X &  \De^*(\omgs{p_{{}_1}}\otimes_{\co_{\X''}}\De_*\omega_f^*) \ar[l]_-\Iso^-{\ttr{{\bf 1},p_1,\De}} 
}
\end{equation*}
Here the map on the bottom row is as in 
\cite[(5.3.2)]{fub-abs}, with $h={\bf 1}_\X$, $i=\De$,
and $f=p_{{}_1}$. It is an isomorphism because $h={\bf 1}_\X$ is an isomorphism.
\end{rem}

\begin{defi} Suppose $f\colon \X\to \Y$ is pseudo-proper and smooth of relative dimension $r$.
The {\em Verdier integral} (or simply the {\em integral})
\stepcounter{thm}
\begin{equation*}\label{map:vin}\tag{\thethm}
\vin{f}\colon \Rp{\X}^r f_*\omega_f \to \co_\Y
\end{equation*}
is the composite 
$\Rp{\X}^r f_*\omega_f \xrightarrow{\bf v} \Rp{\X}^r f_*\omgs{f} \xrightarrow{\tin{f}} \co_\Y$.
If in the above situation, $\X=\Spf{(R,\,J)}$ and $\Y=\Spf{(A,I)}$, then we write
\stepcounter{thm}
\begin{equation*}\label{map:pint}\tag{\thethm}
\vin{R/A}\colon \Hr^r_J(\omega_{R/A})\to A
\end{equation*}
for the global sections of
$\vin{f}$. Here $\omega_{R/A}$ is the $r$-th exterior power of the universally finite module of
differentials for the $A$-algebra $R$. If we wish to emphasise the adic structure on $R$ and
$A$, we will use the inconvenient notation $\vin{(R,J)/(A,I)}$ for $\vin{R/A}$.
\end{defi}

\begin{rem}While we have defined $\vin{f}$ in general, our interest is really in the case 
where $\Rp{\X}^j(\eF)=0$ for every 
$j>r$ and every $\eF\in\Avc(\X)$, for then $(\omega_f, \vin{f})$ represents
the functor $\Hom_\Y({\Rp{\X}^r(\boldsymbol{-}}),\, \co_\Y)$
on coherent $\co_\X$-modules 
(see the result in \cite[Cor.\,(5.1.4)]{fub-abs}). Even here
the notion is most useful in this paper when $\Y$ is an ordinary scheme and either
$\X$ is also ordinary (and hence proper over $\Y$) or else $\Y=\Spec{\,A}$ and
$\X =\Spf{\,R}$ where $R$ is an adic ring, one of whose defining ideals $I$ is
generated by a quasi-regular sequence of length $r$ and such that $R/I$ is finite
and flat over $A$. 
\end{rem}

\subsection{Local description of Verdier's isomorphism}\label{ss:verd-loc}
In the above situation suppose
$\X=\Spf{R}$, $\Y=\Spf{A}$, so that $\X''=\Spf{R''}$  where $R''=R\wid{\otimes}_AR$ is the complete
tensor product of $R$ with itself over $A$.
The diagonal map $\De\colon \X' \hookrightarrow \X''$ corresponds to the surjective
map $R''\to R$ given by $t_1\otimes t_2\mapsto t_1t_2$. Let us assume that the kernel of this map,
i.e., the ideal $I$ of the diagonal immersion, is generated by $r$-elements $\{s_1,\dots, s_r\}$.
Since $R$ is smooth over $A$ of relative dimension $r$, the sequence ${\bf s}=(s_1,\dots, s_r)$
is necessarily a $R''$-sequence. This condition on the diagonal is locally (in $\X$ and $\X''$) 
always achievable.
 
 Let $R_1$ and $R_2$ be the two $R$-algebra structures on $R''$ corresponding to the projections
$p_{{}_i}\colon \X''\to X$, $i=1,2$. For specificity, if $a\in R$, then the $R$-algebra structure on
$R'_1$ is given by $a(b\otimes c)=(ab)\otimes c$ whilst on $R_2$ it is given by 
$a(b\otimes c)=b\otimes(ac)$.
Let $\omgs{R/A}$, $\omgs{R_i/A}$, $\omega_{R/A}$, 
$\omega_{R_i/A}$ be the global sections of $\omgs{f}$, $\omgs{p_{{}_i}}$, $\omega_f$, 
$\omega_{p_{{}_i}}$
respectively, where $i\in\{1,\,2\}$. Similarly, {\em Verdier's isomorphism} in this context
is the isomorphism
\[\bar{\bf v}_{\<\<{}_{R/A}}\colon \omega_{R/A} \iso \omgs{R/A}\]
obtained by taking global sections of $\bar{\bf v}_{\<\<{}_{f}}\colon \omega_f \iso \omgs{f}$.

The isomorphism \eqref{iso:v-pairing} is equivalent to the isomorphism of finitely
generated $R$-modules
obtained by taking global sections:
\stepcounter{thm}
\begin{equation*}\label{iso:v-pairing-aff}\tag{\thethm}
\omgs{R/A}\otimes_R\omega_{R/A}^* \iso R.
\end{equation*}

Here is the promised local description of Verdier's isomorphism.
The module of differentials $\omega_{R/A}=\wedge^r_RI/I^2$ is a free rank one
$R$-module with generator 
\[{\overline{{\mathrm d}{\bf s}}}\set(s_1+I^2)\wedge\dots\wedge(s_r+I^2).\] 
Let ${\bf 1}/{\bf s}$ be the element of  $\wI{R}{I}=\omega_{R/A}^*=\Hom_R(\omega_{R/A},\,R)$ 
which sends ${\overline{{\mathrm d}{\bf s}}}$
to~$1$, i.e., it is the generator of the rank one free module $\wI{R}{I}$ which is dual 
to~${\overline{{\mathrm d}{\bf s}}}$. 


\begin{prop}\label{prop:ds/s-1} In the above situation we have the following:
\begin{enumerate}
\item[(a)] Let $\nu_0({\bf s})\in\omgs{R/A}$  be the unique element
such that  $\nu_0({\bf s})\otimes{\bf 1}/{\bf s}\mapsto 1$  under \eqref{iso:v-pairing-aff}.
Verdier's isomorphism $\bar{\bf v}_{\<\<{}_{R/A}}\colon\omega_{R/A}\iso \omgs{R/A}$ is given
by the formula
\[\bar{\bf v}_{\<\<{}_{R/A}}(r\,\overline{{\mathrm d}{\bf s}}) = r\,\nu_0({\bf s}) \qquad (r\in R).\]
\item[(b)] Suppose further that the adic rings $A$ and $R$ have
discrete topology so that $\Spf{\,A}=\Spec{\,A}$, $\Spf{\,R}=\Spec{\,R}$ and $A\to R$ is
of finite type. The following formula holds:
\[
\ares{\De,\,p_{{}_1}}\begin{bmatrix} {\bar{\bf v}}_{\<\<{{}_{R_1/R}}}(\wdd{s_1}{s_r})\\
s_1,\,\dots,\,s_r \end{bmatrix} = 1.
\]
\end{enumerate}
\end{prop}

\noindent{\em Remarks:} Here $\wdd{s_1}{s_r} \in \omega_{R_1/R}$ and the 
notation $\ares{\De,\,p_{{}_1}}$ is to indicate that the residue
is to be taken for the map $p_{{}_1}$ and not for $p_{{}_2}$. The hypotheses in part\,(b)
regarding the adic topologies on $A$ and $R$ is there because we need the result that Verdier's 
isomorphism
is compatible with base change. This is one of the main results of \cite{cm} (see \cite{cm}, p.740, Theorem 
2.3.5\,(b)]). Unfortunately the results
in \cite{cm} are for maps between ordinary schemes. Since certain special compactifications
are locally used, and these are unavailable for arbitrary formal schemes, we decided it is
best not pursue these issues in this paper, except in the following special case. Suppose the
base change is flat. Then the proof in  \cite{cm} works mutatis mutandis, and we see that  
Verdier's isomorphism is compatible with flat base change whether we are working with ordinary
schemes or formal schemes. See \Tref{thm:verd-flat-bc}.

\proof
Part\,(a) is an immediate consequence of
the definition of $\bar{\bf v}_{\<\<{}_{R/A}}$ in \eqref{def:verd-bar}.
It remains to prove part\,(b).

Let us save on notation and write 
\[\theta\colon \omgs{R/A}\otimes_RR_2 \iso \omgs{R_1/R}\]
for the $R''$-isomorphism corresponding to 
$\theta\colon p_{{}_2}^*\omgs{f}\iso \omgs{p_{{}_1}}$ of \Rref{rem:theta}.
Then the affine version of the commutative diagram \eqref{diag:pairing-tau} is the commutative
diagram
\stepcounter{thm}
\begin{equation*}\label{diag:pairing-tau-aff}\tag{\thethm}
{\xymatrix{
\omgs{R/A}\otimes_R\omega_{R/A}^*
\ar[d]^{\rotatebox{90}{\makebox[0.1cm]{\Iso}}}_{\eqref{iso:v-pairing-aff}} \ar@{=}[r]
& (\omgs{R/A}\otimes_RR_2)\otimes_{R''}\omega_{R/A}^*
\ar[d]_{\,\rotatebox{-90}{\makebox[-0.1cm]{\Iso}}}^{{\text{via $\theta$}}} \\
R &  \omgs{R_1/R}\otimes_{R''}\omega_{R/A}^*
\ar[l]_-\Iso^-{\ttr{R/R, R_1}} 
}}
\end{equation*}
Since $\nu_0({\bf s})\otimes{\bf 1}/{\bf s} \mapsto 1$  under \eqref{iso:v-pairing-aff},
from the commutative diagram \eqref{diag:pairing-tau-aff} we get
$\ttr{R/R, R_1}(\theta(\nu_0\otimes 1)\otimes{\bf 1}/{\bf s})=1$. Moreover, by part\,(a), 
$\nu_0= {\bar{\bf v}}_{\<\<{{}_{R/A}}}({\overline{{\mathrm d}{\bf s}}})$. Thus 
\cite[Prop.\,5.4.4]{fub-abs} gives us
\[
\ares{\De,\,p_{{}_1}}\begin{bmatrix}
\theta(({\bar{\bf v}}_{\<\<{{}_{R/A}}}({\overline{{\mathrm d}{\bf s}}}))\otimes 1))\\
s_1,\,\dots,\,s_r \end{bmatrix} = 1.
\]

Next, if $M$ and $N$ are finitely generated
modules over $R$ and $\varphi\colon M\to N$
a map of $R$-modules, then we denote the map 
$\varphi\otimes_RR_i$ by $p_{{}_i}^*(\varphi)$. One checks easily that in
$N\otimes_RR_2$ we have the following equality for $m\in M$ and $\varphi\in\Hom_R(M,\,N)$:
\[
(p_{{}_2}^*(\varphi))(m\otimes 1)= (\varphi(m))\otimes 1. 
\]
It is immediate that
\[
\ares{\De,\,p_{{}_1}}\begin{bmatrix}
(\theta\smcirc (p_{{}_2}^*({\bar{\bf v}}_{\<\<{{}_{R/A}}})))({\overline{{\mathrm d}{\bf s}}}\otimes 1))\\
s_1,\,\dots,\,s_r \end{bmatrix} = 1.
\]
Now if $s=\sum_i a_i\otimes b_i \in I\subset R''=R\otimes_AR$, then one checks from the definitions
that as elements of  the $R$-module $I/I^2=\Omega^1_{R/A}$, we have the
equality $s+I^2 = \sum_ia_i{\mathrm{d}}b_i$. This means in particular
that in the $R''$-module $\Omega^1_{R/A}\otimes_RR_2=\Omega^1_{R_1/R}$ we have
$(s+I^2)\otimes 1 = {\mathrm{d}}s$.
It is immediate from here that 
${\overline{{\mathrm d}{\bf s}}}\otimes 1= \wdd{s_1}{s_r}.$

We will be done if we can show that 
$\theta\smcirc (p_{{}_2}^*({\bar{\bf v}}_{\<\<{{}_{R/A}}}))={\bar{\bf v}}_{\<\<{{}_{R_1/R}}}$. This
statement about the compatibility of Verdier's isomorphism with base change follows from
\cite[p.740,\,Theorem\,2.3.5]{cm} (see also [\textit{ibid.},\,pp.739--740,\,Remark\,2.3.4]). Incidentally,
this is where we need our hypothesis that our formal schemes are ordinary schemes and our map is
of finite type. 
\qed

\begin{rems}\label{rem:v-local}{\em Two observations are worth making.

1) $\bar{\bf v}=\Hr^{-r}({\bf v})$.

2) If $\U$ is an open subscheme of $\X$, and 
$f_{{}_\U}\colon\U\to \Y$ is the structural morphism on $\U$,
then we have a natural isomorphism $\ush{f}\vert_\U \iso \ush{f_{{}_\U}}$ from the main
results of \cite{pasting}, whence an isomorphism $\omgs{{}_f}\vert_\U\iso \omgs{{}_{f_{{}_\U}}}$.
From the definitions of ${\bf v}_{\<\<\<{{}_f}}$ and ${\bf v}_{\<\<\<{{}_{f_{{}_\U}}}}$ it is easy to see
that the composition of isomorphisms $\omega_{f_{{}_\U}} = \omega_{{}_f}\vert_\U \iso
\omgs{{}_f}\vert_\U\iso \omgs{{}_{f_{{}_\U}}}$ is ${\bf v}_{\<\<\<{{}_{f_{{}_\U}}}}$, where the
first arrow is ${\bf v}_{\<\<\<{{}_f}}\vert_\U$ and the second the just mentioned isomorphism.
}
\end{rems}


\subsection{Compatibility of Verdier's isomorphism with completions} We now wish to show
the compatibility of Verdier's isomorphism with completion. More precisely if $f\colon \X\to \Y$ 
is a smooth map and $\wid{f}\colon \W\to \Y$ its ``completion" along a closed subscheme of $\X$,
then Verdier's isomorphism (i.e., \eqref{def:verdier}) for $\wid{f}$ is the ``completion"
of the Verdier isomorphism for $f$. The formal statement is given in  
\Tref{thm:verd-complete}.

Let $f\colon \X\to \Y$ be a smooth map of relative dimension $r$ between formal schemes.
Suppose $\I$ is a defining ideal of $\X$ and
$\J\subset\co_\X$ a coherent ideal containing $\I$ (so that $\J$ is the ideal of an ordinary 
scheme $Z$ which is a closed subscheme of $\X$). Let $\W$ be the completion
of $\X$ along $\J$ (i.e., along $Z$). Let $\kappa\colon \W\to \X$ be the completion map and 
$\wid{f}\colon \W\to \Y$ the composite $\wid{f}=f\smcirc\kappa$. We wish to show that the Verdier
isomorphism for $\wid{f}$ ``is" $\kappa^*$ of the Verdier isomorphism for $f$.
As before let $\X''=\X\times_{\Y}\X$, and $\De\colon \X\to \X''$ the diagonal
immersion. Let $\W''=\W\times_{\Y}\W$, and let $\delta\colon\W\to\W''$ be the diagonal immersion.

Let $\widetilde\kappa\colon \W''\to \X''$ be the map  $\kappa_{{}_2}=\kappa\times\kappa$. As usual,
we have projections $p_{{}_i}\colon \X''\to \X$ and $\pi_{{}_i}\colon \W''\to \W$ for $i=1,2$. 
The following commutative diagrams may help the reader map the relative
positions of the schemes and maps involved:

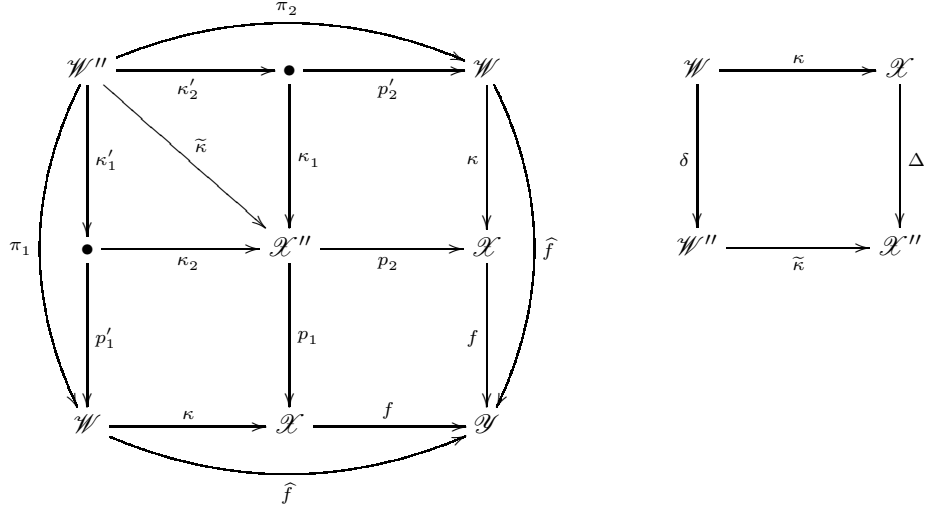
\[
\xymatrix{
\W'' \ar[rr]_{\kappa_2'} \ar[dd]^{\kappa_1'} \ar@/^1.5pc/[rrrr]^{\pi_2^{}} \ar@/_1.5pc/[dddd]_{\pi_1^{}} 
\ar[rrdd]^{\widetilde\kappa} &&
 \bullet \ar[rr]_{p_2'} \ar[dd]^{\kappa_1^{}} && \W \ar[dd]_{\kappa} \ar@/^1.5pc/[dddd]^{\wid{f}} \\
 \\
\bullet \ar[rr]_{\kappa_2^{}} \ar[dd]^{p_1'} && \X'' \ar[rr]_{p_2^{}} \ar[dd]^{p_1^{}} && \X \ar[dd]_{f} \\
\\
\W \ar@/_1.5pc/[rrrr]_{\wid{f}} \ar[rr]^{\kappa} && \X \ar[rr]^{f} && \Y 
}
\qquad \qquad
\xymatrix{
\W \ar[rr]^{\kappa} \ar[dd]_{\delta} & & \X \ar[dd]^{\Delta} \\
\\
\W'' \ar[rr]_{\widetilde\kappa} & & \X''
}
\]

In what follows, let
\stepcounter{thm}
\begin{equation*}\label{iso:kap-wid}\tag{\thethm}
\kappa^*\ush{f} \iso \ush{\wid{f}}
\end{equation*}
be the composite
$\bL\kappa^*\ush{f}\xrightarrow[\eqref{eq:gm}]{\Iso} \ush{\kappa}\ush{f} \iso \ush{\wid{f}}$, and
let $\kappa^*\omega_f[r] \iso \omega_{\wid{f}}[r]$ be the one induced by the canonical isomorphsm
$\kappa^*\omega_f \iso \omega_{\wid{f}}$.
\begin{thm}\label{thm:verd-complete} The following diagram commutes
\[
\xymatrix{
\bL\kappa^*\omega_f[r] \ar[r]^\Iso
\ar[d]^{\rotatebox{90}{\makebox[0.1cm]{\Iso}}}_{\bL\kappa^*{{\bf v}_{{}_f}}}
& \kappa^*\omega_f[r] \ar@{=}[r]
& \omega_{\wid{f}}[r] 
\ar[d]_{\,\rotatebox{-90}{\makebox[-0.1cm]{\Iso}}}^{{\bf v}_{{}_{\wid{f}}}} \\
\bL\kappa^*(\ush{f}\co_\Y) \ar[rr]^{\Iso}_{\eqref{iso:kap-wid}} && \ush{\wid{f}}\co_\Y
}
\]
\end{thm}

\proof 
Since ${\bf v}_{{}_f}$ and ${\bf v}_{{}_{\wid{f}}}$ are isomorphisms, we assume 
$\ush{f}\co_\Y$ and $\ush{\wid{f}}\co_\Y$ are complexes which are zero in all degrees
except at the $(-r)$-th spot, where each is locally free (in fact invertible). This means
we write $h^*(\ush{f}\co_\Y)=\bL h^*(\ush{f}\co_\Y)$ 
(resp.~$h^*(\ush{\wid{f}}\co_\Y)=\bL h^*(\ush{\wid{f}}\co_\Y)$)
for any map of schemes to $\X$ (resp.~$\W$). Similarly $\bL h^*\omega_f[d]
= h^*\omega_f[d]$ etc. Let
\[
\phi\colon \kappa^*(\wnor\De) \iso \wnor\delta
\]
be the canonical isomorphism. 
We have to show that the diagram $\clubsuit$ below
commutes.
\[
\xymatrix{
\ush{\wid{f}\,}\co_\Y\otimes\wnor\delta[-r] 
\ar[dd]^{\,\rotatebox{-90}{\makebox[-0.1cm]{\Iso}}}_{\eqref{iso:pre-verdier}}
\ar@{}[ddrr]|{\clubsuit}
& \kappa^*\ush{f}\co_\Y\otimes\wnor\delta[-r] \ar[l]_-{\Iso}^{\eqref{iso:kap-wid}} 
& \kappa^*(\ush{f}\co_\Y\otimes\wnor\De[-r]) \ar[l]_-{\Iso}^{{\text{via $\phi$}}}
\ar[dd]_{\,\rotatebox{-90}{\makebox[-0.1cm]{\Iso}}}^{{\text{via \eqref{iso:pre-verdier}}}}\\
&& \\
\co_\W \ar@{=}[rr] && \kappa^*\co_\X
}
\]
We expand $\clubsuit$ as follows (with the label (C.4.2), referring to the label in \cite{fub-abs}):
\[
\xymatrix{
\ush{\wid{f}\,}\co_\Y\otimes\wnor\delta[-r] \ar@{=}[d] 
& \kappa^*\ush{f}\co_\Y\otimes\wnor\delta[-r] \ar[l]_-{\Iso}^{\eqref{iso:kap-wid}} \ar@{=}[d]
\ar@{}[ddr]|{\square_2}
& \kappa^*(\ush{f}\co_\Y\otimes\wnor\De[-r]) \ar[l]_-{\Iso}^{{\text{via $\phi$}}}
\ar@{=}[dd] \\ 
\delta^\btrg(\pi_{{}_2}^*\ush{\wid{f}}\co_\Y) 
\ar@{}[ddr]|{\square_1}
\ar[dd]^{\,\rotatebox{-90}{\makebox[-0.1cm]{\Iso}}}_{\eqref{iso:bc-sharp}}
& \delta^\btrg(\pi_{{}_2}^*\kappa^*\ush{f}\co_\Y) \ar[l]_-{\Iso}^{\eqref{iso:kap-wid}} \ar@{=}[d]
& \\ 
& \delta^\btrg(\widetilde\kappa^*p_{{}_2}^*\ush{f}\co_\Y)
\ar[d]^{\,\rotatebox{-90}{\makebox[-0.1cm]{\Iso}}}_{\eqref{iso:bc-sharp}}
& \kappa^*\De^\btrg(p_{{}_2}^*\ush{f}\co_\Y)
\ar[d]_{\,\rotatebox{-90}{\makebox[-0.1cm]{\Iso}}}^{\eqref{iso:bc-sharp}}
\ar[l]_\Iso^{(C.4.2)} \\ 
\delta^\btrg(\ush{\pi_{{}_1}}\co_\W) \ar[d]^{\,\rotatebox{-90}{\makebox[-0.1cm]{\Iso}}}_{\eta'_\delta}
& \delta^\btrg(\widetilde\kappa^*\ush{p_{{}_1}}\co_\X)
\ar@{}[dr]|{\square_3}
\ar[d]^{\,\rotatebox{-90}{\makebox[-0.1cm]{\Iso}}}_{\eta'_\delta} 
\ar[l]_-{\Iso}^{\alpha_1}
& \kappa^*\De^\btrg(\ush{p_{{}_1}}\co_\X)
\ar[l]_\Iso^{(C.4.2)}
\ar[d]_{\,\rotatebox{-90}{\makebox[-0.1cm]{\Iso}}}^{\eta'_\De} \\ 
\ush{\delta}(\ush{\pi_{{}_1}}\co_\W) \ar[d]^{\,\rotatebox{-90}{\makebox[-0.1cm]{\Iso}}}
\ar@{}[drr]|{\square_4}
& \ush{\delta}(\widetilde\kappa^*\ush{p_{{}_1}}\co_\X) \ar[l]_-{\Iso}^{\alpha_2}
& \kappa^*\ush{\De}(\ush{p_{{}_1}}\co_\X)
\ar[l]_\Iso^{\eqref{iso:bc-sharp}}
\ar[d]^{\,\rotatebox{-90}{\makebox[-0.1cm]{\Iso}}} \\ 
\co_\W \ar@{=}[rr] & & \kappa^*\co_\X
}
\]
The maps $\eta'_{\De}$ and $\eta'_{\delta}$
 are the maps defined in \eqref{iso:eta'-i}.
The maps $\alpha_i^{-1}$ are induced by the isomorphism 
$\ush{\pi_{{}_1}}\kappa^*\co_\X \iso \widetilde\kappa^*\ush{p_{{}_1}}\co_\X$ resulting from the following
composite of natural maps
\[
\ush{\pi_{{}_1}}\kappa^* \iso \ush{\kappa_{{}_1}'}\ush{p_{{}_1}'}\kappa^* \xrightarrow[\eqref{eq:gm}]{\Iso}
\ush{\kappa_{{}_1}'}\ush{p_{{}_1}'}\ush{\kappa} \iso \ush{\widetilde\kappa}\ush{p_{{}_1}} 
\xrightarrow[\eqref{eq:gm}]{\Iso} \widetilde\kappa^*\ush{p_{{}_1}}.
\] 
In the above expansion of $\clubsuit$, the unlabelled sub-rectangles clearly commute. 
Sub-rectangle $\square_2$ commutes by 
definition of the isomorphism \cite[(C.4.2)]{fub-abs}. Prop.\,C.4.3 of \cite{fub-abs}
gives the commutativity
of $\square_3$. For $\square_4$ we apply the outer border of the following diagram
on~$\co_\X$.
\[
\xymatrix{
\ush{\delta}\ush{\pi_{{}_1}}\kappa^* \ar[dd] \ar[r] & 
\ush{\delta}\ush{\kappa_{{}_1}'} \ush{p_{{}_1}'}\kappa^* \ar[dd] \ar[r]^{\kappa^* \cong \ush{\kappa}}  &  
\ush{\delta}\ush{\kappa_{{}_1}'} \ush{p_{{}_1}'}\ush{\kappa} \ar[r] \ar[dd] & 
\ush{\delta}\ush{\widetilde\kappa}\ush{p_{{}_1}} \ar[d]^{\qquad \scriptscriptstyle\square\square_4} 
\ar[r]^{\ush{\widetilde\kappa} \cong \widetilde\kappa^*} &  
\ush{\delta}\widetilde\kappa^*\ush{p_{{}_1}} \ar[d]^{\eqref{iso:bc-sharp}} \\
&  &  & \ush{\kappa}\ush{\Delta}\ush{p_{{}_1}} 
\ar[d] \ar[r]_{\ush{\kappa} \cong \kappa^*} & \kappa^*\ush{\Delta}\ush{p_{{}_1}} \ar[d] \\
\kappa^* \ar@{=}[r] &  \kappa^* \ar[r]_{\kappa^* \cong \ush{\kappa}} &  \ush{\kappa} \ar@{=}[r] &
\ush{\kappa} \ar[r]_{\ush{\kappa} \cong \kappa^*} &  \kappa^* 
}
\]
The unlabelled arrows are the obvious ones. The rectangle $\scriptstyle\square\square_4$ commutes 
by \cite[Lemma\,A.1.4]{fub-abs}
while the remaining commute for obvious reasons.

It remains to prove that $\square_1$ commutes. To that end, it suffices to prove that the outer border 
of the
following diagram commutes where the unlabelled arrows are the obvious ones coming from 
pseudofunctoriality 
of $(-)^*$ or $\ush{(-)}$, the ones labelled b-ch are induced by suitable base-change isomorphisms as 
given 
in \eqref{iso:bc-sharp}, the ones labelled $\ush{} = {}^*$ are induced by \eqref{eq:gm} and the ones 
labelled 
$\gamma_i$ are induced by the composite
$\ush{\kappa_{{}_2}'}\ush{\kappa_{{}_1}} \iso \ush{\widetilde\kappa} 
\xrightarrow[\eqref{eq:gm}]{\Iso} \widetilde\kappa^*$.
\[
\xymatrix{
\pi_{{}_2}^*\ush{\wid{f}} \ar[dddd]_{\text{b-ch}} \ar[rr] & & \kappa_{{}_2}'^* p_{{}_2}'^*\ush{\kappa}\ush{f} 
\ar[d]_{\text{b-ch}}^{\hspace{6.7em}\blacktriangle} \ar[rr]^{\ush{} = {}^*} & 
& \kappa_{{}_2}'^* p_{{}_2}'^*\kappa^*\ush{f} \ar[d]  \\ 
& & \kappa_{{}_2}'^*\ush{\kappa_{{}_1}} p_{{}_2}^*\ush{f} \ar[dd]_{\text{b-ch}} 
\ar[rr]^{\ush{} = {}^*} \ar[rd]^{\ush{} = {}^*} & &  \kappa_{{}_2}'^*\kappa_{{}_1}^* p_{{}_2}^*\ush{f} 
\ar[d]_{\triangle\hspace{3.5em}} \\
& \ddag & & \ush{\kappa_{{}_2}'}\ush{\kappa_{{}_1}} p_{{}_2}^*\ush{f} \ar[d]_{\text{b-ch}} \ar[r]_{\gamma_1}
 & \widetilde\kappa^*p_{{}_2}^*\ush{f} \ar[d]^{\text{b-ch}} \\
& \ush{\kappa_{{}_1}'}\kappa_{{}_2}^* \ush{p_{{}_1}}f^* \ar[d]_{\text{b-ch}}^{\hspace{2.5em}\blacktriangle} 
\ar[r]^{\text{b-ch}} \ar[rrd]^{\ush{} = {}^*} & \kappa_{{}_2}'^*\ush{\kappa_{{}_1}} \ush{p_{{}_1}}f^* \ar[r]^{\ush{} = {}^*}
& \ush{\kappa_{{}_2}'}\ush{\kappa_{{}_1}}\ush{p_{{}_1}}f^*  \ar[r]_{\gamma_2} \ar[d]_{\blacktriangle\hspace{2em}}
& \widetilde\kappa^*\ush{p_{{}_1}}f^* \\
\ush{\pi_{{}_1}}\wid{f}^* \ar[r] & \ush{\kappa_{{}_1}'}\ush{p_{{}_1}'}\kappa^*f^* \ar[r] ^{\ush{} = {}^*}
& \ush{\kappa_{{}_1}'}\ush{p_{{}_1}'}\ush{\kappa}f^* \ar[r] & \ush{\kappa_{{}_1}'}\ush{\kappa_{{}_2}} \ush{p_{{}_1}}f^*
}
\]
Now $\ddag$ commutes by transitivity of the base-change isomorphism, 
see \cite[Prop.\,A.1.1]{fub-abs}. 
For the diagrams labelled
$\blacktriangle$, we refer to \cite[Lemma\,A.1.4]{fub-abs}, while $\triangle$ commutes because 
of the pseudofunctorial nature of the isomorphism $\ush{(-)} \cong (-)^*$
of~\eqref{eq:gm} over the category of formally \'etale maps. 
The unlabelled diagrams commute for trivial reasons.
\qed

\setcounter{subsubsection}{\value{thm}} 
\subsubsection{} \stepcounter{thm}
There is a related result. Suppose $f\colon \X\to \Y$ is a smooth map of relative dimension $d$ in
$\bbG$ and suppose $f$ factors as
\[
{\xymatrix{
\X \ar@/_2pc/[rr]_f \ar[r]^{\wid{f}} &\wid{\Y} \ar[r]^\kappa & \Y \\
}}
\]
where $\kappa$ is the completion of $\Y$ along a coherent $\co_\Y$-ideal $\I$, 
and $\wid{f}$ is smooth (necessarily of relative dimension $d$). Note that
\[\wid{\Y}\times_\Y\wid{\Y}=\wid{\Y}.\]

Consider the commutative diagram of cartesian squares:
\[
{\xymatrix{
\X \ar[d]_{\wid{f}} \ar@{=}[r] \ar@{}[dr]|{\square} & \X \ar[d]^{\wid{f}} \ar@/^2pc/[dd]^{f} \\
\wid{\Y} \ar@{=}[r] \ar@{=}[d] \ar@{}[dr]|{\square} & \wid{\Y} \ar[d]^{\kappa} \\
\wid{\Y} \ar[r]_{\kappa} & \Y
}}
\]
From \eqref{iso:bc-sharp} we conclude that we have an isomorphism
\[\ush{f} \iso {\bf 1}_\X^*\ush{f} \iso \ush{\wid{f}}\kappa^*. \leqno{(*)}\]

Now clearly, $\omega_f=\omega_{\wid{f}}$. Call the common $\co_\X$-module $\omega$.
We have two related isomorphisms, namely,
${\bar{\bf v}}_{\<\<{}_{\wid{f}}}\colon \omega[d] \iso \ush{\wid{f}}\co_{\wid{\Y}}$
and 
${\bar{\bf v}}_{\<\<{}_f}\colon \omega[d] \iso \ush{f}\co_\Y$.
With these notations, we have the following Proposition, related to \Tref{thm:verd-complete}:

\begin{prop}\label{prop:verd-complete2} With notations as above,
the following diagram commutes:
\[
{\xymatrix{
\omega[d] \ar[d]_{{\bar{\bf v}}_{\<\<{}_f}}^{\rotatebox{90}{\makebox[0.1cm]{\Iso}}}
 \ar[rr]_-{{\bar{\bf v}}_{\<\<{}_{\wid{f}}}}^-{\Iso} &&
\ush{\wid{f}}\co_{\wid{\Y}} \ar@{=}[d] \\
\ush{f}\co_\Y \ar[rr]_-{(*)}^-{\Iso} && \ush{\wid{f}}\kappa^*\co_\Y .
}}
\]
\end{prop}

\proof
We have the following commutative diagram with all squares cartesian:
\[
{\xymatrix{
\X'' \ar[d]_{p_{{}_1}}  \ar[r]^{p_{{}_2}} \ar@{}[dr]|{\square}
& \X \ar[d]_{\wid{f}} \ar@{=}[r] \ar@{}[dr]|{\square} & \X \ar[d]^{\wid{f}} \ar@/^2pc/[dd]^{f} \\
\X \ar[r]_{\wid{f}} \ar@{=}[d] \ar@{}[dr]|{\square} &
\wid{\Y} \ar@{=}[r] \ar@{=}[d] \ar@{}[dr]|{\square} & \wid{\Y} \ar[d]^{\kappa} \\
\X \ar[r]_{\wid{f}} & \wid{\Y} \ar[r]_{\kappa} & \Y
}}
\]
We claim that diagram $(**)$ below commutes:
\[
{\xymatrix{
p_{{}_2}^*\ush{f} \ar[rr]_{\eqref{iso:bc-sharp}}^{\Iso} 
\ar[d]_{(*)}^{\rotatebox{90}{\makebox[0.1cm]{\Iso}}} && \ush{p_{{}_1}}f^* \ar@{=}[d]\\
p_{{}_2}^*\ush{\wid{f}}\kappa^* \ar[rr]_{\eqref{iso:bc-sharp}}^{\Iso}  && 
\ush{p_{{}_1}}\wid{f}^*\kappa^*
}}
\leqno(**) \]
Indeed, this follows immediately from the horizontal transitivity of the base-change isomorphism 
(see \cite[Prop.\,A.1.1\,(i)]{fub-abs})
corresponding to the ``composite" of base-change diagrams:
\[
{\xymatrix{
\X'' \ar[d]_{p_{{}_1}}  \ar[r]^{p_{{}_2}} \ar@{}[dr]|{\square}
& \X \ar[d]_{\wid{f}} \ar@{=}[r] \ar@{}[dr]|{\square} & \X \ar[d]^f \\
\X \ar[r]_{\wid{f}} & \wid{\Y} \ar[r]_\kappa & \Y
}}
\]
We therefore have the following commutative diagram,
where the square on the left in induced by $(**)$.
\[
\xymatrix{
\De^{\<*}(p_{{}_2}^*\ush{f}\co_\Y\!)\otimes\omega^{-1}[-d] \ar[r]_{\eqref{iso:bc-sharp}}^{\Iso} 
\ar[d]_{(*)}^{\rotatebox{90}{\makebox[0.1cm]{\Iso}}} & 
\De^{\<*}(\ush{p_{{}_1}}\co_\X\!)\otimes\omega^{-1}[-d] \ar@{=}[d]
\ar[r]^-{\Iso} & \ush{\De}(\ush{p_{{}_1}}\co_\X\!) \ar@{=}[d] \ar[r]^-{\Iso} & \co_{\<\X} \ar@{=}[d] \\
\De^{\<*}(p_{{}_2}^*\ush{\wid{f}}\co_{\wid{\Y}}\!)\otimes\omega^{-1}[-d] \ar[r]_{\eqref{iso:bc-sharp}}^{\Iso}  &
\De^{\<*}(\ush{p_{{}_1}}\wid{f}^*\co_{\wid{\Y}}\!)\otimes\omega^{-1}[-d] \ar[r]^-{\Iso} 
& \ush{\De}(\ush{p_{{}_1}}\co_\X\!) \ar[r]^-{\Iso} & \co_{\<\X}
}
\]
In other words
\[
{\xymatrix{
\ush{f}\co_\Y\otimes\omega^{-1}[-d] \ar[d]_{(*)}^{\rotatebox{90}{\makebox[0.1cm]{\Iso}}} 
\ar[r]^-{\Iso} & \co_\X \ar@{=}[d] \\
\ush{\wid{f}}\co_{\wid{\Y}}\otimes\omega^{-1}[-d]
\ar[r]^-{\Iso} & \co_\X 
}}
\]
commutes. This is equivalent to the statement of the Proposition.
\qed

\subsection{Base change and Verdier's isomorphism} As mentioned earlier, according to
\cite[p.740,\,Theorem\,2.3.5\,(a)]{cm}, for any Cohen-Macaulay map between ordinary schemes $f\colon X\to Y$, and
any base change $u\colon Y'\to Y$, with $X'=X\times_YY'$, 
$f'\colon X'\to Y'$ and $v\colon X'\to X$ the base change maps, there is a natural isomorphism
$\theta_u^f\colon v^*\omgs{f}\iso \omgs{f'}$. In the event  $f$ is smooth, then using Verdier's
isomorphisms for $f$ and $f'$ to identify $\omgs{f}$ with $\omega_f$ and $\omgs{f'}$ with
$\omega_{f'}$, the map $\theta_u^f$ corresponds to the obvious canonical map 
(see [Ibid.,  p.740,\,Theorem\,2.3.5\,(b)]).
The difficulty in transferring this statement
to formal schemes is that defining $\theta_u^f$ required certain special local compactifications of $f$
which may or may not be available for general formal scheme maps. However these
difficulties disappear if the map $u$ is flat, and the proof in loc.cit.\:works mutatis mutandis. The
precise statement is:
\begin{thm}\label{thm:verd-flat-bc} Suppose
\[
{\xymatrix{
\X' \ar[d]_{f'} \ar[r]^v \ar@{}[dr]|\square & \X \ar[d]^f \\
\Y' \ar[r]_u & \Y
}}
\]
is a cartesian square with $f$ smooth, in $\bbG$, of relative dimension $d$, and $u$ flat.
Let $\theta\colon v^*\ush{f}\co_\Y \iso \ush{(f')}u^*\co_\Y=\ush{(f')}\co_{\Y'}$ 
be the resulting base change isomorphism (see \eqref{iso:bc-sharp}). Then
the isomorphism ${\bf v}_{f'}^{-1}\smcirc\theta\smcirc v^*({\bf v}_f)\colon 
v^*\omega_f[d] \iso \omega_{f'}[d]$
is the obvious canonical map.
\end{thm} 

\section{\bf Residues}

\subsection{Verdier residue} Let $f \colon X \to Y$ be smooth of relative dimension $r$ and
let $Z \hookrightarrow X$ be a closed subscheme proper over $Y$. Analogous to to the
abstract residue $\ares{Z}$ in \cite[(5.2.2)]{fub-abs} one has
the {\em Verdier residue along $Z$}
\stepcounter{thm}
\begin{equation*}\label{def:vres}\tag{\thethm}
\res{Z}\colon \Rr^rf_*\omega_f \to \co_Y
\end{equation*}
defined as the composite
\[\Rr^r_Zf_*\omega_f \xrightarrow[{\eqref{iso:loc-coh}}]{\Iso} 
\Rp\X^r {\wid{f}}_*\kappa^*\omega_f \iso \Rp\X^r{\wid{f}}_*\omega_{\wid{f}}
 \xrightarrow{\vin{\wid{f}}} \co_Y.\]
where the middle isomorphism is induced by the canonical one $\kappa^*\omega_f \iso \omega_{\wid{f}}$. 
By compatibility of the Verdier isomorphism with completions (\Tref{thm:verd-complete}), the following diagram commutes (where, as before, $\ares{Z}$ is the abstract residue map defined in 
\cite[(5.2.2)]{fub-abs}): 
\stepcounter{thm}
\begin{equation*}\label{diag:res-int-2}\tag{\thethm}
\xymatrix{
\Rr^r_Zf_*\omega_f \ar[r]^-{\Iso}_-{\bar{\bf v}} \ar@/_1.0pc/[dr]_{\res{Z}} & 
\Rr^r_Z{f}_*\kappa^*\omgs{f} \ar[d]^{\ares{Z}} \\
 & \co_Y
}
\end{equation*}

\begin{rem} While we have defined residues in general, our interest is really in the case 
where $\Rp{\X}^j(\eF)=0$ for every 
$j>r$ and every $\eF\in\Avc(\X)$, for then $(\omega_f, \vin{f})$ represents
the functor $\Hom_Y({\Rp{\X}^r(\boldsymbol{-}}),\, \co_Y)$
on coherent $\co_\X$-modules 
(cf.\,\cite[Cor.\,5.1.4]{fub-abs}). Even here
the most useful situation is when $Y=\Spec{\,A}$ and $\X =\Spf{\,R}$ where $R$ is
an adic ring, with a defining ideal $I$ generated by $r$ elements, with $R/I$ finite and
flat over $A$.
\end{rem}

The various relationships between the abstract residue, Verdier residue, the trace
and the Verdier intergal are captured in the following commutative commutative diagram:
\stepcounter{thm}
\[
\begin{aligned}\label{diag:res-int-3}
\xymatrix{
\Rr^r_Zf_*\omega_f \ar[rr]^\Iso_{\eqref{iso:loc-coh}} \ar[d]^{\bf v} \ar@/_3.5pc/[dd]_{\res{Z}} 
& &\Rp\X^r{\wid{f}}_*\omega_{\wid{f}} \ar[d]_{\bf v}
\ar@/^3.5pc/[dd]^{\vin{\wid{f}}}\\
\Rr^r_Zf_*\omgs{f} \ar[rr]^\Iso_{\eqref{iso:iGZ-iGpX}} \ar[d]^{\ares{Z}} && \Rp\X^r{\wid{f}}_*\omgs{\wid{f}}
\ar[d]_{\tin{\wid{f}}} \\
\co_Y \ar@{=}[rr] && \co_Y
}
\end{aligned}\tag{\thethm}
\]

In the event $f\colon X\to Y$ is {\em proper} we have the following commutative diagram
\stepcounter{thm}
\[
\begin{aligned}\label{diag:res-int-4}
\xymatrix{
\Rr^r_Zf_*\omega_f  \ar[rr] \ar[d]^{\bf v} \ar@/_3.5pc/[dd]_{\res{Z}}
 && \Rr^rf_*\omega_f \ar[d]_{\bf v} \ar@/^3.5pc/[dd]^{\vin{f}}
 \\
\Rr^r_Zf_*\omgs{f}  \ar[d]^{\ares{Z}} \ar[rr]
& &\Rr^rf_*\omgs{f}  \ar[d]_{\tin{f}}\\
\co_Y \ar@{=}[rr] && \co_Y
}
\end{aligned}\tag{\thethm}
\]

\subsection{Some residue formulas}

%
Suppose $A\to R$ is a finite type map of rings which is
 {\em smooth}. Set $R''=R\otimes_AR$. As before, the two
$R$-algebra structures on $R''$ will be denoted $R_1$ and $R_2$, with $R_k$
denoting the algebra corresponding to the
projection $p_{{}_k}\colon X''\set X\times_YX \to X$ for $k\in\{1, 2\}$.
The diagonal map $\De\colon X' \hookrightarrow X''$ corresponds to the surjective
map $R''\to R$ given by $t_1\otimes t_2\mapsto t_1t_2$. Suppose the kernel of this map,
i.e., the ideal of the diagonal immersion, is generated by $r$-elements $\{s_1,\dots, s_r\}$.
Since $R$ is smooth over $A$ of relative dimension $r$, the sequence ${\bf s}=(s_1,\dots, s_r)$
is necessarily a $R''$-sequence. By part\,(b) of 
\Pref{prop:ds/s-1} we get the following formula, which is at the heart of much of what we do
in this paper.
\stepcounter{thm}
\begin{equation*}\label{eq:ds/s}\tag{\thethm}
\res{\De,\,p_{{}_1}}\begin{bmatrix}
\wdd{s_1}{s_r}\\ s_1,\,\dots,\,s_r
\end{bmatrix} = 1.
\end{equation*}

\begin{prop}\label{prop:thom-class} Let $X=\Spec{\,A}$ and $Y=\Spec{\,R}$ be affine
schemes, and suppose $f\colon X\to Y$ is a smooth map of relative dimension $r$.
Suppose further that we have a closed subscheme $Z$ of $X$ such that $Z\to Y$ is an
isomorphism and the ideal $J$ of $R$ giving the closed subscheme $Z$ of $X$ is generated by 
$r$-elements
$\{t_1,\,\dots,\,t_r\}$ of $R$. Then
\stepcounter{thm}
\begin{equation*}\label{eq:dt/t}\tag{\thethm}
\res{Z}\begin{bmatrix}
\wdd{t_1}{t_r}\\ t_1,\,\dots,\,t_r
\end{bmatrix} = 1.
\end{equation*}
\end{prop}
\proof First note that since $f$ is smooth, ${\bf t}=(t_1,\,\dots,\,t_r)$ is an $R$-regular
sequence. Next note that the question is local on $Y$ and so we may assume, without
loss of generality, that the diagonal immersion $\De\colon X\to X''$ is cut out by $r$-elements
$\{s_1,\,\dots,\,s_r\}$ in $R''=R\otimes_AR$. As in \Ssref{ss:verd-loc}, we write
$I$ for the ideal of the diagonal and use the notations of that subsection. Let $Z=\Spec{\,B}$.
Let $\sigma\colon Y\to X$ be the section defined by $Z$, and $i\colon Z\hookrightarrow X$
the natural closed immersion. We have a commutative diagram
with all sub-rectangles cartesian.
\[
\xymatrix{
\underset{}{Z} \ar@{^(->}[d]_i \ar[r]^{\sigma_{\<\<\<{}_Z}} \ar@{}[dr]|{\square} 
& \underset{}{X} \ar@{^(->}[d]^{\De}\\
X \ar[r]^{\sigma_{\<\<\<{}_X}} \ar[d]_f \ar@{}[dr]|{\square}
& X'' \ar[d]^{p_{{}_1}} \\
Y \ar[r]_{\sigma} & X
}
\]
We now need some results from \cite{cm} regarding \emph{non-flat} base-change. Since $\sigma$ is a closed
immersion, the usual flat-base-change results do not apply. Nevertheless, we do have the following.
First, there is a base change isomorphism
$\theta=\theta_\sigma^f\colon {\sigma_{\<\<\<{}_X}}^*\omgs{f} \iso \omgs{p_{{}_1}}$
as in [Ibid.,  p.740,\,Theorem\,2.3.5\,(a)]. 
Next, by [Ibid., Prop.\,6.2.2, pp.755--756],
under this isomorphism, residues are compatible. In other words, the diagram
\[
\xymatrix{
\sigma^*\Rr^r_\De {p_{{}_1}}_*(\omgs{p_{{}_1}}) \ar[r]^{\Iso} \ar[d]_{\sigma^*{\ares{\De}}}
& \Rr^r_Zf_*({\sigma_{\<\<\<{}_X}}^*\omgs{p_{{}_1}}) 
\ar[d]_{\,\rotatebox{-90}{\makebox[-0.1cm]{\Iso}}}^{\theta} \\
\co_\Y & \Rr^r_Zf_*\omgs{f} \ar[l]_{\ares{Z}}\\
}
\]
commutes. Finally on replacing $\omgs{f}$ by $\omega_f$ and $\omgs{p_{{}_1}}$ by 
$\omega_{p_{{}_1}}$ via ${\bf v}_f$ and ${\bf v}_{p_{{}_1}}$, according to 
[Ibid.,  p.740,\,Theorem\,2.3.5\,(b)], the map $\theta$ reduces to the standard
identity ${\sigma_{\<\<\<{}_X}}^*\omega_{p_{{}_1}}= \omega_f$. 
Thus it follows that
if $u_i$, $i=1,\,\dots,\,r$, are the images of $s_i$ in $R$ under the map
$R''\to R$ corresponding to $\sigma_{\<\<\<{}_X}\colon X\to X''$, (so that $J$ is generated
by the set $\{u_1,\,\dots,\,u_r\}$) we have (via \eqref{eq:ds/s})
\[\res{Z}\begin{bmatrix}
\wdd{u_1}{u_r}\\ u_1, \dots, u_r
\end{bmatrix} =1.\]
Since
$\bigl [\begin{smallmatrix} \wdd{u_1}{u_r}\\ 
u_1, \dots, u_r \end{smallmatrix}\bigr ] =
\bigl [\begin{smallmatrix} \wdd{t_1}{t_r}\\
t_1, \dots, t_r \end{smallmatrix}\bigr ]$, hence the result.
\qed

\section{\bf Residues along sections}

Let $f\colon X\to Y$ be a smooth map of relative dimension $r$,
which is separated. 
We begin with some  notations and conventions. In general, if we are working over affine schemes
(ordinary or formal) we will use the same notations for maps between modules as the corresponding
sheaves. For example if $A\to R$ is smooth map of rings of relative dimension $r$, and
$I$ an $R$ ideal generated by a regular sequence $\{t_1,\dots,t_r\}$ such that $A\to B\set R/I$
is finite, then with $X=\Spec{\,R}$, $Y=\Spec{\,A}$ and $Z=\Spec{\,B}$, and $f\colon X\to Y$
the map given by $A\to R$, we will write
$\res{Z}\colon \Hr^r_Z(X,\,\omega_f)\to A$ instead of $\Gamma(Y,\,\res{Z})$. As another illustration
of this principle, in the above situation,
if $\omega_{R/A}$ is the $A$-module given by $\omega_{R/A}=\Gamma(X,\,\omega_f)$,
then we will make no distinction between $\Hr^r_Z(X,\,\omega_f)$ and
$\Hr^r_I(\omega_{R/A})$. 

\subsection{}\label{ss:affine-section} Suppose $Y=\Spec{\,A}$ and $Z\hookrightarrow X$ is a closed subscheme such that $Z\to Y$ is an
isomorphism and $Z$ lies in an open affine subscheme $U=\Spec{\,R}$ of $X$ such that 
$Z$ is given in $U$ by an ideal $I$ which is generated by $r$ elements $t_1, \dots, t_r $
of $R$. We have a map
\[\res{{\bf t}}\colon \Hr^r_Z(X,\,\omega_f)\to A\]
defined by the formula
\stepcounter{thm}
\begin{equation*}\label{eq:st-res}\tag{\thethm}
\res{{\bf t}}\begin{bmatrix}
\wdd{t_1}{t_r}\\
t_1^{\alpha_1},\dots, t_r^{\alpha_r} 
\end{bmatrix}
= \begin{cases}
1 & \text{when $\alpha_1=\dots = \alpha_r=1$} \\
0 & \text{otherwise.}
\end{cases}
\end{equation*}
This map depends \textit{a priori} on the choice of ${\bf t}=(t_1, \dots, t_r)$, but as we will see later,
it is independent of it. It should be pointed out that if $Z$ is also defined (in $U$) by the vanishing
of $s_1,\,\dots,\,s_r$, then by \cite[Thm.\,C.7.2\,(iii)]{fub-abs}
\stepcounter{thm}
\begin{equation*}\label{eq:dt/t=ds/s}\tag{\thethm}
\begin{bmatrix}
\wdd{t_1}{t_r}\\
t_1, \dots, t_r
\end{bmatrix}
=
\begin{bmatrix}
\wdd{s_1}{s_r}\\
s_1, \dots, s_r
\end{bmatrix}.
\end{equation*}
Moreover, there is an $A$-module direct sum decomposition 
\[\Hr^r_Z(X,\,\omega_f)= \Hr^r_I(\omega_{R/A})=
\bigoplus_{\underline{\alpha}} A\bigl [\begin{smallmatrix}
\wdd{t_1}{t_r}\\
t_1^{\alpha_1},\dots, t_r^{\alpha_r} 
\end{smallmatrix}
\bigr ]
\]
with ${\underline{\alpha}}=(\alpha_1,\dots,\alpha_r)$ running over $r$-tuples of
positive integers. The summands are a free $A$-modules.
While this decomposition depends on ${\bf t}$, the summand generated
by $\bigl [\begin{smallmatrix}
\wdd{t_1}{t_r}\\
t_1,\dots, t_r
\end{smallmatrix}
\bigr ]$
is independent of ${\bf t}$ by \eqref{eq:dt/t=ds/s}. In what follows, let
\[ 
\theta_{{}_Z}= \begin{bmatrix}
\wdd{t_1}{t_r}\\
t_1, \dots, t_r
\end{bmatrix}.
\]

\subsection{Relative projective space} Let ${\bbP}={\bbP}^r_Y$, 
the relative projective space of relative dimension $r$
over an ordinary scheme $Y$. We regard $\bbP={\bf Proj}(\co_Y[T_0,\dots,T_r])$. Let
$\pi\colon \bbP\to Y$ be the structure map and
\[\int_{\bbP/Y}\colon \Rr^r\pi_*\omega_\pi \iso \co_Y\]
be the standard trace map (known to be an isomorphism) defined, for example in 
\cite[${\rm{III}}_1$, 2.1.12]{ega} or \cite[p.152,\,Theorem 3.4]{RD}.
The generating section $\mu=\mu_{{}_\bbP}$ of
$\Rr^r\pi_*\omega_\pi$ corresponding to the standard section $1$
of $\co_Y$ is described as follows.  Let 
${\mathscr U}=\{U_i\mid i=0,\dots, r\}$ be the open cover of $\bbP$
given by $U_i=\{T_i\ne 0\}$. On $U_0\cap\ldots\cap U_r$
we have inhomogeneous coordinates $t_i=T_i/T_0$, $i=1,\dots,r$
whence a section
$$
\check{\mu}_T\set \frac{{\mathrm{d}}t_1\wedge\ldots\wedge{\mathrm{d}}
t_r}{t_1\dots t_r}
\in \Gamma(U_0\cap\dots\cap U_r, \omega_\pi).
$$
We have an isomorphism
$$
\Hr^r(\pi_*\check{\mathcal C}^\bullet({\mathscr U},\,\,\omega_\pi))  \iso \Rr^r\pi_*\omega_\pi
$$
and $\check{\mu}_T$ has a natural image in the left side as a
\v{C}ech cohomology class. Let $\mu$ be the corresponding
element on the right side.  The section $\mu$ does not depend 
on the choice of homogeneous coordinates $T_0,\dots,T_r$ of $\bbP$ 
(cf. \cite[p.34,\,Lemma\,2.3.1]{conrad}) and is the sought after
section.

Let $Z_0$ be the closed subscheme of $\bbP$ defined
by $\{T_i=0\mid i=1,\dots, r\}$, i.e., the intersection of the relative hyperplanes
$H_i=\{T_i=0\}$, $i=1,\dots, r$. Then $Z_0\to Y$ is an isomorphism. 
The section $\sigma_0\colon Y\to \bbP$
defined by $Z_0$ is the $Y$-valued point of the $Y$-scheme $\bbP$ given by the 
``homogeneous co-ordinates" $(1, 0, 0,\dots, 0)$. 

Now suppose $Y=\Spec{\,A}$.
It is well known 
(see for example \cite[p.74,\,Prop.\,(8.4)]{ast117}, the proof of which generalizes to our situation)
that the following diagram commutes.
\stepcounter{thm}
\begin{equation*}\label{diag:res-thm-p}\tag{\thethm}
\xymatrix{
\Hr^r_{Z_0}(\bbP,\,\omega_\pi) \ar@/_1.0pc/[dr]_{\res{{\bf t}}} \ar[r] & 
\Hr^r(\bbP,\,\omega_\pi) \ar[d]^{\int_{\bbP/Y}}\\
& A
}
\end{equation*}

We now indicate how the commutativity of \eqref{diag:res-thm-p} is proved in \cite{ast117}.
For an $n$-tuple of positive integers $\underline{\alpha}=(\alpha_1,\dots,\alpha_r)$ 
one can regard fractions of the form 
${\mathrm{d}}t_1\wedge\ldots\wedge{\mathrm{d}}t_r/t_1^{\alpha_1}\dots t_r^{\alpha_r}$
as $r$-cocycles in the \v{C}ech complex 
$\check{C}^\bullet(\U,\,\omega_\pi)=
\Gamma(\bbP,\, \check{\mathcal C}^\bullet({\U},\,\,\omega_\pi))$. Let us write
$\nu(\underline{\alpha})$ for the image of this fraction in $\Hr^r(\bbP,\,\omega_\pi)$.
(Note that $\nu(1,\dots,1)=\mu$.) 
According to \cite[pp.79--80,\,Lemma\,(8.6)]{ast117} the 
natural map
\[\Hr^r_{Z_0}(\bbP,\,\omega_\pi)\to \Hr^r(\bbP,\,\omega_\pi)\]
is described by
\stepcounter{thm}
\begin{equation*}\label{map:loc-glob}\tag{\thethm}
\begin{bmatrix}
 {\mathrm{d}}t_1\wedge\ldots\wedge{\mathrm{d}}t_r\\
 t_1^{\alpha_1}\dots t_r^{\alpha_r}
 \end{bmatrix} \mapsto \nu(\underline{\alpha})
 \end{equation*}
In particular $\theta_{{}_{Z_0}}\mapsto \mu=\nu(-1,\dots,-1)$.
It is well known that if $\underline{\alpha}\neq (-1,\dots,-1)$ the \v{C}ech $r$-cocycle
${\mathrm{d}}t_1\wedge\ldots\wedge{\mathrm{d}}t_r/t_1^{\alpha_1}\dots t_r^{\alpha_r}$
for the complex $\check{C}^\bullet(\U,\,\omega_\pi)$ is a coboundary, whence
in this case $\nu(\underline{\alpha})=0$.  This establishes the commutativity of
\eqref{diag:res-thm-p}. 

If $K_{Z_0}$ is the kernel of $\res{{\bf t}}$, we have a split short exact sequence of
$A$-modules, with
$\mu\mapsto \theta_{{}_{Z_0}}$ giving the splitting:
\[
0 \lra K_{Z_0} \xrightarrow{\phantom{\,\,\res{{\bf t}}\,\,}} 
\Hr^r_{Z_0}(\bbP,\,\omega_\pi) \xrightarrow{\text{canonical}} \Hr^r(\bbP,\,\omega_{\pi})
\lra 0.
\]

\begin{prop}\label{prop:TrS=int} With the above notations we have:
\begin{enumerate}
\item[(i)] The Verdier integral for $\pi$  equals the standard trace for the relative projective space
$\bbP^r_Y$, i.e.,
\[\vin{\pi}=\int_{\bbP/Y}.\]
\item[(ii)] Let $A$ be a ring, ${\bf t}=(t_1,\dots,t_r)$ analytically independent variables over $A$,
and $J\subset A[[{\bf t]}]$ the ideal of $A[[{\bf t}]]$ generated by ${\bf t}$. Then the Verdier
integral 
$\vin{A[[{\bf t}]]/A}\colon \Hr^r_J(\omega_{A[[{\bf s}]]/A})\to A$  defined in 
\eqref{map:pint} is given by
\[
\begin{bmatrix}
\wdd{t_1}{t_r}\\
t_1^{\alpha_1},\dots, t_r^{\alpha_r} 
\end{bmatrix}
\mapsto \begin{cases}
1 & \text{when $\alpha_1=\dots = \alpha_r=1$} \\
0 & \text{otherwise.}
\end{cases}
\]
\end{enumerate}
\end{prop}

\proof For part (i) without loss of generality we may assume $Y=\Spec{\,A}$. Let
$\mu$ be the canonical generator of the free rank one $A$-module $\Hr^r(\bbP,\,\omega_\pi)$
corresponding to $1\in A$ under the isomorphism 
$\int_{\bbP/Y}\colon\Hr^r(\bbP,\,\omega_\pi)\iso A$ . It is enough to show
that $\vin{\pi}(\mu)=\int_{\bbP/Y}(\mu)$, i.e., it is enough to show that $\vin{\pi}(\mu)=1$.
According to \eqref{map:loc-glob},
the image of $\theta_{{}_{Z_0}}\in \Hr^r_{Z_0}(\bbP,\,\omega_\pi)$ in $\Hr^r(\bbP,\,\omega_\pi)$ is
$\mu$. We have
\begin{align*}
\vin\pi(\mu) & =\res{{Z_0}}(\theta_{{}_{{Z_0}}}) \qquad \text{(by \eqref{diag:res-int-4})}  \\
& = 1  \qquad \phantom{XXXX} \,\,\,\<\text{(via \Pref{prop:thom-class})}
\end{align*} 
and hence we are done for part (i).

For part (ii), let us agree to write $Y=\Spec{\,A}$.  Let us write 
$\eP=\eP^r_A$ for $\Spf{\,A[[{\bf t}]]}$, and $\wid{\pi}\colon \eP\to Y$ for the structure map.
With $\bbP$, $\pi$, $Z_0$ as above, we can identify $\eP$ with
the completion of $\bbP$ along $Z_0$. We thus have a completion map $\kappa\colon \eP \to \bbP$,
which factors through the open subscheme $U_0$ of $\bbP$ where $T_0\neq 0$ as
$\eP \to U_0 \subset \bbP$. Moreover, if $U_0$ is identified in the usual way with 
$\Spec{\,A[t_1,\dots, t_r]}$ (via $t_i=T_i/T_0$), then the first map in the factorization arises
from the inclusion of the polynomial ring $A[{\bf t}]$ into the power series ring $A[[{\bf t}]]$.
Now, by part (a),  \eqref{diag:res-int-4}, and \eqref{diag:res-thm-p}, we have
$\res{{\bf t}}=\res{{Z_0}}$. Since the composite
\[\Rr^r_{Z_0}\pi_*\omega_\pi \iso \Rp{\eP}^r{\wid{\pi}}_*\omega_{\wid{\pi}} 
\xrightarrow{\vin{\wid{\pi}}} \co_Y\]
is $\res{{Z_0}}$ by \eqref{diag:res-int-2}, and $\res{{Z_0}}=\res{{\bf t}}$, by taking global
sections we are done.
\qed

\subsection{} Let us return to our smooth map $f\colon X\to Y$ of relative dimension $r$,
and suppose $Y=\Spec{\,A}$ and $Z\hookrightarrow X$ as before a closed subscheme
such that $Z\to Y$ is an isomorphism, $Z$ lies in affine open set $U=\Spec{\,R}$ of $X$,
and $Z$ is cut out in $U$ by the vanishing or $r$ elements $t_1, \dots, t_r$ in $R$.

\begin{prop}\label{prop:s-res-indep} In the above situation $\res{{\bf t}}=\res{Z}$. In particular,
if ${\bf s}=(s_1, \dots, s_r)$ is another sequence in $R$ generating the ideal defining $Z$, then
$\res{{\bf t}}=\res{{\bf s}}$.
\end{prop}

\proof
Let $I$ be the ideal generated by ${\bf t}$. Suppose $I$ is also generated by 
${\bf s}=(s_1, \dots, s_r)$.
The completion of $R$ in the $I$-adic topology is $A[[{\bf t}]]=A[[{\bf s}]]$
and both are the completion $\wid{R}$ of $R$ in the $I$-adic topology.
It follows that $\vin{A[[{\bf t}]]/A}=\vin{\wid{R}/A}=\vin{A[[{\bf s}]]/A}$. Part (ii) of \Pref{prop:TrS=int}
and the relationship between $\res{Z}$ and $\vin{\wid{R}/A}$ then proves out assertion.
\qed

\medskip

Consider again the $A$-module decomposition
\[\Hr^r_Z(X,\,\omega_f)= \Hr^r_I(\omega_{R/A})=
\bigoplus_{\underline{\alpha}} A\bigl [\begin{smallmatrix}
\wdd{t_1}{t_r}\\
t_1^{\alpha_1},\dots, t_r^{\alpha_r} 
\end{smallmatrix}
\bigr ]
\]
with ${\underline{\alpha}}=(\alpha_1,\dots,\alpha_r)$ running over $r$-tuples of
positive integers. Each summand is a free $A$-module.
While this decomposition depends on ${\bf t}=(t_1,\dots,t_r)$, we have seen that
summand generated
by $\theta_Z=\bigl [\begin{smallmatrix}
\wdd{t_1}{t_r}\\
t_1,\dots, t_r
\end{smallmatrix}
\bigr ]$
is independent of ${\bf t}$ by \eqref{eq:dt/t=ds/s}. Moreover, since the
 the sum of the remaining summands in the direct sum is the kernel $K_Z$ of
 $\res{Z}$,  it too is independent of ${\bf t}$. Thus, we have
a canonical decompostion of $A$-modules
\stepcounter{thm}
\begin{equation*}\label{decomp:thom1} \tag{\thethm}
\Hr^r_Z(X,\,\omega_f)=K_Z\oplus A\<\<\<\cdot\<\<\theta_{{}_Z}
\end{equation*}
which is independent of ${\bf t}$ with $K_Z=\ker(\res{Z})$.

\begin{rem}\label{rem:dirac} 
Let $A$ and ${\bf t}=(t_1, \dots, t_r)$ be as in \Pref{prop:TrS=int}\,(ii). Then a little thought
shows that for $f\in A[[{\bf t}]]$, with $\mu(i_1,\dots i_r)$ the coefficient of $t_1^{i_1}\dots t_r^{i_r}$
in the power series expansion of $f$, one has the formula:
\[
\vin{A[[{\bf t}]]/A}\begin{bmatrix}
f\cdot \wdd{t_1}{t_r}\\
t_1^{\alpha_1},\dots, t_r^{\alpha_r} 
\end{bmatrix} = \mu(\alpha_1-1, \dots, \alpha_r-1).
\]
In particular, we have
\[
\vin{A[[{\bf t}]]/A}\begin{bmatrix}
f\cdot \wdd{t_1}{t_r}\\
t_1,\dots, t_r
\end{bmatrix} 
=f(0, \dots, 0).
\]
Similarly, if $A$, ${\bf t}$, $Z$, $R$ are as in \Pref{prop:s-res-indep}, then
\[
\res{Z}\begin{bmatrix}
f\cdot \wdd{t_1}{t_r}\\
t_1,\dots, t_r
\end{bmatrix} 
=\bar{f}
\]
where $\bar{f}\in A$ is the image of $f$ in $A$ under the natural surjection $R\to R/I\cong A$.
More generally, given positive integers $\alpha_1$, \dots, $\alpha_r$ one can write
\[f = \sum_{i_1,\dots,i_r}\mu(i_1, \dots, i_r)t_1^{i_1}\dots t_r^{i_r} + g,\] 
where $i_k$ are non-negative integers such that
$\sum_j i_j < \alpha_1+\dots +\alpha_r$, $\mu(i_1, \dots, i_r)\in A$, and
$g\in I^{\alpha_1+\dots +\alpha_r}$. In this case we have
\[
\res{Z}\begin{bmatrix}
f\cdot \wdd{t_1}{t_r}\\
t_1^{\alpha_1},\dots, t_r^{\alpha_r} 
\end{bmatrix} = \mu(\alpha_1-1, \dots, \alpha_r-1).
\]

\end{rem}

\subsection{The Verdier residue for sections of smooth maps}\label{ss:std-res}
 Now suppose $Y$ is not necessarily affine, and as above we have a closed subscheme
$Z\hookrightarrow X$ is such that $Z\to Y$ is an isomorphism. Let $z\in Z$ be a point.
Pick affine open subschemes $U'$ in $X$ and $V'$ in $Y$ such that $z\in U'$ and $f(U')\subset V'$ 
and such that $U'\cap Z$ is given in $U'$ by the vanishing of $r$-elements 
$t_1,\dots, t_r\in \Gamma(U',\co_X)$. Let $V=f(U'\cap Z)$. Then $V$ is an affine open subscheme
of $V'$ (since it is isomorphic to $U'\cap Z$ which, being a closed subscheme of $U'$ is affine).
Moreover $U'\to V'$ is affine, whence $U\set f^{-1}(V)\cap U'$ is affine. Note that 
$U'\cap Z = U\cap Z$,  $f(U)=V$, $Z\cap U$ is given by the vanishing of $t_1,\dots,t_r$ and
$Z\cap U\to V=f(Z\cap U)$ is an isomorphism. Thus locally we can reduce to the situation in
\Ssref{ss:affine-section}. If $Z_U=U\cap Z$, then from \eqref{eq:dt/t=ds/s}, it is clear that
$\theta_{{}_{Z_U}}$ glue to give a section $\theta_{{}_Z}$ of $\Rr^r_Zf_*\omega_f$:
\[\theta_{{}_Z}\in \Gamma(X,\,\Rr^r_Zf_*\omega_f).\]
Moreover, the $A$-module $K_Z$ in \eqref{decomp:thom1} being independent of ${\bf t}$ means
that its construction globalizes to give a quasi-coherent submodule $\eK_Z$ of
$\Rr^r_Z f_*\omega_f$. Finally, since the decomposition \eqref{decomp:thom1} is
canonical, it globalizes to give a decompostion:
\stepcounter{thm}
\begin{equation*}\label{decomp:thom2} \tag{\thethm}
\Rr^r_Z f_*\omega_f= \eK_Z\oplus (\co_Y\<\<\<\<\cdot\<\<\theta_{{}_Z}).
\end{equation*}

\begin{thm}\label{thm:sres=res} Let $Z$ be a closed subscheme of $X$ such that $Z\to Y$ is an isomorphism. Then $\res{Z}$ is the composite
\[ 
\Rr^r_Zf_*\omega_f = \eK_Z\oplus (\co_Y\<\<\<\<\cdot\<\<\theta_{{}_Z}) 
\xrightarrow{\text{projection}}
\co_Y\<\<\<\<\cdot\<\<\theta_{{}_Z} \iso \co_Y
\]
where the direct sum decomposition is \eqref{decomp:thom2} and the last isomorphism 
is $\theta_{{}_Z}\mapsto 1$.
\end{thm}

\proof Without loss of generality we may assume $X=\Spec{\,R}$, $Z=\Spec{\,R/I}$ where
$I$ is an ideal of $R$ generated by $r$ elements $\{t_1,\dots,t_r\}$ and $Y=\Spec{\,A}$.
The result then follows from \Pref{prop:s-res-indep} and the explicit description of $\res{\bf t}$.
\qed

\smallskip

Before stating the next theorem we need some notation.
If $\psi\colon \omega_f[r] \to \fs\co_Y$ is a map of $\co_X$-modules, then
 $\bar\psi \colon \omega_f\to \omgs{f}$ will denote the map $\bar\psi =\Hr^{-r}(\psi)$. 
 We remind the reader that ${\bf v}={\bf v}_{\<\<\<{}_f}$ denotes the Verdier isomorphism
 $\omega_f[r]\iso \fs\co_Y$. We alert the reader to one notational issue.
 In this subsection, for good book-keeping purposes we
 will write $\bar{\bf v}\colon \omega_f\iso \omgs{f}$ for $\Hr^{-r}({\bf v})$. For most of the
 paper we do not put the ``bar" over ${\bf v}$ for this map, as that abuse of notation is usually 
 harmless. (Cf.\:also \Rref{rem:v-local}.)
 
 \begin{lem}\label{lem:section-iso} Let $Z$ be a closed subscheme of $X$ such that $Z\to Y$ is finite
 and flat. Suppose we have an isomorphism $\psi\colon \omega_f[r]\iso \fs\co_Y$ 
 such that the composite
 \[\Rr^r_Zf_*\omega_f \xrightarrow[\text{{\em via $\bar{\psi}$}}]{\Iso} \Rr^r_Zf_*\omgs{f} 
 \xrightarrow{\ares{Z}} \co_Y\]
 is the residue map $\res{Z}$. Then there is an open neighbourhood $U$ of $Z$ 
 in~$X$ such that $\psi\vert_U = {\bf v}\vert_U$.
 \end{lem}
\proof It is enough to prove that there is an open neighbourhood $U$ of $Z$ such
that $\bar{\psi}\vert_U=\bar{\bf v}\vert_U$. 
Let $\varphi\colon \omgs{f} \iso \omgs{f}$ be the automorphism given by
$\varphi = \bar{\bf v} \smcirc  \bar{\psi}^{-1}$. Let $\kappa\colon \X=X_{/Z} \to X$ be the
completion of $X$ along $Z$ and $\wid{f}=f\smcirc\kappa$. By the hypothesis
we have $\ares{Z}\smcirc\Rr^r_Zf_*(\varphi)=\ares{Z}$,
and hence by definition of $\tin{\wid{f}}$ we get
$\tin{\wid{f}}\smcirc\Rp{\X}^r{\wid{f}}_*(\kappa^*(\varphi))=\tin{\wid{f}}$.
Thus by local duality \cite[Cor.\,5.1.4]{fub-abs}, 
we see that $\kappa^*(\varphi)$ is the identity map, whence there is an open neighbourhood
$U$ of $Z$ such that $\varphi\vert_U$ is the identity map.
\qed

\medskip

We need a little more notation in order to state the next Lemma. Consider a cartesian diagram
\stepcounter{thm}
\[
\begin{aligned}\label{diag:res-psi}
\xymatrix{
X' \ar@{}[dr]|{\square} \ar[d]_{f'} \ar[r]^v & X \ar[d]^{f} \\
Y' \ar[r]_{u} & Y
}
\end{aligned}\tag{\thethm}
\]
where $f$ is smooth (and hence Cohen-Macaulay) of relative dimension $r$.
We will use the notation of \cite{cm} and denote by
\[\theta_u^{{}_f}\colon v^*\omgs{f} \iso \omgs{f'}\]
the corresponding base change isomorphism (see \cite[p.740,\,Theorem\,2.3.5\,(a)]{cm}). 
Now suppose we have
an isomorphism $\psi\colon \omega_f[r] \iso \fs\co_Y$ and suppose
$Z\hookrightarrow X'$ is a closed subcheme such that $Z\to Y'$ is an isomorphism. We write
\stepcounter{thm}
\begin{equation*}\label{map:res-psi}\tag{\thethm}
\res{{\psi,Z}}\colon \Rr^r_Zf'_*\omega_{f'} \to \co_{Y'}
\end{equation*}
for the composite:
\[\Rr^r_Z f'_*\omega_{f'} =\Rr^r_Z f'_*v^*\omega_{f} 
\xrightarrow{\text{via $\bar{\psi}$}} \Rr^r_Z f'_* v^*\omgs{f} \xrightarrow[\text{via }\theta_u^{{}_f}]{\Iso} 
 \Rr^r_Z f'_* \omgs{f'}
\xrightarrow{\ares{Z}} \co_{Y'}.\]
In other words
$\res{{\psi,Z}} = \ares{Z}\smcirc \Rr^r_Zf_*(\theta_u^{{}_f}\smcirc v^*({\bar \psi}))$.

\begin{lem}\label{lem:et-section-iso} Let $u\colon Y'\to Y$ be an \'etale map, and let
$X'$, $f'$, $v$, $\theta_u^{{}_f}$ be as above.
Suppose we have an isomorphism
$\psi\colon \omega_f[r] \iso \fs\co_Y$, and a 
closed subscheme $Z$ of $X'$ such that $Z\to Y'$ is finite and flat and such that 
$\res{{\psi,Z}}=\res{Z}$. Then there is an open neighbourhood $U$ of the locally
closed set $v(Z)$ such that $\psi\vert_U={\bf v}_{\<\<\<{}_f}\vert_U$.
\end{lem}

\proof
By the hypothesis on $\psi$ and
by \eqref{lem:section-iso} we can find
an open neighbourhood $V$ of $Z$ in $X'$ such that
 $(\theta_u^{{}_f}\smcirc v^*({\bar{\psi}}))\vert_{V}=\bar{\bf{v}}_{\<\<\<{{}_{f'}}}\vert_{V}$.
On the other hand, by \cite[p.740,\,Theorem\,2.3.5\,(b)]{cm}, we have
$\theta_u^{{}_f}\smcirc v^*\bar{\bf v}_{\<\<\<{}_f}=\bar{\bf{v}}_{\<\<\<{{}_{f'}}}$.
It follows that $v^*({\bar{\psi}}))\vert_{V}= v^*\bar{\bf v}_{\<\<\<{}_f}\vert_V$.
Set $U=v(V)$. Since $v$ is \'etale, $U$ is open, and $V\to U$ is faithfully flat, whence
$\bar{\psi}\vert_U = \bar{\bf v}_{\<\<\<{}_f}\vert_U$.
\qed

\begin{rem}\label{rem:dense} Since $f$ is smooth, if
$x$ is an associated point of $X$, then $y=f(x)$ is an associated point of $Y$, and
$x$ is a generic point of the fibre $f^{-1}(y)$. This means that if an open subscheme
$V$ of $X$ is such that $V\cap f^{-1}(s)$ is dense in $f^{-1}(s)$ for every associated point
$s$ of $Y$, then $V$ is scheme theoretically dense in $X$, since it contains every
associated point of $X$. We use this fact in what follows.
\end{rem}

\begin{thm}\label{thm:main} Let $\psi\colon \omega_f[r] \iso \fs\co_Y$. A necessary
and sufficient condition that $\psi$ is the Verdier isomorphism ${\bf v}_{\<\<\<{}_f}$
is the following:\\
For every \'etale map $u\colon Y'\to Y$ and every closed subscheme
$Z$ of $X'$ such that $Z\to Y'$ is an isomorphism, we have
$\res{{\psi,Z}}=\res{Z}$. Here $X'$, $f'$, $v$ are as in diagram \eqref{diag:res-psi}.
\end{thm}

\proof For  $u\colon Y'\to Y$, $f'\colon X'\to Y'$, $v\colon X'\to X$ as above,
according to \cite[p.740,\,Theorem\,2.3.5\,(b)]{cm} we have
$\theta_u^{{}_f}\smcirc v^*\bar{\bf v}_{\<\<\<{}_f}=\bar{\bf{v}}_{\<\<\<{{}_{f'}}}$.
The necessity part of the theorem then follows from \Tref{thm:sres=res}.

Conversely, suppose we have an isomorphism $\psi\colon \omega_f[r] \iso \fs\co_Y$ satisfying
the condition stated in the theorem. We have to show that $\bar{\psi}=\bar{\bf v}_{\<\<\<{}_f}$.
Fix $y\in Y$. Since $f$ is smooth, the set $W_y$ of points $x\in f^{-1}(y)$ such that $k(x)$ is finite and separable over $k(y))$, is dense in $f^{-1}(y)$ by \cite[p.\,42, \S2.2, Cor.\,13]{blr}. 
Let $x$ be such a point. We can find an
\'etale map $u\colon Y'\to Y$ such that (with the usual notations) there is a section of $f'$ passing
through a point $x'$  satisfying $v(x')=x$ \cite[p.\,43, \S2.2, Prop.\,14]{blr}. 
Let $Z$ be the image of this section. Then $Z$ is closed,
and $Z\to Y'$ is an isomorphism, whence by our hypotheses on $f$ and by
\Lref{lem:et-section-iso} there is an open neighbourhood $U$ of $v(Z)$ on which 
$\bar{\psi}=\bar{\bf v}_{\<\<\<{}_f}$. Since $x\in v(Z)$, this equality holds in an open neighbourhood
of $x$. Varying $x$ over $W_y$, and varying $y$ over $Y$, by \Rref{rem:dense}
the equality holds in a scheme theoretically
dense open subset of $X$ and hence everywhere, for $\omega_f$ and $\omgs{f}$ are
invertible $\co_X$-modules.
\qed

Recall that given a point $x\in X$, closed in its fibre, with $k(x)$  separable over $k(f(x))$,
 since $f$ is smooth we can find an \'etale 
neighbourhood $Y' \to Y$ of $f(x)$ and a section of $f'$ (with the usual notations for base change
that we have been following) passing through one of the points of $v^{-1}(x)$. It is immediate that
one can find an open cover $\{U_\alpha\}$ of $Y$, \'etale surjective maps 
$u_\alpha\colon Y_\alpha \to U_\alpha$, such that (with $X_\alpha\set X\times_YY_\alpha$, and 
$f_\alpha\colon X_\alpha\to Y_\alpha$, $v_\alpha\colon X_\alpha\to X$ the projections) there is
a closed subscheme $Z_\alpha$ of $X_\alpha$ which maps isomorphically on to $Y_\alpha$.
Let $Y'=\coprod_\alpha Y_\alpha$, $X'=\coprod_\alpha X_\alpha$,
$f'=\coprod_\alpha f_\alpha$,  $u=\coprod_\alpha u_\alpha$.
Then we have a closed subscheme $Z$ of $X'$ such that $Z\to Y'$
is an isomorphism (take $Z=\coprod_\alpha Z_\alpha$). Note that $u\colon Y'\to Y$ is
\'etale and {\em surjective}, whence it is {\em faithfully flat}.

\begin{prop} Let $\psi\colon \omega_f[r] \iso \fs\co_Y$ be an isomorphism.
\begin{enumerate}
\item[(a)] If the fibres of $f$ are connected, and $Z$ is a closed subscheme of $X$ such that 
$Z\to Y$ is an isomorphism and $\ares{Z}\smcirc \Rr^r_Zf_*(\bar{\psi}) = \res{Z}$, then 
$\psi={\bf v}_{\<\<\<{}_f}$.
\item[(b)] Suppose the fibres of $f$ are geometrically connected. Then $\psi={\bf v}_{\<\<\<{}_f}$
if and only if there is an \'etale surjective map $u\colon Y'\to Y$ and
(with the usual notation) a closed subscheme $Z$ of $X'=X\times_YY'$ 
with $Z\to Y'$ an isomorphism such that
$\res{{\psi,Z}}=\res{Z}$.
\end{enumerate}
\end{prop}
\proof For part (a), we note that if $\kappa\colon \X\to X$ is the completion of $Z$ along $X$,
then $\kappa^*\bar{\psi}=\kappa^*\bar{{\bf v}_{\<\<\<{}_f}}$. We therefore have an open sunscheme
$V$ containing $Z$ such that $\bar{\psi}\vert_V=\bar{{\bf v}_{\<\<\<{}_f}}\vert_V$. Since $f$ is smooth,
it has a (locally) a factorization $f=\pi\smcirc h$, where $h$ is \'etale and 
$\pi$ is the structural map ${\mathbb A}^r_Y\to Y$. 
Since the fibres of $f$ are connected, and $f^{-1}(y)\cap V\supset f^{-1}(y)\cap Z
\neq \emptyset$, it follows that $V\cap f^{-1}(y)$ is dense in $f^{-1}(y)$. Thus $V$ is scheme-theoretically dense in $X$ by \Rref{rem:dense}.
Now ${\bar{\psi}}^{-1}\smcirc \bar{{\bf v}_{\<\<\<{}_f}}$ is the
identity automorphism on $\omega_f$ on $V$, which is scheme theoretically
dense on $X$, and $\omega_f$ is invertible on $X$. It follows that
it ${\bar{\psi}}^{-1}\smcirc \bar{{\bf v}_{\<\<\<{}_f}}$ is the identity automorphism on all of $X$.

For part (b), first suppose $\psi={\bf v}_{\<\<\<{{}_f}}$. By the remarks made above the statement
of the theorem, there is an \'etale surjective map $u\colon Y'\to Y$, and (with the usual
meaning attached to $X'$, $f'$ and $v$) a closed subscheme $Z$ of $X'$ such that $Z\to Y'$
is an isomorphism. Now 
$\res{{\psi,Z}}= \ares{Z}\smcirc \Rr^r_Zf_*(\theta_u^{{}_f}\smcirc v^*({\bar \psi}))
=\ares{Z}\smcirc \Rr^r_Zf_*(\theta_u^{{}_f}\smcirc v^*({\bar{\bf v}}_{\<\<\<{}_f}))$. On the other hand,
by \cite[p.740,\,Theorem\,2.3.5\,(b)]{cm}, ${\bf v}$ behaves well with respect to base change, i.e.,
$\theta_u^{{}_f}\smcirc v^*({\bar{\bf v}}_{\<\<\<{}_f})=\bar{\bf v}_{\<\<\<{}_{f'}}$. Thus
$\res{{\psi,Z}}=\ares{Z}\smcirc\Rr^r_Zf_*(\bar{\bf v}_{\<\<\<{}_{f'}})=\res{Z}$. 

Conversely, suppose we have an \'etale surjective map $u\colon Y'\to Y$ and a closed subscheme
$Z$ of $X=X\times_YY'$, 
with $Z\to Y'$ an isomorphism satisfying $\res{{\psi,Z}}=\res{Z}$.  Let $f'\colon X'\to Y'$
and $v\colon X'\to X$ be the projections. Since the fibres
of $f'$ are connected, by part (a) we have 
$\theta_u^{{}_f}\smcirc v^*(\bar{\psi})=\bar{\bf v}_{\<\<\<{{}_{f'}}}$. Now,
$\bar{\bf v}_{\<\<\<{{}_{f'}}}= \theta_u^{{}_f}\smcirc v^*(\bar{\bf v}_{\<\<\<{{}_f}})$ (by
\cite[p.740,\,Theorem\,2.3.5\,(b)]{cm} again) from which it is immediate that
$v^*(\bar{\psi})=v^*(\bar{\bf v}_{\<\<\<{{}_{f}}})$. The map $v\colon X'\to X$ is
\'etale surjective, and hence faithfully flat, giving the result.
\qed

\section{\bf Regular Differential Forms}\label{s:reg-diffs} The results in this section do
not affect the results
in the rest of the paper, and so may be skipped on first reading. 
These results are  here to give a non-trivial application of the 
characterisation of Verdier's map in the previous 
section. The main results of this section, connecting the Kunz-Waldi regular differentials with Verdier's 
isomorphism, are proved
again in \Ssref{ss:reg-diff-again} without making use of the results in \cite{kw} or the results in
this section. There is
a fleeting reference to definition of the map \eqref{map:can} of this section in \Ssref{ss:reg-diff-again}
(see in \S\S\S\,\ref{sss:reg-ver}).

There is a fleeting reference to
the definition of the map \eqref{map:can} of this section in \Ssref{ss:reg-diff-again}. 
 All schemes in this section, unless otherwise stated, are
ordinary schemes. The aim is to relate the concrete form of Grothendieck duality via Kunz's
regular differential forms to Verdier's isomorphism. In somewhat greater detail, regular differential
forms defined for certain types of maps $f\colon X\to Y$ are concrete representations of many
aspects of Grothendieck duality. A well-known special case is that of
Rosenlicht's differentials on singular cuves \cite{rosen}. Kunz defined generalization of these 
to more general situations (higher dimensions) in a series
of papers, and in \cite{kw}, Kunz and Waldi defined the sheaf of (relative) regular differentials
for dominant finite type equidimensional  maps $f\colon X\to Y$ between excellent schemes
which are do not have embedded components. When such an $f$ is generically smooth,
this was related to duality theory by Kunz, Lipman, H\"ubl, Sastry (see 
\cite{ast117}, \cite{hk1}, \cite{hk2}, \cite{ajm}). All the papers just mentioned work within the
framework of a simpler version of Grothendieck duality (one eschewing derived categories)
due to Kleiman \cite{kl}. We now review this, taking a slightly revisionist view, in that we interpret
the principal objects ($r$-dualizing pairs) in terms of the full blown duality theory of Grothendieck. 

\subsection{Overview of Kleiman's functor}\label{ss:kleiman} 
Regular differentials are a vast generalisation of the differentials Rosenlicht used for describing
describing duality for singular curves \cite{rosen}. To put the theory in context one we give a quick account
of Kleiman's theory of $r$-dualizing pairs given in  \cite{kl}.
Let $f\colon X\to Y$ be a proper map such that $\dim({X\otimes k(y)})\le r$ for every $y \in Y$. For
any scheme $Z$, let $Z_{qc}$ denote the category of quasi-coherent $\co_Z$-modules. 
According to [{\em loc.\:cit.,} pp.\,41--42, Definition\,(1)], an {\em $r$-dualizing pair} $(f^K,\,t_f)$ 
consists of a covariant functor $f^K\colon Y_{qc}\to X_{qc}$ and a natural transformation
$t_f\colon \Rr^rf_*f^K\to {\bf 1}{\<\<{}_{Y_{qc}}}$ inducing a bifunctorial isomorphism of quasi-coherent 
sheaves,
\[f_*\sHom_X(\eF,\,f^K\eG)\iso \sHom_Y(\Rr^rf_*\eF,\,\eG)\] 
for each $\eF\in X_{qc}$ and each $\eG\in Y_{qc}$. Kleiman explicitly eschewed derived categories
in his paper, and shows the existence of an $r$-dualizing pair (for $f$ of the kind we are
considering) using the special adjoint functor theorem. From our point of view $f^K$ can
identified with $\Hr^{-r}(\fs(\boldsymbol{-}))$. Our hypotheses on $f$ ensure that
$\Hr^j(\fs(\boldsymbol{-}))=0$ for $j< -r$, whence we get a map
(of functors from $Y_{qc}$ to $\Dqc(X)$) $f^K(\boldsymbol{-})[r]\to \fs$. The map $t_f$
is then the composite
\[\Rr^rf_*f^K=\Hr^0(\Rfs f^K(\boldsymbol{-})[r])\lra
\Hr^0(\Rfs\fs\co_Y) \xrightarrow{\Tr{f}} \Hr^0(\co_Y)=\co_Y.\]

When $f$ is not proper, $f^K$ still makes sense (even if $f$ is not
compactifiable, i.e., even if $f$ is not separated), since $\Hr^n(\fs)$ makes sense even for every
integer $n$, even if $\fs$ is not defined (see comment above 
\cite[Def.\,4.1.2]{fub-abs}), and hence one
can set $f^K=\Hr^{-r}(\fs)$. If $f$ is Cohen-Macaulay of relative dimension $r$,
$f^K\co_Y=\omgs{f}$, and if further $f$ is smooth we have, via Verdier's isomorphism
$\omega_f\iso f^K\co_Y$. In the proper, Cohen-Macaulay case we have
 $(\omgs{f},\,\tin{f})=(f^K\co_Y,\,t_f(\co_Y))$. If $f$ is in addition smooth, we have a unique
 isomorphism of pairs $(\omega_f,\,\vin{f}) \iso (f^K\co_Y,\,t_f(\co_Y))$. 

\subsection{Regular Differentials}\label{ss:reg-diffs} Let $f\colon X\to Y$
be a finite type map. Following Kunz in \cite[B.17]{kd} say it is {\em equidimensional of dimension $r$}
if 
\begin{itemize}
\item the generic point of $X$ are mapped to the generic points of $Y$ ,and 
\item the non-empty fibres
of $f$ are such that irreducible component of these fibres are all of dimension $r$.
\end{itemize}

Now suppose the map $f\colon X\to Y$ satisfies the following conditions
\begin{itemize}
\item $X$ and $Y$ are excellent schemes, and neither have embedded points amongst their
associated points;
\item  $f$ is equidimensional of dimension $r$, and
\item the smooth locus of $f$ is scheme-theoretically dense in $X$ (which, given our hypotheses,
means that the smooth locus of $f$ contains all the generic points of $X$).
\end{itemize}
Next let $X_0$ be the artinian scheme
\[X_0=\coprod_s\Spec{\,\co_{X,s}}\]
where $s$ runs through the set of associated (= maximal in this case) points and 
$i_X\colon X_0\to X$ the natural affine map. Similarly, we have the artinian scheme $Y_0$
constructed out of the generic points of $Y$, and an affine map $i_Y\colon Y_0\to Y$.
We write
\[k(X)={i_X}_*\co_{X_0}\]
where as before $s$ runs over generic points of $X$.
The sheaf of relative meromorphic $r$-forms $\Omega^r_{k(X)/k(Y)}$
on $X$ is then the quasi-coherent $\co_X$-module
given by the formula
\[\Omega^r_{k(X)/k(Y)}={i_X}_*\Omega^r_{X_0/Y_0} = \omega_f\otimes_{\co_X}k(X).\]

Under our hypotheses on $f$ the $\co_X$-module of {\em regular differentials}
$\oreg{f}$ (denoted
$\omega^r_{X/Y}$ in \cite{hk1}, \cite{hk2}, and \cite{ajm})
 is defined in
\cite[\S\,3,\,\S\,4]{kw}. It is coherent and is an $\co_X$ submodule
of the module of meromorphic $r$-differentials
$\Omega^r_{k(X)/k(Y)}$, and hence is torsion-free.
On the smooth locus $X^s$ of $f$, and writing $f^s\colon X^s\to Y$
for the smooth map obtained by restricting $f$, we have 
$\oreg{f}\vert_{X^s}= \oreg{f^s}=\omega_{f^s}$.

When $f$ is {\em proper} we have a trace map (denoted $\int_{X/Y}$ in \cite{hk1}, \cite{hk2},
\cite{ajm})
\[\rin{f}\colon \Rr^rf_*\oreg{f}\to \co_Y.\]
This map is defined when $f$ is {\em projective} in \cite{hk1}, 
and is generalized to proper $f$ in \cite{ajm}. One of the main results of \cite{ajm} is that
the resulting map $\oreg{f}\to f^K\co_Y$ is an isomorphism (a fact proved in \cite{hk2} for
projective maps $f$). There is also a notion of a residue map $\Rr^r_Zf_*\oreg{f}\to\co_Y$
(denoted $\int_{X/Y,Z}$ in \cite{hk1}, \cite{hk2}, and \cite{ajm}) for certain special closed
 subschemes $Z$ of $X$ which are finite over $Y$ (see 
 \cite[pp.77--78, Assumption\,4.3 and Theorem\,4.4]{hk1}).
 
 To avoid notational confusion we denote this
 \[\res{Z}^{\mathrm{reg}}\colon \Rr^r_Zf_*\oreg{f}\to \co_Y.\]
 
 \subsection{Summary of the main result of \cite{ajm}}\label{ss:ajm}
 The complete statement concerning $\oreg{f} (=\omega^r_{X/Y})$, $\rin{f} (=\int_{X/Y})$
 and $\res{Z}^{\mathrm{reg}} (=\int_{X/Y,Z})$ can be found in \cite[pp.750--752, Theorem]{ajm}.
 In brief, here are the main points of this result:
 \begin{enumerate}
 \item One has a canonical isomorphism $\varphi=\varphi_{\<{}_f}\colon \oreg{f}\iso f^K\co_Y$ 
 such that when $f$ is proper $\varphi$ is the unique isomorphism for which the diagram
 \[
 \xymatrix{
 \Rr^rf_*\oreg{f} \ar@/_0.75pc/[drr]_{\rin{f}} \ar[rr]^-{\text{via $\varphi$}} && \Rr^rf_*f^K\co_Y 
 \ar[d]^{t_f(\co_Y)}\\
 && \co_Y
 }
 \]
 commutes \cite[pp.750--751, (i) (The Duality Theorem)]{ajm}.
 \item The isomorphism $\varphi_{\<{}_f}$ is compatible with open immersions into $X$. In greater
 detail, if $j\colon U\to X$ is an open immersion, as submodules of the $\co_X$ module
 $\Omega_{k(U)/k(Y)}$, $i^*\oreg{f}=\oreg{fi}$ and the diagram
 \[
 \xymatrix{
 i^*\oreg{f} \ar@{=}[d] \ar[rr]^{\Iso}_{i^*\varphi_{\<{}_f}} && i^*f^K\co_Y
 \ar[d]^{\,\rotatebox{-90}{\makebox[-0.1cm]{\Iso}}}\\
 \oreg{fi} \ar[rr]^{\Iso}_{\varphi_{\<{}_{fi}}} && (fi)^K\co_Y
 }
 \]
 commutes [Ibid, pp.750--751, (i) and (ii)].
 \item If $Z$ is a closed subscheme of $X$ satisfying Assumptiom\,4.3 of \cite[p.77]{hk1}
 then the diagram
 \[
 \xymatrix{
 \Rr^r_Zf_*\oreg{f} \ar@/_.75pc/[drr]_{\res{Z}^{\mathrm{reg}}} 
 \ar[rr]^-{\text{canonical}} && \Rr^rf_*\oreg{f} 
 \ar[d]^{\rin{f}}\\
 && \co_Y
 }
 \]
 commutes \cite[p.752, (iii) (The Residue Theorem)]{ajm}.
 \item If $Z$ is a closed subscheme of $X$ such that $Z$ lies in the smooth locus of $f$
 and $Z\to Y$ is an isomorphism, then $\res{Z}^{\mathrm{reg}}=\res{Z}$ (see 
 \cite[p.62, Cor.\:1.13]{hk1} and \cite[p.78, 4.4]{hk1} as well as the formulae in \Rref{rem:dirac}).
 \item The map $\varphi$ is compatible with flat base change to excellent schemes without
 embedded associated points \cite[3.13]{kw} and \cite[pp.751--752, (ii) and (iv)]{ajm}. 
 In particular $\varphi$ is compatible with \'etale base change. 
 \end{enumerate}
 
\subsection{Regular Differentials and Verdier}\label{ss:reg-ver}
Now suppose $f$ is smooth. Then $f^K\co_Y=\omgs{f}$ and $\oreg{f}=\omega_f$.
 Let $\psi=\varphi[r]$. Identifying $\fs\co_Y$ with $\omgs{f}[r]$ we have an isomorphism
 \[\psi\colon \omega_f[r] \iso \fs\co_Y.\]
 Then using the notations of \Ssref{ss:std-res}, we have $\varphi=\bar{\psi}$.
 In light of the properties listed above for $\varphi_{\<{}_f}$ and $\oreg{f}$ we see that
 if $u\colon Y'\to Y$ is an \'etale map and $Z$ is a closed subscheme of $X'=X\times_YY'$
 such that $Z\to Y'$ is an isomorphism, then $\res{Z}^{\mathrm{reg}}=\res{Z}$. However,
 the left side is the map $\res{{\psi,Z}}$ of \eqref{map:res-psi}, whence we conclude from
 \Tref{thm:main} that $\psi={\bf v}_{\<\<\<{}_f}$, where the right side is the Verdier isomorphism
 of \eqref{def:verdier}.
 
 Our next observation is one that was made by J.\:Lipman in pp.\,33--34 of \cite{ast117} for varieties
 over fields in his discussion leading to Lemma (2.2) of \textit{ibid}. Suppose $f$ is as in the previous
 subsection, and $U$ is the smooth locus of $f$. Let $j\colon U\to X$ be the open immersion
 and $g=f\smcirc j\colon U\to Y$ the resulting smooth map. By our hypotheses, $U$ contains all the associated
 points of $X$, whence it is scheme theoretically dense.
 Without getting into the notions of canonical structures and dualizing structures, we have
 a composition
 \stepcounter{thm}
 \begin{equation*}\label{map:can}\tag{\thethm}
 f^K\co_Y \hookrightarrow j_*g^K\co_Y \xrightarrow{j_*{\bar{\bf{v}_{\<\<\<{}_g}}}^{-1}} 
 j_*\omega_g \hookrightarrow \Omega^r_{k(X)/k(Y)}
 \end{equation*}
with every arrow an inclusion since $f^K\co_Y$, $g^K\co_Y$ and $\omega_g$ 
are torsion free and $j^*f^K\co_Y \iso g^K\co_Y$. The image of $f^K\co_Y$ in 
$\Omega^r_{k(X)/k(Y)}$ must be $\oreg{f}$ since $\bar{\bf{v}_{\<\<\<{}_g}}$ is
$\varphi_{\<\<\<{}_g}$ of item (1) of \S\S\,\ref{ss:ajm}. In greater detail, if $\bar{\omega}$
is the image of $f^K\co_Y$ in $\Omega^r_{k(X)/k(Y)}$ under \eqref{map:can}, 
and $\alpha\colon f^K\co_Y\iso \bar{\omega}$ the resulting isomorphism, then we have
an isomorphism $\beta\colon \oreg{f}\iso \bar{\omega}$ such that 
$\alpha=\beta\smcirc\varphi_{\<\<\<{}_f}$. 
Now $j^*\beta={\bf 1}_{\omega_g}$ since $\varphi_{\<\<\<{}_g}= \bar{\bf v}_{\<\<\<{}_g}$.
Since $U$ is scheme theoretically dense in $X$ and the sheaves involved are torsion free,
the assertion follows. In other words Verdier's isomorphism gives us the regular differential
forms of Kunz and Waldi, as well as the dualizing structure on them.

Here is the formal statement of the result(s) we just proved.

\begin{thm}\label{thm:reg-ver} Let $f\colon X\to Y$ be a finite type map between excellent
schemes such that $X$ and $Y$ have no embedded points, $f$ is equidimensional of dimension
$r$, and the smooth locus of $f$ contains all the associated points of $X$ (i.e., the smooth locus
of $X$ is scheme-theoretically dense in $X$). 
\begin{enumerate}
\item[(a)] If $f$ is smooth then the map $\varphi_{\<\<\<{}_f}$ of item (1) in \Ssref{ss:ajm}
is the Verdier isomorphism $\bar{\bf v}_{\<\<\<{}_f}$ defined in \eqref{iso:pre-verdier}.
\item[(b)] If $j\colon U\to X$ is the open immersion from the smooth locus of $f$ to $X$, and
$g\colon U\to Y$ is the composite $g=f\smcirc i$, then the module of regular differential $r$-forms
$\oreg{f}$ of Kunz and Waldi \cite[\S\,3, \S\,4]{kw} is the image of $f^K\co_Y$ under injective
composite \eqref{map:can}. Moreover the resulting isomorphism $f^K\co_Y\iso \oreg{f}$ is inverse of
the map $\varphi_{\<\<\<{}_f}$.
\end{enumerate}
\end{thm}

\section{\bf Transitivity for smooth maps}\label{s:transitivity}

\subsection{The map $\zeta_{g,f}$ between differential forms} 
Suppose $f\colon\X\to\Y$ and $g\colon\Y\to\Z$ are maps in $\bbG$, with $f$ a smooth map of relative
dimension $e$ , and $g$ a smooth map of relative dimension $d$. We have a map of
differential forms
\stepcounter{thm}
\begin{equation*}\label{map:zeta-gf}\tag{\thethm}
\zeta_{g,f}\colon f^*\omega_g[d]\otimes_{\co_\X}\omega_f[e] \lra \omega_{gf}[d+e]
\end{equation*}
defined by the commutativity of the following diagram
\[
{\xymatrix{
f^*\omega_g[d]\otimes_{\co_\X}\omega_f[e] \ar[rr]^-{\zeta_{g,f}} 
\ar[d]^-{\,\rotatebox{-90}{\makebox[-0.1cm]{\Iso}}}_-{f^*{{\bf{v}}}_{\<\<{}_{g}}\otimes 
{{\bf{v}}}_{\<\<{}_{f}}} 
&& \omega_{gf}[d+e] \ar[d]_{\,\rotatebox{-90}{\makebox[-0.1cm]{\Iso}}}^{{\bf{v}}_{\<\<{}_{gf}}} \\
f^*\ush{g}\co_{\Z}\overset{\bL}{\otimes}_{\co_\Y}\ush{f}\co_\Y \ar[rr] _-{\chi_{{}_{[g,f]}}}
&& \ush{(gf)}\co_\Z
}}
\]
where $\chi_{{}_{[g,f]}}\colon f^*\ush{g}\co_{\Z}\overset{\bL}{\otimes}_{\co_\Y}\ush{f}\co_\Y
\to \ush{(gf)}\co_\Z$ is the map defined in 
\cite[Def.\,7.2.16]{fub-abs}.

\begin{prop}\label{prop:zeta-gf} The following hold:
\begin{enumerate}
\item[(a)] {\em (Flat Base Change)} Suppose 
\[
{\xymatrix{
\U \ar[r]^u \ar[d]_p \ar@{}[dr]|\square & \X \ar[d]^f\\
\V \ar[r]_v \ar[d]_q \ar@{}[dr]|\square & \Y \ar[d]^g \\
\W \ar[r]_w & \Z
}}
\]
is a cartesian square with $u$ flat, $f$ and $g$ smooth and in $\bbG$. Then
\[u^*\zeta_{g,f} = \zeta_{g,f}.\]
\item[(b)] Suppose $f\colon \X\to \Y$ and $g\colon\Y\to\Z$ are smooth maps and in $\bbG$.
Let $\kappa\colon \X^*\to \X$ be the completion of $\X$ with respect to an open coherent ideal.
Then 
\[\zeta_{{}_{g,f\kappa}}=\kappa^*\zeta_{{}_{g,f}}.\]
\item[(c)] Suppose $\X\xrightarrow{f} \Y_1 \xrightarrow{\kappa} \Y_2 \xrightarrow{g} \Z$ is
a sequence of maps in $\bbG$ with $f$ and $g$ smooth, and $\kappa$ a completion map with
respect to an open coherent ideal.
Then
\[\zeta_{g,\kappa f} = \zeta_{g\kappa,f}.\]
\item[(d)] Suppose
\[
{\xymatrix{
\wid{\X} \ar[d]_{\wid{f}} \ar[r]^{\kappa_2}  & \X \ar[d]^{f}\\
\wid{\Y}  \ar[dr]_{\wid{g}} \ar[r]^{\kappa_1} & \Y \ar[d]^g \\
& \Z
}}
\]
is a commutative diagram in $\bbG$ with $f$ and $g$ smooth, and
 $\kappa_1$ and $\kappa_2$ completions
with respect to open coherent ideals of $\co_\Y$ and $\co_\X$ respectively. Then
\[\kappa_2^*\zeta_{{g,f}}=\zeta_{\wid{g},\wid{f}}.\]
\end{enumerate}
\end{prop}

\proof Follows from the properties for $\chi_{{}_{g,f}}$ listed in 
\cite[\S\S\,7.2]{fub-abs}, 
\Tref{thm:verd-complete},  \Tref{prop:verd-complete2}, and the fact that the Verdier
isomorphism is compatible with flat base change.
\qed
\subsection{The map $\varphi_{g,f}$ between differential forms}
For a smooth map between ordinary schemes $f\colon X\to Y$
of relative dimension $d$, $\omega_f\set \wedge^d_{\co_X}\Omega^1_{X/Y}$.
Let $X$, $Y$, and $Z$ be ordinary schemes.
Suppose $f\colon X\to Y$ is a smooth map of schemes
of relative dimension $d$ and $g\colon Y\to Z$ is smooth of relative dimension $e$. Let
\stepcounter{thm}
\begin{equation*}\label{iso:phi-gf-bar}\tag{\thethm}
\bar{\varphi}_{g,f}\colon f^*\omega_g\otimes\omega_f \iso \omega_{gf}
\end{equation*}
be the map which is locally given by 
\[
f^*(\wdd{t_1}{t_e})\otimes\wdd{s_1}{s_d}\mapsto
\wdd{s_1}{s_d}\wedge\wdd{t_1}{t_e}. 
\]
Here ${\bf t}=(t_1,\dots,t_e)$ and ${\bf s}=(s_1,\dots, s_d)$ are local relative ``co-ordinates", i.e.,
${\bf t}$ gives an \'etale map $U\to {\mathbb A}^e_Z$ on an open subscheme  $U$ of $Y$, and on 
an open sunscheme $V$ of $f^{-1}(U)$, ${\bf s}$ gives an \'etale map $V\to {\mathbb A}^d_Y$.
The local map given above (i.e., $f^*({\mathrm{d}}{\bf t})\otimes {\mathrm{d}}{\bf s} \mapsto
{\mathrm{d}}{\bf s}\wedge {\mathrm{d}}{\bf t}$) is independent of these local relative co-ordinates and 
hence globalises to give $\bar{\varphi}_{g,f}$.

Using the recipe that gives us $\psi$ in \cite[(8.1.2)]{fub-abs}
 from $\bar{\psi}$ 
we get a well defined isomorphism
in $\Dc(X)$
\stepcounter{thm}
\begin{equation*}\label{iso:phi-gf}\tag{\thethm}
\varphi_{g,f} \colon f^*\omega_g[e]\otimes_{\co_X} \omega_f[d] \iso \omega_{gf}[d+e]. 
\end{equation*}
Note that 
\stepcounter{thm}
\begin{equation*}\label{iso:phi-bar-phi}\tag{\thethm}
\Hr^{-(d+e)}(\varphi_{g,f})=\bar{\varphi}_{g,f}
\end{equation*}
and hence one can go back and forth between $\bar{\varphi}_{g,f}$ and $\varphi_{g,f}$.

Here is the main theorem:

\begin{thm}\label{thm:fubini} Let $f\colon \X\to \Y$ and $g\colon \Y\to \Z$ be maps in $\bbG$
which are smooth. Then
\[\zeta_{g,f}=\varphi_{g,f}.\]
\end{thm}

\proof 
We divide the proof into cases.

{\bf{Case 1.}} Let $A$ be a noetherian ring, $u_1,\ldots, u_d$, $v_1,\ldots, v_e$ analytically independent
variables over $A$, and consider the $A$-algebras $R$, $S$, and $T$ given by
$R=A[[u_1,\,\ldots,\,u_d]]$, $S=R[[v_1,\ldots, v_e]]=A[[u_1,\,\ldots,\,u_d, v_1,\ldots, v_e]]$. 
Let $I={\bf u}R$ be the $R$-ideal generated by ${\bf u}=(u_1,\ldots,u_d)$, $J={\bf v}S$ the
$S$-ideal generated by ${\bf v}=(v_1,\ldots,v_e)$, and $L=IS+J$. 
In other words, $J$ is the $S$-ideal generated by $({\bf u},\,{\bf v})$.

Suppose $\X=\Spf{(S,L)}$, $\Y=\Spf{(R,\,I)}$,  $\Z=Z=\Spec{\,A}$, and that our smooth maps
$f\colon \X\to \Y$,  $g\colon \Y\to \Z=Z$ are the natural maps corresponding to the maps of
 adic rings $(R,I)\to (S,L)$ and $(A,0) \to (R,I)$. 

We have additional schemes, namely $Y=\Spec{\,R}$,  and  $\V=\Spf{(S,\,J)}$. The natural maps 
between the adic rings involved give us a commutative diagram with the square on top being a 
cartesian square:
 \[
{\xymatrix{
\X \ar[d]_{f} \ar[r]^{\kappa'} \ar@{}[dr]|\square & \V \ar[d]^p \\
\Y \ar[dr]_{g} \ar[r]^\kappa & Y^* \ar[d]^q\\
& Z
}}
\]
The maps $\kappa$ and $\kappa'$are completion maps and $f$, $g$, $p$, $q$ are the
obvious maps. Note that $p$, $f$, and $g$ are smooth and pseudo-proper (however this is not
true for the map $q$, which is not of pseudo-finite-type unless $d=0$).

The rank one free $\co_\X$-modules $\omega_f$ and $\omega_{gf}$ correspond to the universal finite
$S$-module of differentials $\omega_{S/R}\set\widehat{\Omega}_{S/R}^e$ and 
$\omega_{S/A}\set\widehat{\Omega}_{S/A}^{d+e}$. The rank one free $\co_\Y$-module $\omega_g$
corresponds to the universal finite $R$-module of degree $d$ differentials 
$\omega_{R/A}\set \widehat{\Omega}_{R/A}^d$. Thus
\begin{align*}
\omega_{S/R}&=S\,\<\wdd{v_1}{v_e}, \\
\omega_{S/A}&=S\,\<\wdd{u_1}{u_d}\wedge\wdd{v_1}{v_e}, \\
\omega_{R/A}&=R\,\<\wdd{u_1}{u_d}.
\end{align*}

The $S$-module $\omega_{S/R}$ gives us a rank one free $\co_\V$-module. A little thought
shows us that this module is in fact $\omega_p$. Define $\omega_q$ as the rank one free 
$\co_Y$-module corresponding to $\omega_{R/A}$. The equations $\Gamma(\X,\,\omega_f)=
\omega_{S/R}=\Gamma(\V,\,\omega_p)$ and $\Gamma(\Y,\,\omega_g)=\omega_{R/A}
=\Gamma(Y,\,\omega_q)$ can be re-written as 
\[
\omega_f=(\kappa')^*\omega_p\,\,{\text{and}}\,\,\, 
\omega_g=\kappa^*\omega_q.
\]

%

Write $\bar{\varphi}$ and $\varphi$ for (the global sections of) the maps $\bar{\varphi}_{g,f}$
and $\varphi_{g,f}$. Then the $S$-module isomorphism
\[\bar{\varphi}\colon \omega_{R/A}\otimes_R\omega_{S/R} \iso \omega_{S/A}\]
is given by 
$\bar{\varphi}(\mathrm{d}{\bf u}\otimes\mathrm{d}{\bf v})=\mathrm{d}{\bf v}\wedge\mathrm{d}{\bf u}$.
We have the following formula, where
$\vin{A[[{\bf u}, {\bf v}]]/A}$ and $\vin{R[[{\bf v}]]/R}$ are as in \eqref{map:pint}.
\[
\vin{A[[{\bf u}]]/A}\begin{bmatrix}
\vin{R[[{\bf v}]]/R}\begin{bmatrix}
{\mathrm{d}}{\bf v}\\
v_1^{\beta_1},\ldots,v_e^{\beta_e}
\end{bmatrix} {\mathrm{d}}{\bf u}\\
u_1^{\alpha_1},\ldots, u_e^{\alpha_e}
\end{bmatrix}
=
\vin{A[[{\bf{u,\,v}}]]/A}\begin{bmatrix}
\bar{\varphi}({\mathrm{d}}{\bf u}\otimes{\mathrm{d}}{\bf v})\\
v_1^{\beta_1},\ldots,v_e^{\beta_e}, u_1^{\alpha_1},\ldots,u_d^{\alpha_d}
\end{bmatrix}\\. \leqno{(*)}
\]
Indeed, if any of the $\alpha_l$'s or $\beta_k$'s is not equal to $1$, then both sides
equal zero. If $\alpha_l=\beta_k=1$ for $l=1,\ldots,d$, $k=1,\ldots,e$, both sides equal $1$. This
means that the following diagram commutes.
\[
{\xymatrix{
\Hr^{d+e}_L(\omega_{R/A}\otimes_R\omega_{S/R}) 
\ar[d]_{{\text{\cite[(8.2.4)]{fub-abs}}}}^{\,\rotatebox{-90}{\makebox[-0.1cm]{\Iso}}}
\ar[r]^-{\Iso}_-{\bar{\varphi}} & \Hr^{d+e}_L(\omega_{S/A}) \ar[dd]^{\vin{A[[{\bf{u,\,v}}]]/A}}\\
\Hr^d_I(\omega_{R/A}\otimes_R\Hr^e_{J}(\omega_{S/R})) \ar[d]_{{\text{via}}\,\vin{R[[{\bf{v}}]]/R}} &  \\
\Hr^d_I(\omega_{R/A}) \ar[r]_{\vin{A[[{\bf u}]]/A}} & A
}}
\] 
If, in the above  diagram, we replace ${\bar{\varphi}}$  by $\bar{\zeta}_{g,f}$, then by 
\cite[Prop.\,8.3.1\,(b)]{fub-abs},
the resulting diagram commutes.  (See also 
\cite[Remark 8.3.4]{fub-abs}.)
By the universal property of the 
pair $(\omgs{S/A},\vin{S/A})$ we see that $\bar{\varphi}=\bar{\zeta}_{g,f}$, i.e., 
$\varphi_{g,f}=\zeta_{g,f}$.

{\bf Case 2.} Suppose we have a section $\sigma\colon \Z\to \X$ and $\tau\set f\smcirc\sigma$,
and $\X$ and $\Y$ are the completions of $\X$ and $\Y$
along the closed subschemes given by the closed
immersions $\sigma\colon \Z\hookrightarrow \X$ and $\tau\colon \Z\hookrightarrow \Y$ respectively.
More precisely, if $\I_1\subset \co_\X$ and $\I_2\subset \co_\Y$ are the coherent ideals giving
the embeddings of $\Z$ into $\X$ and $\Y$ (via $\sigma$ and $\tau$), and $\I\subset \co_\Z$
is an ideal of definition of  $\X$, then $\I\co_\X+\I_1$ and $\I\co_\Y+\I_2$ are ideals of definition
of $\X$ and $\Y$ respectively. Since the source and target of $\zeta_{g,f}$ and $\varphi_{g,f}$
are concentrated in one degree, the question of their equality is a local question on $\X$
and hence, without loss of generality, we may assume that the schemes involved are
affine,  say $\X=\Spf{(S,L)}$, $\Y=\Spf{(R,I)}$ and 
$\Z=\Spf{(A,I_0)}$ respectively. In fact we may assume that $\tau$ and $\sigma$ are given
by regular sequences $(u_1,\ldots, u_d)$ and $(u_1,\ldots, y_d, v_1,\ldots, v_e)$ respectively,
and ${\bf u}$ is analytically independent over $A$, and ${\bf v}$ is analytically independent
over $R$. We then have a cartesian diagram (where the power series rings $A[[u_1,\ldots, u_d]]=
A[[{\bf u}]]$ and $R[[v_1,\ldots, v_e]]=R[[{\bf v}]]$ are given the adic topologies from the ideals
$(u_1,\ldots, u_d)$ and $(v_1,\ldots, v_e)$ repsectively)
\[
{\xymatrix{
\X \ar[rr]^v \ar[d]_f \ar@{}[drr]|\square && \Spf{\,R[[\bf v}]] \ar[d]^p \\
\Y \ar[rr]^u \ar[d]_g \ar@{}[drr]|\square && \Spf{\,A[[\bf u}]] \ar[d]^q \\
\Z \ar[rr]_w && \Spec{\,A}
}}
\]
with the horizontal arrows being the natural ones. Note that $w$ is flat being a completion map.
Therefore flat base change applies (see \Pref{prop:zeta-gf}\,(a)) and we have
 $\zeta_{g,f} =v^*\zeta_{q,p}$. Clearly $\varphi_{g,f}=v^*\varphi_{q,p}$ from the explicit description
 of $\varphi_{q,p}$ and $\varphi_{g,f}$.  By Case 1, we have $\zeta_{q,p}=\varphi_{q,p}$.
Applying $v^*$ to both sides, we get the result for this case.

{\bf Case 3 (The General Case).} In the general case, let $\Y\times_\Z\X=\eQ$, $\X\times_\Z\X=\eP$,
and let $p\colon \eP\to \eQ$, $q\colon \eQ\to \X$ be the base changes of $f$ and $g$, and
let $\pi_i\colon \eP\to \X$ and $\pi_2$ be the projections $\X\times_\Z\X\to \X$, with 
$\pi_1=q\smcirc p$. It $\De\colon \X\to \eP$ is the diagonal immersion, then let 
$\kappa\colon \U \to \X$ be the completion of $\X$ with respect to $\De(\X)$, and
let $\kappa'\colon \V\to \Y$ be the completion of $\Y$ along $(p\smcirc\De)(\X)$.
We have a natural map $\wid{p}\colon \U \to \V$ such that
$\kappa'\smcirc \wid{p}=p\smcirc\kappa$. Let  $\wid{q}=q\smcirc\kappa'$ and let
$\delta\colon \X\to \U$ be the natural closed immersion.
We then have a commutative diagram with the two rectangles on the right being cartesian:
\[
{\xymatrix{
\X \ar[dr]^\delta \ar@/^1.5pc/[drrr]^{{\bf 1}_\X} \ar@/_3pc/[dddrr]_{{\bf 1}_\X}& & & \\
&\U \ar[r]^\kappa \ar[d]_{\wid{p}} & \eP \ar[d]^p \ar[r]^{\pi_2} \ar@{}[dr]|\square & \X \ar[d]^f \\
& \V \ar[r]^{\kappa'} \ar[dr]_-{\wid{q}} & \eQ \ar[d]^{q} \ar[r] \ar@{}[dr]|\square & \Y \ar[d]^g \\
& & \X \ar[r]_f & \Z 
}}
\]
Now, using the explicit formula for $\varphi_{g,f}$, $\varphi_{q,p}$, and $\varphi_{\wid{q},\wid{p}}$, we
see that $\pi_2^*\varphi_{g,f}=\varphi_{q,p}$ and $\kappa^*\varphi_{q,p}=\varphi_{\wid{q},\wid{p}}$.
Thus $\kappa^*\pi_2^*\varphi_{g,f}=\varphi_{\wid{q},\wid{p}}$.

A similar relationship holds for the $\zeta_{\bullet,\bullet}$ maps. Indeed, by \Pref{prop:zeta-gf}\,(a)
we have $\pi_2^*\zeta_{g,f}=\zeta_{q,p}$ and by \Pref{prop:zeta-gf}\,(d) we have
$\kappa^*\zeta_{q,p}=\zeta_{\wid{q},\wid{p}}$, giving
$\kappa^*\pi_2^*\zeta_{g,f}=\zeta_{\wid{q},\wid{p}}$.  On the other hand, by Case\,2 considered 
above, we have $\zeta_{\wid{q},\wid{p}}=\varphi_{\wid{q},\wid{p}}$. Thus
\[\kappa^*\pi_2^*\zeta_{g,f}=\kappa^*\pi_2^*\varphi_{g,f}.\]
Applying $\delta^*$ to both sides of this equation, and noting that 
$\pi_2\smcirc\kappa\smcirc\delta={\bf 1}_\X$, we get the result.
\qed

\section{\bf Applications of Transitivity}
\subsection{Iterated residues} Suppose $f\colon X\to Y$ is smooth of relative
dimension $e$, $g\colon Y\to Z$ smooth of relative dimension $d$, $W_1\hookrightarrow X$
a closed subscheme, finite and flat over $Y$, $W_2\hookrightarrow Y$ a closed subscheme
which is finite and flat over $Z$. Let $W=W_1\cap f^{-1}(W_2)\hookrightarrow X$. Suppose
further that $W_1$ is cut out by a quasi-regular sequence ${\bf v}=(v_1,\dots, v_e)$ in $S$
and $W_2$ is cut out by a quasi-regular sequence ${\bf u}=(u_1,\dots, u_d)$ in $R$.

\begin{thm}\label{thm:res-res} In the above situation, for $\nu\in\Gamma(\co_Y,\,\omega_g)$,
$\mu\in\Gamma(\co_X,\,\omega_f)$, we have
\[
\res{{W_2}}\begin{bmatrix}
\res{{W_1}}\begin{bmatrix} \mu \\ 
v_1^{\beta_1},\,\ldots,\,v_e^{\beta_e}
\end{bmatrix} \nu \\ 
u_1^{\alpha_1},\,\ldots,\,u_d^{\alpha_d}
\end{bmatrix}= \res{W}\begin{bmatrix} \mu\wedge f^*\nu \\ 
v_1^{\beta_1},\,\ldots,\,v_e^{\beta_e},\,u_1^{\alpha_1},\,\ldots,\,u_d^{\alpha_d}
\end{bmatrix}
\]
where, for notational simplicity, we denote the image of $u_i$ in $S$ also by $u_i$.
\end{thm}

\proof Recall from the definition of $\zeta_{g,f}$ in \eqref{map:zeta-gf} that $\zeta_{g,f}$ is the 
transform of $\chi_{[g,\,f]}$ after applying Verdier's isomorphism to $f^!\co_Y$, $g^!\co_Z$
and $(gf)^!\co_Z$. From \Tref{thm:fubini} and 
\cite[(8.3.2)]{fub-abs} we get
\[
\res{{W_2}}\begin{bmatrix}
\res{{W_1}}\begin{bmatrix} \mu \\ 
v_1^{\beta_1},\,\ldots,\,v_e^{\beta_e}
\end{bmatrix} \nu \\ 
u_1^{\alpha_1},\,\ldots,\,u_d^{\alpha_d}
\end{bmatrix}= \res{W}\begin{bmatrix} \bar{\varphi}_{g,f}(\nu\otimes\mu)\\ 
v_1^{\beta_1},\,\ldots,\,v_e^{\beta_e},\,u_1^{\alpha_1},\,\ldots,\,u_d^{\alpha_d}
\end{bmatrix}.
\]
The result then follows from the definition of $\bar{\varphi}_{g, f}$ in \eqref{iso:phi-gf-bar}.
\qed

 \subsection{The Restriction Formula} \label{ss:restriction}
 An important application of our transitivity result is the so-called {\emph{Restriction Formula}}, namely
 the formula in \Cref{cor:res-thm} below.
The formula is related to the following problem. Suppose
\stepcounter{thm}
 \[
 \begin{aligned}\label{diag:restr-thm}
 {\xymatrix{
 X\, \ar@{^(->}[r]^i \ar[dr]_f & P \ar[d]^\pi \\
 & Y
 }}
 \end{aligned}\tag{\thethm}
 \]
 is a commutative diagram of ordinary schemes, with $\pi$ and $f$ smooth and separated and
 $i$ a closed immersion. Let the
 relative dimension of $\pi$ be $n=d+e$ and the relative dimension of $f$ be $e$. As usual, let
 $\eN_i^d$ be the $d$-th exterior power of the normal bundle $\eN_i$ of $X$ in $P$. 
 We have, via  Verdier's isomorphism and the isomorphism 
 $\eta_i'$ of \eqref{iso:eta'-i}, an isomorphism
 \stepcounter{thm}
 \begin{equation*}\label{map:aXP}\tag{\thethm}
 a_{{}_{X/P}}\colon i^*\omega_\pi[n]\otimes_{\co_X}\eN_i^d[-d] \iso \omega_f[e],
 \end{equation*}
 defined as the composite
 \[
 \begin{aligned}\label{def:aXP}
i^*\omega_\pi[n]\otimes_{\co_X}\eN^d_i[-d] &\xrightarrow[\phantom{\Iso}]{\eta_i'} i^! \omega_\pi[n] \\
 & \xrightarrow[\phantom{\Iso}]{{\bf v}_{\<\<{}_\pi}} i^!\pi^!\co_Y= i^!\omgs{\pi}[n]\\
 & \xrightarrow{\Iso} f^!\co_Y =\omgs{f}[e]\\
 & \xrightarrow{{\bf v}_{\<\<{}_f}^{-1}} \omega_f[e].
 \end{aligned}\tag{\thethm}
 \]
 The question then is, what is the concrete form of $a_{{}_{X/P}}$ in terms of local relative coordinates?
 We answer the question in \Tref{thm:res-thm} below. 
 
 We leave it to the reader to check that $a_{{}_{X/P}}$ is compatible with open immersions into $Y$ and
 $P$. Indeed every map in the composition defining $a_{{}_{X/P}}$ is well behaved with respect
 to open immersions into $Y$ and $P$. Thus we may assume that $P$ and $Y$ are affine and
 that $X$ is defined by an ideal generated by a quasi-regular sequence, which is part of a relative
 system of co-ordinates for $\pi\colon P \to Y$. 
 
 We now assume that $Y=\Spec{\,A}$, $P=\Spec{\,S}$, and $X=\Spec{\,R}$ where $R=S/I$,
 and $I$ is generated by a quasi-regular sequence ${\bf t}=(t_1,\dots, t_d)$ in $S$, 
 and there is an \'etale map
 $A[T_1,\ldots,\,T_d,\,V_,\ldots,\,V_e] \to S$,  (where $T_l$, $l=1,\ldots,d$, and $V_k$, $k=1,\ldots,e$
 are algebraically independent variables), and $t_i$ is the image of $T_i$ for $i=1,\ldots,d$. This
 can always be achieved by shrinking $Y$ and $P$ (see \cite[pp.\,39--40, Prop.\,7(c)]{blr}). 
 Now every $\mu\in \omega_{S/A}$ can be
 written uniquely as
 \stepcounter{thm}
 \begin{equation*}\label{eq:mu=nu-dt}\tag{\thethm}
 \mu= \wdd{t_1}{t_d}\wedge\nu
 \end{equation*}
 with $\nu\in \wedge^e_S\Omega^1_{S/A}$.
 Define
 \stepcounter{thm}
 \begin{equation*}\label{eq:b-X/P}\tag{\thethm}
 b_{{}_{X/P}}\colon i^*\omega_\pi\otimes_{\co_X}\eN_i^d \iso \omega_f
 \end{equation*}
 by the formula 
\[
 \mu\otimes 1/{\bf t} \mapsto i^*\nu
 \]
 where $\mu$ and $\nu$ are related by 
 \eqref{eq:mu=nu-dt}. We should clarify that $i^*\nu\in\omega_{R/A}$ is the pull-back of
 $\nu$ as a differential form.
 In other words, $i^*\nu$ is the image of $\nu$ under the composite of maps
 $\wedge^e_S\Omega^1_{S/A}\to R\otimes_S \wedge^e_S\Omega^1_{S/A} \to 
 \wedge^e_R\Omega^1_{R/A}=\omega_{R/A}$.
 
 In what follows, let 
 \stepcounter{thm}
 \begin{equation*}\label{eq:a-X/P}\tag{\thethm}
 \bar{a}_{\<\<{}_{X/P}}\colon i^*\omega_\pi\otimes_{\co_X}\eN^d_i\to \omega_f
 \end{equation*}
 be the map 
 \[\bar{a}_{\<\<{}_{X/P}}=\Hr^0(a_{\<\<{}_{X/P}}).\]
The notation follows the conventions we have been using throughout,
 and as observed earlier,
 $a_{{}_{X/P}}$ can be recovered from $\bar{a}_{{}_{X/P}}$ 
 (see \cite[\S\S\,8.1]{fub-abs}).
 
 \begin{thm}\label{thm:res-thm} Under the above assumptions on $i$, $\pi$, $f$, $A$, $S$, and $R$,
 we have 
 \[\bar{a}_{{}_{X/P}}=b_{{}_{X/P}}.\]
 \end{thm}
 
 \proof 
Since we will be using results from \cite{fub-abs} which
have been stated in terms of the abstract dualizing sheaves of the form
$\omgs{f}$, rather than in terms of $\omega_f$, it is convenient for us to have an analogue
of $a_{{}_{X/P}}$ taking values in $\omgs{f}[e]$. To that end, suppose
$k\colon W\hookrightarrow P$ is a regular immersion of codimension $m\le n$, and that
$g=\pi\smcirc k$ is flat over $Y$, 
so that $g\colon W\to Y$  is Cohen-Macaulay of relative dimension
$n-m$. Define
\[\ush{a}_{\<\<\<\<\<\<\<{}_{W/P}}\colon k^*\omega_\pi[n]\otimes_{\co_W}\eN^m_k[-m] \iso 
\omgs{g}[n-m] \]
as the composite:
\begin{align*}
k^*\omega_\pi[n]\otimes_{\co_W}\eN^m_k[-m] 
&\xrightarrow[\phantom{\Iso}]{\eta_k'} k^! \omega_\pi[n] \\
 & \xrightarrow[\phantom{\Iso}]{{\bf v}_{\<\<{}_\pi}} k^!\pi^!\co_Y= k^!\omgs{\pi}[n]\\
 & \xrightarrow{\Iso} g^!\co_Y =\omgs{g}[n-m]
 \end{align*}
 (Similarly we have an isomorphism $\ush{a}_{\<\<\<\<\<\<\<{}_{W/X}}$ for a regular
 immersion $W\hookrightarrow X$ of $Y$-schemes such that $W\to Y$ is flat.)
 If $W\to Y$ is
 finite (in addition to being flat), so that $m=n$,
 and $W$ is given by the vanishing of ${\bf v}=(v_1,\ldots,v_n)$, then by the definition of
 the map $\ttr{g,\pi,k}$ in  \cite[(5.3.2)]{fub-abs}, we have:
 \[\tin{g}\smcirc h_*(\ush{a}_{\<\<\<\<\<\<\<{}_{W/X}})= \ttr{g,\pi,k}. \leqno{(\dag)}\]

 By definition of $a_{{}_{X/P}}$, it is clear that 
 ${\bf v}_{\<\<{}_f}\smcirc a_{{}_{X/P}}=\ush{a}_{\<\<\<\<\<\<\<{}_{X/P}}$.  
 
 We first prove that the map $\bar{a}_{{}_{X/P}}$  is compatible
 with base change. In greater detail, suppose $u\colon Y'\to Y$ is a map, $X'\set X\times_YY'$, 
 $P'\set P\times_YY'$, and let $f'\colon X'\to Y'$, $\pi'\colon P'\to Y'$, $i'\colon X'\to P'$,
 $w\colon P'\to P$, $v\colon X'\to X$ be the resulting maps obtained from base change. We then
 clearly have $v^*( i^*\omega_\pi\otimes_{\co_X}\eN_i^d)=  
 {i'}^*\omega_{\pi'}\otimes_{\co_{X'}}\eN_{i'}^d$ and $v^*\omega_f=\omega_{f'}$. We claim
 that $v^*\bar{a}_{{}_{X/P}}=\bar{a}_{{}_{X'/P'}}$.
 In what
 follows, $Y'=\Spec{\,A'}$, $R'=R\otimes_AA'$, $S'=S\otimes_AA'$.
 
 Consider the base change isomorphism $\theta=\theta_u^f\colon v^*\omgs{f} \iso \omgs{f'}$ 
 of part (a) of
 \cite[p.740,\,Theorem\,2.3.5\,(a)]{cm}. According to loc.cit.\,(b) we have a commutative
 diagram
 \[
 {\xymatrix{
 v^*\omega_f \ar@{=}[d] \ar[rrr]^{\Iso}_{v^*(\bar{\bf{v}}_{\<\<{}_f})} &&& v^*\omgs{f} 
 \ar[d]^{\theta}_{\,\rotatebox{-90}{\makebox[-0.1cm]{\Iso}}} \\
 \omega_{f'} \ar[rrr]^{\Iso}_{\bar{\bf{v}}_{\<\<{}_{f'}}} &&& \omgs{f'}
 }}
 \]
 Since, ${\bf v}_{\<\<{}_f}\smcirc a_{{}_{X/P}}=\ush{a}_{\<\<\<\<\<\<\<{}_{X/P}}$, from the above
 diagram we see that it is enough to show that 
 $\theta\smcirc v^*\ush{a}_{\<\<\<\<\<\<\<{}_{X/P}}=\ush{a}_{\<\<\<\<\<{}_{X'/P'}}$ in order
 to show that $v^*a_{{}_{X/P}}=a_{{}_{X'/P'}}$.
 
 Let $\eta_i$ and $\eta'_i$ be the maps defined in (C.2.11) and (C.2.13) of \cite{fub-abs}.
 For the next few lines, all labels of the form (B.x.y) or (C.x.y) refer to the labels in 
 \cite[Appendix]{fub-abs}. By definition, $\eta_i'=(\textup{B}.1.2)\smcirc \eta_i$. It follows that
 ${\ush{\bar a}}_{\<\<\<\<\<{}_{X/P}}$ is the composite of isomorphisms:
 \[i^*\omega_\pi\otimes_{\co_X}\eN_i^d \xrightarrow[\text{(C.2.7)}]{\Iso}  \Ext^d_{\co_P}(\co_X,\,\omega_\pi)
 \xrightarrow[(\textup{B}.1.2)]{\Iso} \Hr^0(i^!\omgs{\pi}[n]) \iso \omgs{f}.\]
 Let the composite of the last two maps in the above composition be denoted
 $c_{{}_{X/P}}\colon  \Ext^d_{\co_P}(\co_X,\,\omega_\pi) \iso \omgs{f}$. Consider the
 diagram
  \[
 {\xymatrix{
  v^*(i^*\omega_\pi\otimes_{\co_X}\eN_i^d) \ar@{=}[d] \ar[r]^-{\Iso}_-{\text{(C.2.7)}} &
  v^*\Ext^d_{\co_P}(\co_X,\,\omega_\pi) \ar[d]^{\,\rotatebox{-90}{\makebox[-0.1cm]{\Iso}}} 
  \ar[r]^-{\Iso}_-{c_{{}_{X/P}}} &  v^*\omega_{f}
  \ar[d]^{\theta}_{\,\rotatebox{-90}{\makebox[-0.1cm]{\Iso}}}\\
 {i'}^*\omega_{\pi'}\otimes_{\co_{X'}}\eN_{i'}^d \ar[r]^-{\Iso}_-{\text{(C.2.7)}}  &
 \Ext^d_{\co_{P'}}(\co_{X'},\,\omega_{\pi'}) \ar[r]^-{\Iso}_-{c_{{}_{X'/P'}}} & \omega_{f'} 
 }} \leqno{(\ddag)}
 \]
 where the isomorphism in the middle is the natural one, which we now describe. Let
$Q^\bullet \to R$ be a projective resolution of the $S$-module $R$. Then
$Q^\bullet\otimes_RR'=Q^\bullet\otimes_AA' \to R\otimes_AA'=R'$ is an $S'$-projective resolution
of the $S'$-module $R'$. Now 
\begin{align*}
\Homb_S(Q^\bullet,\,\omega_{S/A}[d])\otimes_AA' & =
\Homb_{S'}(Q^\bullet\otimes_AA',\,\omega_{S/A}[d]\otimes_AA')\\
& = \Homb_{S'}(Q^\bullet\otimes_AA',\,\omega_{S'/A'}[d]).
\end{align*}
Since ${\bf t}$ is a quasi-regular sequence in 
$S$,  we can (and will) pick $Q^\bullet$ to be the version of the
Koszul homology complex on ${\bf t}$ such that $\Homb_S(Q^\bullet,\,S)=K^\bullet(\bf{t})$, and
the equality $\Homb_S(Q^\bullet,\,\omega_{S/A}[d])\otimes_AA' = 
\Homb_{S'}(Q^\bullet\otimes_AA',\,\omega_{S'/A'}[d])$ reduces to the well-known equality
$\omega_{S/A}[d]\otimes_SK^\bullet({\bf t})\otimes_SS' 
= \omega_{S'/A'}[d]\otimes_{S'}K^\bullet({\bf t}')$,
where ${\bf t}'={\bf t}\otimes 1$. By right-exactness of tensor products, we get:
\begin{align*}
\Hr^0(\omega_{S/A}[d]\otimes_SK^\bullet({\bf t}))\otimes_SS' &= 
\Hr^0(\omega_{S/A}[d]\otimes_SK^\bullet({\bf t})\otimes_SS') \\
& = \Hr^0(\omega_{S'/A'}[d]\otimes_{S'}K^\bullet({\bf t}')).
\end{align*}
The isomorphism $v^*\Ext^d_{\co_P}(\co_X,\,\omega_\pi) \iso 
\Ext^d_{\co_{P'}}(\co_{X'},\,\omega_{\pi'})$ then follows from the isomorphism in
\cite[(C.2.3)]{fub-abs}.
(See also the proof of Lemma\,1 of \cite[pp.39--40]{lip-dp} as well as \cite[p.762,\,(8.9)]{cm} 
for the case when $X\hookrightarrow P$ is not necessarily a regular immersion, but $R$ is relatively
Cohen-Macaulay over $A$). The description of the isomorphism we have given also shows
that the rectangle on the left in diagram $(\ddag)$ above commutes.
 The rectangle on the right commutes by \cite[p.741,\,Theorem\,2.3.6]{cm}. Thus 
 diagram $(\ddag)$ commutes and $\theta_u^f\smcirc v^*\ush{\bar{a}}_{\<\<\<{}_{X/P}}=
 \ush{\bar{a}}_{\<\<\<{}_{X'/P'}}$, i.e.,  
 $v^*\bar{a}_{\<\<{}_{X/P}}=\bar{a}_{\<\<{}_{X'/P'}}$.  
 
 Next suppose we have a closed subscheme $j\colon Z\hookrightarrow X$ such that (a)
 $h=f\smcirc j\colon Z\to Y$ is an isomorphism, and (b), if ${\bar L}\subset R$ is the ideal of $R$ which
 defines $Z$, then ${\bar L}$ is generated by a quasi-regular sequence 
 $\bar{\bf u}=(\bar{u}_1,\ldots,\bar{u}_e)$. Let $u_i\in S$ be lifts of ${\bar u}_i\in R$ for
 $i=1,\ldots, e$. Let $L$ be the ideal generated by $({\bf t,u})$. Then $L$ is the ideal defining
 the closed immersion $ij\colon Z\hookrightarrow P$. Let $B=\Gamma(Z,\,\co_Z)$. 
 For a sequence of positive integeres ${\bf m}=(m_1,\ldots, m_e)$, let 
 ${\bf u^m}=(u_1^{m_1},\,\ldots,\,u_e^{m_e})$, 
 ${\bf\bar{u}^m}=(\bar{u}_1^{m_1},\,\ldots,\,\bar{u}_e^{m_e})$, 
 $L_{\bf m}$ the $R$-ideal generated by ${\bf\bar{u}^m}$, $B_{\bf m}=R/L_{\bf m}$,
 $Z_{\bf m}=\Spec{\,B_{\bf m}}$ and $j_{\bf m}\colon Z_{\bf m} \hookrightarrow X$ the
 natural closed immersion. Let $h_{\bf m}\colon Z_{\bf m}\to Y$ be the finite flat
 map $h_{\bf m}=f\smcirc j_{\bf m}$. Finally let $\kappa\colon \X\to X$ be the completion
 of $X$ along $Z$.
 
 For $\mu\in \omega_{S/A}$, and positive integers $m_i$, $i=1,\ldots,e$, it is easy to see that
 \[
 \res{{Z,\pi}}\begin{bmatrix} 
 \mu \\
 t_1,\,\ldots,\,t_d, u_1^{m_1},\,\ldots,\,u_e^{m_e}
 \end{bmatrix}
 =\res{{Z,f}}\begin{bmatrix} 
 b_{\<\<{}_{X/P}}(\mu\otimes 1/{\bf t}) \\
 {\bar u_1}^{m_1},\,\ldots,\,{\bar u}_e^{m_e}
 \end{bmatrix}. \leqno{(*)}
 \]
 Indeed, we can write $\mu$ in a unique manner as 
 $\mu=f\wdd{t_1}{t_d}\wedge\wdd{u_1}{u_e}$, with $f\in S$. Then 
 $b_{\<\<{}_{X/P}}(\mu\otimes 1/{\bf t})={\bar f}\wedge\wdd{\bar{u}_1}{\bar{u}_e}$, where
 $\bar{f}$ is the image of $f$ in $R$.
 Both sides of $(*)$ are then
 realised as the coefficient of $u_1^{m_1-1}u_2^{m_2-1}\ldots u_e^{m_e-1}$ in the power
 series expansion of $f$, whence$(*)$ holds.
 On the other hand, according to 
 \cite[Prop.\,C.6.6]{fub-abs}, 
 \[
 \ush{\bar{a}}_{\<\<{}_{Z_{\bf m}/P}}
 (\mu\otimes 1/({\bf t, u^m})) =
 \ush{\bar{a}}_{\<\<{}_{Z_{\bf m}/X}}(\bar{a}_{\<\<{}_{X/P}}
 (\mu\otimes 1/{\bf t}))\otimes 1/{\bf \bar{u}^m}).
 \]
 Apply $\tin{h_{\bf m}}\smcirc {h_{\bf m}}_*$ to both sides. 
 By $(\dag)$ and \cite[Prop.\,5.4.4]{fub-abs}, this yields,
 \[
 \res{{Z, \pi}}\begin{bmatrix} 
 \mu \\
 t_1,\,\ldots,\,t_d, u_1^{m_1},\,\ldots,\,u_e^{m_e}
 \end{bmatrix}
 =\res{{Z, f}}\begin{bmatrix} 
 \bar{a}_{\<\<{}_{X/P}}(\mu\otimes 1/{\bf t}) \\
 {\bar u_1}^{m_1},\,\ldots,\,{\bar u}_e^{m_e}
 \end{bmatrix}. \leqno{(**)}
 \]
 From $(*)$ and $(**)$ we conclude that 
 $\res{Z}\bigl[\begin{smallmatrix} \bar{a}_{\<\<{}_{X/P}}(\mu\otimes 1/{\bf t})\\ 
 {\bf \bar{u}^m} \end{smallmatrix}\bigr] = 
 \res{Z}\bigl[\begin{smallmatrix} b_{\<\<{}_{X/P}}(\mu\otimes 1/{\bf t})\\ {\bf \bar{u}^m} \end{smallmatrix}
 \bigr]$.
 Now apply local duality, i.e. \cite[Cor.\,5.1.4]{fub-abs}, 
 to conclude that 
 $\kappa^*\bar{a}_{\<\<{}_{X/P}}=\kappa^*b_{\<\<{}_{X/P}}$. This means that on a Zariski
 open neighbourhood of $Z$, $\bar{a}_{\<\<{}_{X/P}}=b_{\<\<{}_{X/P}}$. We point out that
 the hypothesis that $Z$ be defined globally by the vanishing of  a quasi-regular sequence is
 not necessary to reach this conclusion, since $j$ is a regular immersion and locally, one
 can arrange this. In other words, if we have a section of $f$, then in an open neighbourhood 
 $U$ of the image of the section, $\bar{a}_{\<\<{}_{X/P}}\vert_U=b_{\<\<{}_{X/P}}\vert_U$.
 
 In the general case, let $X''=X\times_YX$, $P'=P\times_YX$,
 and consider the cartesian square
 \[
 {\xymatrix{
 X'' \ar@{}[dr]|{\square} \ar[d]_{p_1} \ar[r]^{p_2} & X \ar[d]^f \\
 X \ar[r]_f & Y
 }}
 \] 
 We know that $p_2^*\bar{a}_{\<\<{}_{X/P}}= \bar{a}_{\<\<{}_{X''/P'}}$.
 It is clear from the description of $b_{{}_{X/P}}$ that it is compatible with arbitrary
 base change and hence  $p_2^*b_{\<\<{}_{X/P}}= b_{\<\<{}_{X''/P'}}$
 Then by what we have proven, there is a Zariski open subscheme $V$ of $X''$ containing
  the diagonal such that 
 \[p_2^*\bar{a}_{\<\<{}_{X/P}}\vert_V=p_2^*b_{\<\<{}_{X/P}}\vert_V.\]
 Let $\De\colon X\hookrightarrow V$ be the map induced by the diagonal immersion $X
 \hookrightarrow X''$. Applying $\De^*$ to both sides
of the displayed equation above, we see that $\bar{a}_{\<\<{}_{X/P}}= b_{\<\<{}_{X/P}}$.
  \qed

\begin{cor}\label{cor:res-thm} Let ${\bf v}=(v_1,\dots,v_e)\in\Gamma(P,\,\co_P)=S$,
$J$ the ideal in $S$ generated by $({\bf t}, {\bf v})$, $Z=\Spec{\,S/J}$, and
$v'_i$ the restriction of $v_i$ to $Z$ for $i=1,\,\dots, e$.  If
$Z\to Y$ is finite and flat, then
 \[
 \res{{Z, \pi}}\begin{bmatrix} 
 \wdd{t_1}{t_d}\wedge\nu \\
 t_1,\,\ldots,\,t_d, v_1,\,\ldots,\,v_e
 \end{bmatrix}
 =\res{{Z, f}}\begin{bmatrix} 
 i^*\nu \\
 v'_1,\,\ldots,\, v'_e
 \end{bmatrix}. 
 \]
 for $\nu\in \wedge^e\Omega^1_{S/A}$.
\end{cor}

\proof
\Tref{thm:res-thm} together with \cite[Prop.\,C.6.6]{fub-abs}
yields
\[
 \ush{\bar{a}}_{\<\<{}_{Z/P}}
 (\wdd{t_1}{t_d}\wedge\nu\otimes 1/({\bf t, v})) =
 \ush{\bar{a}}_{\<\<{}_{Z/X}}(i^*\nu\otimes 1/{\bf v'})
 \]
where $\ush{\bar{a}}_{\<\<{}_{Z/P}}$ and $\ush{\bar{a}}_{\<\<{}_{Z/X}}$ are as in the proof of \Tref{thm:res-thm} and ${\bf v'}$ is $(v_1',\dots, v_e')$. Let $h\colon Z\to Y$ be the
composite $Z\hookrightarrow X \xrightarrow{f} Y$. 
Applying $\tin{h}\smcirc h_*$ to both sides, we get the result. (See 
\cite[Prop.\,5.4.4]{fub-abs}.)
\qed

\setcounter{subsubsection}{\value{thm}} 
\subsubsection{Quasi-finite maps}\label{sss:q-finite} \stepcounter{thm}
 Suppose the map $\pi \colon P \to Y$ in \eqref{diag:restr-thm}
factors as $P\xrightarrow{p} W \xrightarrow{g} Y$, with $p$ smooth of relative
dimension $d$, and $g$ smooth
 of relative dimension $e$, and 
assume $h=p\smcirc \pi$ is \emph{quasi-finite}.
In other words we have
a commutative diagram of ordinary schemes
\[
{\xymatrix{
X \,\ar[d]_{h} \ar@{^(->}[r]^i & P \ar[dl]^{p} \ar[d]^\pi\\
W  \ar[r]_g & Y
}}
\]
with $h$ quasi-finite and $f=g\smcirc h=p\smcirc i$, and with $p$, $\pi$, $g$ and $f$
smooth of relative dimensions $d$, $d+e$, $e$ and $e$, respectively. 
To lighten notation, we write 
\[\eN=\eN_i^d=(\wedge^d_{\co_X}\I/\I^2)^*\]
where $\I$ is the quasi-coherent ideal sheaf in $\co_P$ defining $i\colon X\hookrightarrow P$.

Since $h$ is quasi-finite and flat
over $W$ (the latter because $p$ is smooth, and $i$ is a local complete intersection map), 
for quasi-coherent $\co_W$-module $\eF$, $h^!\eF$ can be identified with $\Hr^0(h^!\eF)$
in the standard way, and we will do so in what follows. With this convention,
we have three isomorphisms which we now describe. First, we clearly have
\stepcounter{sth}
\begin{equation*}\label{iso:h!}\tag{\thesth}
h^!\omega_g \iso \omega_f
\end{equation*}
via the isomorphism $h^!g^! \iso f^!$, and Verdier's isomorphisms for 
$f$ and $g$. 

Next, for a quasi-coherent $\co_W$-module $\eF$, we have
the transitivity isomorphism 
\stepcounter{sth}
\begin{equation*}\label{iso:chi-h}\tag{\thesth}
\chi^h(\eF, \co_W)\colon h^*\eF\otimes_{\co_X} h^!\co_W 
\xrightarrow{\>\Iso\>} h^!\eF
\end{equation*}
 of \cite[(7.2.1)]{fub-abs}. Since we are dealing
with ordinary schemes, taking account of our choice of order of tensor product,
this is the same as the map $\chi^h_{\eF, \co_W}$ of \cite[p.\,231, (4.9.1.1)]{notes}. 
The map $\chi^h(\eF,\co_W)$ is an isomorphism since $h$ is flat and hence perfect 
\cite[pp.\,234--235, Thm.\,4.9.4]{notes}.

Finally, we have an isomorphism
\stepcounter{sth}
\begin{equation*}\label{iso:p-i!}\tag{\thesth}
i^*\omega_p\otimes\eN \iso h^!\co_W
\end{equation*}
given by $\eta'_i(\omega_p[d])\colon i^\btrg(\omega_p[d])\iso i^!\omega_p$ of
\cite[(C.2.13)]{fub-abs},
the isomorphism $i^!p^!\iso h^!$, and Verdier's isomorphism 
$\bar{\bf v}\colon \omega_p[d]\iso p^!\co_W$.

These three isomorphisms are related in the following way.

\begin{sprop}\label{prop:h!} The following diagram of isomorphisms commutes
\[
{\xymatrix{
i^*(p^*\omega_g\otimes_{\co_P}\omega_p)\otimes_{\co_X}\eN \ar@{=}[d] 
\ar[r]^-{{\bar{\varphi}_{g,p}}}_-{\rotatebox{180}{\makebox[-0.1cm]{\Iso}}} 
& i^*\omega_\pi\otimes_{\co_X}\eN \ar[r]^-{\bar{a}_{X/P}}_-{\rotatebox{180}{\makebox[-0.1cm]{\Iso}}} 
 & \omega_f  \\
h^*\omega_g\otimes_{\co_X}(i^*\omega_p\otimes_{\co_X}\eN) 
\ar[d]^-{\rotatebox{90}{\makebox[0.1cm]{\Iso}}}_-{{\bf 1}\otimes\eqref{iso:p-i!}} & & \\
h^*\omega_g\otimes_{\co_X}h^!\co_W \ar[rr]^\Iso_{\eqref{iso:chi-h}} & & h^!\omega_g 
\ar[uu]^{\,\rotatebox{-90}{\makebox[-0.1cm]{\Iso}}}_-{\eqref{iso:h!}}
}}
\]
where $\bar{\varphi}_{g,\,p}$ is the explicit map described in \eqref{iso:phi-gf-bar}
and $\bar{a}_{X/P}$ is the map  \eqref{eq:a-X/P} described locally, via \Tref{thm:res-thm}, by
the explicit map $b_{X/P}$ in \eqref{eq:b-X/P}.
\end{sprop}

\proof
The essential point is that other than \eqref{iso:h!}, all other maps in the diagram
are various avatars of transitivity maps. The map $\eta'_i$ which is used
in the definition of \eqref{iso:p-i!} is 
\[\chi^i(\boldsymbol{-}, \co_P)\colon \bL i^*(\boldsymbol{-})\overset{\bL}{\otimes}_{\co_X}i^!\co_P 
\iso i^!\]
with $\eN[-d]$ substituted for $i^!\co_P$, via the canonical isomorphism 
\[\eta'(\co_P)\colon\eN[d]\iso i^!\co_P\]
(see also 
\cite[C.2.14.1]{fub-abs} for another way of looking at 
this). 

Next, according to \Tref{thm:fubini},  and the definition in
\eqref{map:zeta-gf}, the map $\bar{\varphi}_{g,\,p}$ is $\Hr^{d+e}$ of
the composite (after substituting $g^!\co_Y$, $\pi^!\co_Y$, $p^!\co_W$ with
$\omega_g[e]$, $\omega_\pi[d+e]$, and $\omega_p[d]$ respectively, via
Verdier's isomorphisms):
\[
p^*g^!\co_Y\otimes_{\co_P}p^!\co_W \iso
 p^!g^!\co_Y \iso \pi^!\co_Y,
\]
where the first arrow is the transitivity map $\chi^p(g^!\co_Y, \co_W)$, which is an isomorphism since
$p$ is flat and hence perfect.

The map $\bar{a}_{X/P}$ is, according to \eqref{map:aXP} and \eqref{def:aXP}
(after the usual Verdier substitutions and the substitution $\eN[-d]\iso i^!\co_P$), $\Hr^d$
applied to the composite
\[
\bL i^*\pi^!\co_Y \overset{\bL}{\otimes}_{\co_X}i^!\co_P \iso i^!\pi^!\co_Y \iso f^!\co_Y
\]
where the first arrow is the transitivity map $\chi^i(\pi^!\co_Y, \co_P)$. This is an isomorphism
since a regular immersion is a perfect map.

Finally, \eqref{iso:h!} is by definition $\chi^h(\omega_g, \co_W)$.

Consider the diagram below, in which the arrows are either natural ones
arising from the pseudofunctorial nature of $\boldsymbol{-^!}$ or from abstract
transitivity maps, and in which:
\[i^*=\bL i^* \, {\text{and $\otimes=\overset{\bL}{\otimes}$}}.\] 

\[
{\xymatrix{
 i^*(p^*g^!\co_Y\otimes p^!\co_W)\otimes i^!\co_P
\ar[r]^-\Iso \ar@{=}[d] &  i^*(p^!g^!\co_Y)\otimes i^!\co_P \ar[r]^-\Iso  & i^*\pi^!\co_Y\otimes i^!\co_P 
\ar[d]^-{\rotatebox{90}{\makebox[0.1cm]{\Iso}}} \\
h^*(g^!\co_Y)\otimes (i^*p^!\co_W\otimes i^!\co_P) \ar[dd]_{\,\rotatebox{-90}{\makebox[-0.1cm]{\Iso}}} 
&& i^!\pi^!\co_Y \ar[d]^--{\rotatebox{90}{\makebox[0.1cm]{\Iso}}}\\
&& f^!\co_Y \\
h^*(g^!\co_Y)\otimes i^!p^!\co_W \ar[r]^\Iso & h^*(g^!\co_Y) \otimes h^!\co_W \ar[r]^-\Iso 
& h^!g^!\co_Y \ar[u]_{\rotatebox{90}{\makebox[0.1cm]{\Iso}}}
}}
\]
The diagram commutes by
\cite[Prop.-Def.\,7.2.4\,(ii)]{fub-abs} and 
\cite[p.\,238]{notes}.  The Proposition follows.
\qed

\section{\bf Traces of differential forms for finite maps}
\subsection{Tate traces} Let $A$ be a ring, and $C$ an $A$-algebra
which is finite and free as an $A$-module. We have the {\emph{canonical trace}}
\stepcounter{thm}
\begin{equation*}\label{map:can-tr}\tag{\thethm}
\trc{C/A}\colon C\to A
\end{equation*}
given by the composite
\[C \lra \mathrm{End}_A(C,\,C) \lra A\]
where the first arrow is the map $c\mapsto (x\mapsto cx)$ and the second the standard
trace of an endomorphism of a finite free $A$-module.  

If the $C$-module $\Hom_A(C,\,A)$ is a free
$C$-module of rank one (this happens if and only if, in addition to $C$ being a finite free $A$-module,
its fibres are Gorenstein) 
then, following Kunz in \cite{kd}, we regard any free generator of 
$\Hom_A(C,\,A)$ as a ``trace" for the $A$-algebra $C$ (cf.~\cite[F8\,(b), pp.\,362--363]{kd}).
If there is one, then clearly we have exactly as many as the units of $C$. We point
out that the canonical trace, $\trc{C/A}$, need not be a trace in this sense on $C$. 
Indeed if $A$ and $C$ are
fields and $C$ is a purely inseparable extension of $A$, then $\trc{C/A}=0$ and hence cannot be 
a free generator of $\Hom_A(C,A)$.

Tate studies the existence and characterisation of traces in an important situation which
includes the case of $C$ being a complete intersection algebra over $A$.

In the rest of this sub-section we make the following assumptions and use the following
notations. The $A$-algebra $C$ (which is free of finite rank as an $A$-module) is such that the
 canonical map $A\to C$ factors as 
\[A\lra B \xrightarrow{\pi} C\]
with $\pi$ a surjective map, the kernel $I$ of $\pi$ generated by a regular $B$-sequence
${\bf f}=(f_1,\dots, f_n)$, and the kernel $J$ of the canonical map
\[s\colon B\otimes_AC \lra C\]
is generated by a $B\otimes_AC$-sequence ${\bf g}=(g_1,\dots, g_n)$. In somewhat greater
detail, if $m\colon C\otimes_AC\to C$ is the $A$-algebra map $c\otimes c'\mapsto cc'$, then 
$s$ is the composition
\[B\otimes_AC \xrightarrow{\phantom{X}\pi\otimes {\bf 1}_C\phantom{X}} C\otimes_AC 
\xrightarrow{\phantom{X}m\phantom{X}} C.\]
Note that $f_i\otimes 1\in J$ and hence we have $h_{ij}\in B\otimes_AC$ such that
$f_i\otimes 1 = \sum_{j=1}^n h_{ij}g_j$ for $i=1,\dots, n$.
Let
\stepcounter{thm}
\begin{equation*}\label{def:det-gf}\tag{\thethm}
\De = \det{(h_{ij})},
\end{equation*}
and
\stepcounter{thm}
\begin{equation*}\label{def:det-bar}\tag{\thethm}
\overline{\De} = (\pi\otimes {\bf 1}_C)(\De).
\end{equation*}
Set $\overline{J}=\ker{m}= (\pi\otimes {\bf 1}_C)(J)$. We have the following commutative diagram
with 
\[\check{s}\colon C\otimes_AB \to C\]
being the composite $m\smcirc (1_C\otimes\pi)$.
\stepcounter{thm}
\[\begin{aligned}\label{diag:ABC}
{\xymatrix{
& & & & &C\\ 
C \ar[rr] && C\otimes_AB \ar[rr]^{1_C\otimes\pi} \ar@/^1.5pc/[urrr]^{s^{\<{}_{\vee}}}&& C\otimes C \ar[ur]^{m} & \\
B \ar[u]^{\pi}\ar[rr] && B\otimes_AB \ar[u]^{\pi\otimes 1_B} \ar[rr]_{1_B\otimes\pi} 
&& B\otimes_AC \ar[u]^{\pi\otimes 1_C}  \ar@/_1.5pc/[uur]_{s}& \\
A \ar[u] \ar[rr] & &B \ar[u] \ar[rr]_{\pi} & & C \ar[u] &
}}
\end{aligned}\tag{\thethm}
\]

In the above situation it is shown in \cite[Appendix]{tate} that traces exist (i.e., $\Hom_A(C,\,A)$ is
a rank one free $C$-module)
and there is a canonical free generator  (i.e., a trace)
$\lambda=\lambda({\bf f}, {\bf g})$ of $\Hom_A(C,\,A)$. We summarise the results of Tate
as given in \cite[Appendix]{tate} in the following two theorems in which we make the standard
indentifications $B\otimes_A\Hom_A(C,\,A)=\Hom_B(B\otimes_AC,\,B)$
and $C\otimes_A\Hom_A(C,\,A)=\Hom_C(C\otimes_AC,\,C)$. Under these identifications it is clear
that
\stepcounter{thm}
\begin{equation*}\label{eq:phi-pi-d}\tag{\thethm}
\pi\smcirc(1_B\otimes\phi) = (1_C\otimes\phi)\smcirc(\pi\otimes 1_C) \qquad (\phi\in\Hom_A(C,\,A)).
\end{equation*}

\begin{thm}\label{thm:tate1} 
{\emph{(Tate) \cite[p.231, Lemma\,(A.10)]{tate}}} The map 
\[t\colon \Hom_A(C,\,A) \lra C\]
given by 
\[\phi \mapsto \pi((1_B\otimes\phi)(\De))=(1_C\otimes\phi)(\overline{\De})\]
is an isomorphism of $C$-modules.
\end{thm}

In \textit{loc.cit.}~the description of $t$ is $\phi \mapsto \pi((1_B\otimes\phi)(\De))$. Using
\eqref{eq:phi-pi-d} it is clear that $t$ can also be described as 
$\phi \mapsto (1_C\otimes\phi)(\overline{\De})$.

The results in \cite[Appendix]{tate} are perhaps more useful when stated in the following way.

\begin{thm}\label{thm:tate2} {\emph{(Tate)}} Let $\lambda=\lambda({\bf f},\,{\bf g})$ 
be the free $C$-module generator of $\Hom_A(C,\,A)$ given by
\[\lambda=t^{-1}(1)\]
where $t\colon \Hom_A(C,\,A)\iso C$ is the isomorphism in \Tref{thm:tate1}.
\begin{enumerate}
\item[(a)] If $\phi\in \Hom_A(C,\,A)$, then the constant of proportionality $c\in C$
such that $\phi=c\lambda$, is given by 
\[c=\pi((1_B\otimes\phi)(\De) = (1_C\otimes \phi)(\overline{\De}).\]
\item[(b)] If $\psi\in \Hom_B(B\otimes_AC,\,B)$ and $\phi\in\Hom_A(C,\,A)$
are such that $1_B\otimes\phi -\psi\in J\Hom_B(B\otimes_AC,\,B)$, then
\[\pi((1_B\otimes\phi)(\De))=\pi\psi(\De).\]
\item[(c)]  If $\psi\in \Hom_C(C\otimes_AC,\,C)$  and $\phi\in\Hom_A(C,\,A)$ are
such that $1_C\otimes\phi -\psi\in \overline{J}\Hom_C(C\otimes_AC,\,C)$, then
\[(1_C\otimes\phi)(\overline{\De}))=\psi(\overline{\De}).\]
\item[(d)] If $\trc{C/A}\colon C\to A$ is the canonical trace given in \eqref{map:can-tr}, then
\[\trc{C/A}=m(\overline{\De})\lambda.\] 
\end{enumerate}
\end{thm}

\proof These are all results in \cite[Appendix]{tate}, stated in perhaps a different way.
Part\,(a) is [\textit{ibid},\,pp.\,229--230,\, 3.\,of Theorem\,(A.3)] (together with \eqref{eq:phi-pi-d}).
Part (b) is an immediate consequence of  [\textit{ibid},\,p.\,230,\, Lemma\,(A.9)] and (c) is the same, together with \eqref{eq:phi-pi-d}.
Part\,(d) is  [\textit{ibid},\,pp.\,229--230,\, 4.\,of Theorem\,(A.3)].
\qed

The first application of Tate's result we give is the following (this is (R6) of \cite[p.\,198]{RD} but
for our version of residues).

\begin{thm}\label{thm:R6} In the above situation, suppose $B$ is smooth of relative dimension $n$ 
over $A$, $f\colon X\to Y$ the corresponding smooth map from $X=\Spec{\,B}$ to $Y=\Spec{\,A}$,
and $Z=\Spec{\,C}$. Then
\[\res{{Z,f}}\begin{bmatrix} b\,\wdd{f_1}{f_n}\\
f_1, \dots,f_n
\end{bmatrix}
= \trc{C/A}(b\vert_Z).
\]
\end{thm}
\proof It is important to keep diagram\,\eqref{diag:ABC} in mind when following this proof.
There is an annoying issue that $\De$ is defined in terms of $f_i\otimes 1_C$ and $g_i$,
but in dealing with the base change $1_C\otimes\phi$, for $\phi\in\Hom_A(C,\,A)$, the
natural elements that show up are $1_C\otimes f_i\in C\otimes_AB$. One has to do somewhat
careful book-keeping to avoid confusion.
Since $C\otimes_AB$ and $B\otimes_AC$ play different roles, 
let us agree to write $x^{{}_\vee}$ for the element of $C\otimes_AB$
corresponding to $x\in B\otimes_AC$ under the standard isomorphism between
$B\otimes_AC$ and $C\otimes_AB$.

In what follows, the $C$-algebra structures on $C\otimes_AB$ and $C\otimes_AC$ are
$c\mapsto c\otimes 1_B$ and $c\mapsto c\otimes 1_C$ respectively. Let
$N=(I/I^2)^*$ and $N_C=C\otimes_AN$. Let $h\colon Z\to Y$ be the
natural finite flat map corresponding to $A\to C$ and $i\colon Z\hookrightarrow X$
the natural closed immersion, with normal bundle $\eN$.
If $\vttr{C/A}\colon
\Omega^n_{B/A}\otimes_B\wedge^n_CN\to A$ 
is the map arising from the composite (all isomorphisms
being the obvious ones, e.g., the fundamental local isomorphism,  Verdier's isomorphism, \dots)
\[h_*(i^*(\Omega^n_{X/Y})\otimes_Z\wedge^n\eN) \iso \Hr^0(h_*(i^!f^!\co_Y)) \iso \Hr^0(h_*h^!\co_Y)
\xrightarrow{\Hr^0(\Tr{h})} \co_Y\]
then $(\Omega^n_{B/A}\otimes_B\wedge^n_CN, \vttr{C/A})$ represents the functor
$M\mapsto \Hom_A(M,\,A)$ from
finite $C$-modules to finite $A$-modules, whence we have an isomorphism of $C$-modules
\[\Phi\colon \Omega^n_{B/A}\otimes_C\wedge^n_CN \iso \Hom_A(C,\,A)\]
with $\vttr{C/A}$ corresponding to ``evaluation at 1" under this isomorphism. According
to \cite[Prop.\,5.4.4]{fub-abs}, we have
\[\vttr{C/A}(\mu\otimes {\bf 1/f})= \res{Z}\begin{bmatrix} \mu \\ f_1, \dots f_n\end{bmatrix}
\qquad (\mu\in \Omega^n_{B/A}).\] 
Thus
\[\Phi(\mu\otimes {\bf 1/f})(c) = \res{Z}\begin{bmatrix} b\cdot \mu \\ f_1, \dots, f_n\end{bmatrix} 
\qquad (c\in C)\]
where $b\in B$ is any pre-image of $c$. If $b\in I$, then $b \, \wdd{f_1}{f_n}\otimes{1/{\bf f}}=0$
in $\Omega^n_{B/A}\otimes_B \wedge^nN$ and hence the right side of the
above displayed formula is well-defined as a function of $c\in C$.

Similarly we have an isomorphism of $C\otimes_AC$-modules
\[\Phi'\colon \Omega^n_{(C\otimes_AB)/C}\otimes \wedge^n N_C \iso \Hom_C(C\otimes_AC, C)\]
given by
\[\Phi'(\nu\otimes {\bf 1/}(1_C\otimes {\bf f}))(x) = 
\res{Z\times_YZ, p}\begin{bmatrix} \wit{x}\cdot \mu \\ (1_C\otimes f_1), \dots, (1_C\otimes f_n)\end{bmatrix} 
\qquad (x\in C\otimes_AC)\]
where $\wit{x}\in C\otimes_AB$ is any pre-image of $x$ and $p\colon Z\times_YX\to Z$
is the natural projection.

Let $s^{{}_\vee}\colon C\otimes_AB \to C$ be as in \eqref{diag:ABC}, i.e., 
$s^{{}_\vee}=m\smcirc (1_C\otimes\pi)$. Then $J^{{}_\vee}\set \ker{s^{{}_\vee}}$ is generated by
$g^{{}_\vee}_1, \dots, g^{{}_\vee}_n$.

Let $\phi\in \Hom_A(C,\,A)$ and $\psi\in \Hom_C(C\otimes_AC,\,C)$ be the maps defined by
\[\phi =\Phi(\wdd{f_1}{f_n}\otimes {\bf 1/f}),\]
and 
\[\psi = \Phi'((\De^\vee\cdot \wdd{g_1^{{}_\vee}}{g_n^{{}_\vee}})\otimes {\bf 1/}(1_C\otimes {\bf f})).\]
We have to show that $\phi=\trc{C/A}$. By \Tref{thm:tate2}\,(d), this is equivalent to
showing that $(1_C\otimes\phi)(\overline{\De})= m(\overline{\De})$. It easier
to show that $\psi(\overline{\De})=m(\overline{\De})$, and we can reduce to this
via \Tref{thm:tate2}\,(c). The details are as follows.
First, we claim that 
$1_C\otimes\phi -\psi \in \overline{J}\Hom_C(C\otimes_AC,\,C)$ so that \Tref{thm:tate2}\,(c) applies.
Before we prove the claim, we point out that
\[1_C\otimes \phi = \Phi'((\wdd{(1_C\otimes f_1)}{1_C\otimes f_n})\otimes {\bf 1/}(1_C\otimes {\bf f}).\]
Since 
$1_C\otimes f_i=\sum_j h_{ij}^{{}_\vee}g_j^{{}_\vee}$ we have 
\[\wdd{(1_C\otimes f_1)}{(1_C\otimes f_n)} = \mu + {\De}^{\vee}\,\wdd{g^{{}_\vee}_1}{g^{{}_\vee}_n}\]
where $\mu\in J^{\vee}\Omega^n_{(C\otimes_AB)/C}$ (for $h_{ij}^{{}_\vee}\in J^\vee$). It follows
that 
\[1_C\otimes\phi -\psi = \Phi'(\mu\otimes ({\bf 1/}(1_C\otimes {\bf f})))\in 
\overline{J}\Hom_C(C\otimes_AC,\,C)\] 
as claimed.

We then have, with $\delta\in C\otimes_AB$ a lift of $\overline{\De}\in C\otimes_AC$,
\[
\begin{aligned}
(1_C\otimes\phi)(\overline{\De})  = \psi(\overline{\De})
& = \res{{Z\times_YZ, p}}
\begin{bmatrix}\delta\,\De^\vee\,\wdd{g^{{}_\vee}_1}{g^{{}_\vee}_n}\\
(1_C\otimes f_1),\dots, (1_C\otimes f_n) \end{bmatrix} \\
& = \res{{Z,p}}
\begin{bmatrix}\delta\,\wdd{g^{{}_\vee}_1}{g^{{}_\vee}_n}\\
g^{{}_\vee}_1,\dots, g^{{}_\vee}_n \end{bmatrix}\\
& = s^{{}_\vee}(\delta) \\
&= m(\overline{\De}).
\end{aligned}\tag{{\text{$*$}}}
\]
In the above sequence, the first equality is from \Tref{thm:tate2}\,(c), the one in the second
line from \cite[Thm.\,5.4.5]{fub-abs}, the third from the
fact that the composite $Z\xrightarrow{\text{via $s^{{}_\vee}$}} Z\times_YX\xrightarrow{p} Z$ is
an isomorphism, which means the formulae in \Rref{rem:dirac} apply. The last equality
is from the definition of $s^{{}_\vee}$ as $m\smcirc (1_C\otimes\pi)$. From $(*)$ and 
\Tref{thm:tate2}\,(d) we get that $\phi=\trc{C/A}$, and from this the Theorem follows.
\qed

\begin{rems}\label{rems:R6}\label{rems:lambda'} {\emph{ 1) The above proof would be easier if one could show that $\De^\vee$
is a pre-image of $\overline{\De}$ under $C\otimes_AB \xrightarrow{1_C\otimes\pi} C\otimes_AC$.
But there is no guarantee it is so. However, in the special case where $B$ is a polynomial
ring over $A$, something like this be arranged as the proof \Pref{prop:lambda-res} below shows.}}

{\emph{2) If $B$ is flat over $A$, then $\lambda$ is stable under any base change of $A$. In 
somewhat greater detail, if $A\to A'$ is a map of rings, $B'$, $C'$, ${\bf f}'$, and ${\bf g}'$ the obvious base changes of $B$, $C$, ${\bf f}$, and ${\bf g}$, then, under the identification
$\Hom_{A'}(C'\,A')=A'\otimes_A\Hom_A(C,\,A)$, we have 
$\lambda({\bf f}',\,{\bf g}' )=1\otimes\lambda({\bf f},\,{\bf g})$. This is because, if $B$ is flat over
$A$, then ${\bf f}'$ and ${\bf g}'$ are regular sequences.
}}
\end{rems}

\begin{prop}\label{prop:lambda-res} Let $q\in A[T_1, \dots, T_n]=A[{\bf T}]$.
Suppose $B$ is the $A$-algebra 
$B=A[{\bf T}]_q$.
For $i=1,\dots,n$ let $\gamma_i=\pi(T_i)$ and 
\[g_i=T_i\otimes 1_C - 1_B\otimes\gamma_i.\]
Let $Z=\Spec{\,C}$. Then $\lambda=\lambda({\bf f},\,{\bf g})\in\Hom_A(C,\,A)$ is given by
\[\lambda(c) = \res{Z}\begin{bmatrix} b\,\wdd{T_1}{T_n}\\ f_1,\dots, f_n\end{bmatrix} \qquad(c\in C)\]
where $b\in B=A[T_1,\dots, T_n]$ is any pre-image of $c$.
\end{prop}

\proof It is straightforward to see that the $g_i$, as defined in the Proposition, 
generate $J=\ker{s}$, and form a regular $B\otimes_AC$-sequence. As before, let
$X=\Spec{\,B}$, $Y=\Spec{\,A}$, $Z=\Spec{\,C}$, and let $p\colon Z\times_YX\to Z$
be the projection map. As we did earlier, we need to distinguish between $B\otimes_AC$ and
$C\otimes_AB$, and so between $X\times_YZ$ and $Z\times_YX$, and $p$ corresponds
to the map $C\to C\otimes_AB$ given by $c\mapsto c\otimes 1$.

For the proof of the theorem, it is simpler to regard the two copies of $B$ in Diagram\,
\eqref{diag:ABC}, the one in the
middle of the bottom row, and the one in the middle of the left column, as two different copies
of $A[{\bf T}]_q$, say $A[X_1,\dots, X_n]_{q(X_1,\dots, X_n)}=A[{\bf X}]_{q({\bf X})}$ 
and $A[Y_1,\dots, Y_n]_{q(Y_1,\dots,Y_n)}=A[{\bf Y}]_{q({\bf Y})}$ respectively. 
Then $B\otimes_AB$ can be regarded as $A[{\bf X},\,{\bf Y}]_{q({\bf X})q({\bf Y})}$. 
Moreover, $B\otimes_AC$ is then
identified with $C[{\bf Y}]_{q({\bf Y})}$ and $C\otimes_AB$ with $C[{\bf X}]_{q({\bf X})}$.
Diagram\,\eqref{diag:ABC} translates to
\[\begin{aligned}
{\xymatrix{
& & & & &C\\ 
C \ar[rr] && C[{\bf X}]_{q({\bf X})} \ar[rr]^{\pi_{{}_1}''} 
\ar@/^1.5pc/[urrr]^{s^{\<{}_\vee}}&& C\otimes C \ar[ur]^{m} & \\
A[{\bf Y}]_{q({\bf Y})} \ar[u]^{\pi_{{}_2}}\ar[rr] && A[{\bf X},\,{\bf Y}]_{q({\bf X})q({\bf Y})} \ar[u]^{\pi'_{{}_2}} \ar[rr]_{\pi'_{{}_1}} 
&& C[{\bf Y}]_{q({\bf Y})} \ar[u]^{\pi''_{{}_2}}  \ar@/_1.5pc/[uur]_{s}& \\
A \ar[u] \ar[rr] & & A[{\bf X}]_{q({\bf X})} \ar[u] \ar[rr]_{\pi_{{}_1}} & & C \ar[u] &
}}
\end{aligned}
\]
Here $\pi_{{}_1}$ is the map $X_i\mapsto \gamma_i$, and $\pi_{{}_2}$ is $Y_i\mapsto \gamma_i$.
The maps $\pi'_{{}_1}$ and $\pi''_{{}_1}$ are the base changes of $\pi_{{}_1}$, and 
$\pi'_{{}_2}$, $\pi''_{{}_2}$ the base changes of $\pi_{{}_2}$. We point out that
\[\pi''_{{}_1}\Biggl(\sum_{\underline{i}} c_{\underline{i}}{\bf X}^{\underline{i}}\Biggr) = 
\sum_{\underline{i}} c_{\underline{i}}\otimes{\boldsymbol{\gamma}}^{\underline{i}}\]
and
\[\pi''_{{}_2}\Biggl(\sum_{\underline{i}} c_{\underline{i}}{\bf Y}^{\underline{i}}\Biggr) = 
\sum_{\underline{i}} {\boldsymbol{\gamma}}^{\underline{i}}\otimes c_{\underline{i}}.\]
For any $h\in B=A[{\bf T}]_q$, the element $h\otimes 1_B$ (resp.~$h\otimes 1_C$) is identified
with the element  $h({\bf Y})$ of $A[{\bf X},\,{\bf Y}]_{q({\bf X})q({\bf Y})}=B\otimes_AB$ (resp.~the 
element $h({\bf Y})$ of $C[{\bf Y}]_{q({\bf Y})}=B\otimes_AC$),
whereas $1_B\otimes h$ and $1_C\otimes h$ are identified with
$h({\bf X})$ (regarded as elements of $A[{\bf X},\,{\bf Y}]_{q({\bf X})q({\bf Y})}$ and of 
$C[{\bf X}]_{q({\bf X})}$ respectively).
Finally $s(\sum_{\underline{i}} c_{\underline{i}}{\bf Y}^{\underline{i}})= 
\sum_{\underline{i}} c_{\underline{i}}{\boldsymbol{\gamma}}^{\underline{i}}$ and
$s^{{}_{\vee}}(\sum_{\underline{i}} c_{\underline{i}}{\bf X}^{\underline{i}})=
\sum_{\underline{i}} c_{\underline{i}}{\boldsymbol{\gamma}}^{\underline{i}}$.
It follows that
\[g_i= Y_i-\gamma_i, \quad {\text{and}} \quad g_i^{{}_\vee} = X_i-\gamma_i \qquad (i=1,\dots n).\]
Now there exist $h_{ij}({\bf X}, {\bf Y})\in A[{\bf X},\,{\bf Y}]$ such that
\[
f_i({\bf X}) - f_i({\bf Y}) = \sum_j h_{ij}({\bf X}, {\bf Y})(X_j-Y_j).
\]
Then $f_i({\bf Y}) = \sum_{j} h_{ij}({\bf {\boldsymbol{\gamma}}}, {\bf Y})(Y_j-\gamma_j)$
and
$f_i({\bf X}) = \sum_{j} h_{ij}({\bf X}, {\boldsymbol{\gamma}})(X_j-\gamma_j)$.
Let
\[\delta({\bf X},\,{\bf Y})=\det{(h_{ij}({\bf X}, {\bf Y}))}.\]
If $\De$ is defined as in \eqref{def:det-gf}, then 
\[\De=\delta({\boldsymbol{\gamma}},\,{\bf Y}). \]
Note that 
\[\overline{\De} = \pi''_{{}_2}(\De)= \pi''_{{}_1}(\delta({\bf X},\,{\boldsymbol{\gamma}})).\leqno{(*)}\]
On the other hand, since 
$f_i({\bf X}) = \sum_{i,j} h_{ij}({\bf X}, {\boldsymbol{\gamma}})(X_i-\gamma_i)$,
according to \cite[Thm.\,5.4.5]{fub-abs} we have
\[\res{{Z\times_YZ,\,p}}\begin{bmatrix}\delta({\bf X},\,{\boldsymbol{\gamma}}) \mu\\ 
f_1({\bf X}),\dots, f_n({\bf X})\end{bmatrix}=
\res{Z}\begin{bmatrix} \mu\\ g_1^{{}_\vee},\dots, g_n^{{}_\vee}\end{bmatrix} 
\qquad \Bigl(\mu\in \Omega^n_{C[{\bf X}]/C}\Bigr). \leqno{(**)}\]
Let $\phi\colon C\to A$ be defined by
\[\phi(c)=\res{Z}\begin{bmatrix} b\,\wdd{T_1}{T_n}\\ f_1,\dots, f_n\end{bmatrix} \qquad (c\in C)\]
where $b\in B=A[{\bf T}]$ is any element in $\pi^{-1}(c)$. Since 
$\mathrm{d}X_i=\mathrm{d}(X_i-\gamma_i)$,  therefore for $x\in C\otimes_AC$ and
$\wit{x}\in C[{\bf X}]$ such that $\pi''_{{}_1}(\wit{x})=x$, we have
\[
(1_C\otimes\phi)(x)=\res{{Z\times_YZ,p}}\begin{bmatrix}
\wit{x}\,\,\wdd{(X_1-\gamma_1)}{(X_n-\gamma_n)}\\
f_1({\bf X}),\dots, f_n({\bf X})
\end{bmatrix}.
\]
By $(*)$ and $(**)$ we get
\[(1_C\otimes\phi)(\overline{\De})=\res{Z}\begin{bmatrix} \wdd{g_1^{{}_\vee}}{g_n^{{}_\vee}}\\
g_1^{{}_\vee},\dots, g_n^{{}_\vee}
\end{bmatrix} = 1.
\]
\Tref{thm:tate2}\,(a) then gives $\phi=\lambda$.
\qed

\subsection{Traces of differential forms} Suppose we have a commutative diagram of 
ordinary schemes
\stepcounter{thm}
\[
\begin{aligned}\label{diag:f=gh}
{\xymatrix{
X \ar[d]_h \ar@/^1.0pc/[rrd]^f& & \\
Y \ar[rr]^g && Z
}}
\end{aligned}\tag{\thethm}
\]
with $f$ and $g$ smooth of relative dimension $n$ and $h$ a finite map (necessarily flat).
The composite $h_*f^!\co_Z \iso h_*h^!g^!\co_Z \xrightarrow{\Tr{h}} g^!\co_Z$ gives
a $\co_Y$-map map (after applying Verdier's isomorphism to $f^!\co_Z$ and $g^!\co_Z$ and
applying $\Hr^{-n}({\boldsymbol{-}})$)
\stepcounter{thm}
\begin{equation*}\label{map:tr-h-fin}\tag{\thethm}
\vin{h}\colon h_*\omega_f \lra \omega_g.
\end{equation*}
A note of caution. We have used the symbol $\vin{p}$ earlier for the trace map 
$\Rr^mp_*\omega_f \to \co_W$ for a smooth proper map $p\colon V\to W$ of
relative dimension $m$. The context will make the meaning of the symbol clear. 

\begin{prop}\label{prop:R10} Let $W$ be a closed subscheme of $Y$, proper over $Z$,
and let $W'=h^{-1}(W)$. Assume $g$ (and hence $f$) is separated.
Then the following diagram commutes:
\[
{\xymatrix{
\Rr^n_{W'}f_*\omega_f \ar[d]_{\res{{W'}}}  &&
\Rr^n_{W}g_* h_*\omega_f \ar[ll]_{\Iso} \ar[d]^{\mathrm{via}\,\, \vin{h}} \\
\co_Z & & \Rr^n_Wg_*\omega_g \ar[ll]^{\res{W}}  
}}
\]
\end{prop}
\proof
By Nagata's compactification \cite{nagata} we have an open immersion $u\colon Y\to \overline{Y}$
together with a proper map $\bar{g}\colon \overline{Y} \to Z$  such that $\bar{g}\smcirc u = g$.
By Zariski's Main Theorem the quasi-finite map $u\smcirc h\colon X\to \overline{Y}$ can be
completed to a finite map, i.e., we can find an open immersion $v\colon X\to\overline{X}$
and a finite map $\bar{h}\colon \overline{X}\to \overline{Y}$ such that $u\smcirc h = \bar{h}\smcirc v$.
Moreover, we may assume $X$ is scheme-theoretically dense in $\overline{X}$ so that
$\bar{h}^{-1}(u(Y)) = v(X)$. Let $\bar{f}= \bar{g}\smcirc\bar{h}$. We have a
composite
\[\bar{h}_*{\bar f}^! \iso \bar{h}_*\bar{h}^!\bar{g}^! \xrightarrow{\Tr{\bar{h}}} \bar{g}^!.\tag{\dag}\]

Consider the
commutative diagram
\[
{\xymatrix{
\R g_* \R\iG{W}\omega_g[n] \ar[r]^{\Iso}
& \R \bar{g}_* \R\iG{u(W)}\bar{g}^!\co_Z \ar[r]  
& \R \bar{g}_* \bar{g}^!\co_Z \ar[r]^{\Tr{\bar{g}}} 
&  \co_Z \ar@{=}[ddd] 
\\
 \R g_*\R\iG{W}h_*\omega_f[n] \ar[r]^{\Iso} \ar[u]^{\vin{h}} 
 \ar[d]_{\,\rotatebox{-90}{\makebox[-0.1cm]{\Iso}}}
& \R \bar{g}_*\R\iG{u(W)}\bar{h}_*\bar{f}^!\co_Z \ar[r] \ar[u]^{(\dag)}
 \ar[d]_{\,\rotatebox{-90}{\makebox[-0.1cm]{\Iso}}}
& \R \bar{g}_*\bar{h}_*\bar{f}^!\co_Z  \ar[u]^{(\dag)} \ar@{=}[d]
&
\\
\R g_*\R h_* \R\iG{W'}\omega_f[n] \ar[d]_{\,\rotatebox{-90}{\makebox[-0.1cm]{\Iso}}} \ar[r]^{\Iso} 
&  \R \bar{g}_*\R \bar{h}_* \R\iG{v(W')}\bar{f}^!\co_Z
 \ar[d]_{\,\rotatebox{-90}{\makebox[-0.1cm]{\Iso}}} \ar[r] 
& \R \bar{g}_*\R \bar{h}_*\bar{f}^!\co_Z \ar[d]_{\,\rotatebox{-90}{\makebox[-0.1cm]{\Iso}}}
&
\\
 \Rfs\R\iG{W'}\omega_f[n] \ar[r]^{\Iso}
 & \R\bar{f}_*\R\iG{v(W')}\bar{f}^!\co_Z \ar[r]
 &\R\bar{f}_*\bar{f}^!\co_Z \ar[r]_{\Tr{\bar{f}}}
 &\co_Z
}}
\]
The rectangle on the right commutes by definition of $(\dag)$ (especially of the isomorphism
$\bar{h}^!\bar{g}^! \iso \bar{f}^!$ which drives $(\dag)$).

Applying $\Hr^0(\boldsymbol{-})$ to the above diagram we get the asserted result.
\qed 
\medskip

\begin{prop}\label{prop:kunz-sigma} Let $f$, $g$, $h$ be as above, and
suppose $u\colon Z'\to Z$ is a map of ordinary
schemes. Let
\[
{\xymatrix{
X' \ar[r]^w \ar[d]_{h'} \ar@{}[dr]|\square \ar@/_2.5pc/[dd]_{f'}
 & X \ar[d]^h \ar@/^2.5pc/[dd]^{f}\\
Y' \ar[r]^v \ar[d]_{g'} \ar@{}[dr]|\square & Y \ar[d]^g \\
Z' \ar[r]_u & Z
}}
\]
be the corresponding base change diagram.
Then $v^*\vin{h}=\vin{h'}$.
\end{prop}
\proof
By \cite[${\rm{IV}}_3$, (13.3.2)]{ega}, $Y$ can be covered by open subschemes $U$ such
that $U\to Y$ is the composite of a quasi-finite map $U\to {\mathbb P}^n_Z$ followed
by the structural map ${\mathbb P}^n_Z\to Z$. Since the question is local on $Y$, we replace
$Y$ by $U$ if necessary, and assume we have a quasi-finite map $Y\to {\mathbb P}^n_Z$.
Using Zariski's Main Theorem we can
find a finite map $\overline{Y} \to {\mathbb P}^n_Z$ such that $Y$ is an open 
${\mathbb P}^n_Z$-subscheme of $\overline{Y}$. 

Since $h$ is finite, the composite
$X \to Y \hookrightarrow \overline{Y}$ is quasi-finite, and another application of
Zariski's Main Theorem tells us that $X \to \overline{Y}$ factors as an open
immersion $X \hookrightarrow \overline{X}$ followed by a finite map 
$\bar{h}\colon \overline{X} \to \overline{Y}$. Replacing $\overline{X}$ by the scheme
theoretic closure of its open subscheme $X$ if necessary, we may assume that $X$ is
scheme theoretically dense in $\overline{X}$. This forces $X=\bar{h}^{-1}(Y)$. We thus
have a cartesian diagram, with horizontal arrows being open immersions
\[
{\xymatrix{
X\>\> \ar[d]_h \ar@{^(->}[r] \ar@{}[rd]|\square& \overline{X} \ar[d]^{\bar h} \\
Y\>\> \ar@{^(->}[r] & \overline{Y}
}}
\]
We write $\bar{g}\colon \overline{Y}\to Z$ for the composite $\overline{Y}\to {\mathbb P}^n_Z\to Z$,
and set $\bar{f}=\bar{g}\smcirc\bar{h}$. 
The important point is that $\bar{f}\colon \overline{X}\to Z$ and $\bar{g}\colon \overline{Y}\to Z$ are proper over $Z$ and 
\textit{the fibres of $\bar{f}$ and $\bar{g}$ have dimension $\le n$}.
This means $\Hr^j(\bar{f}^!\co_Z)=\Hr^j(\bar{g}^!\co_Z)=0$ for
$j< -n$. It follows that if $\omgs{\bar{f}}\set \Hr^{-n}(\bar{f}^!\co_Z)$, and if
$\tin{\!\!\!\bar{f}}\colon \Rr^n\bar{f}_*\omgs{\bar{f}}\to \co_Z$
is the map induced by $\Tr{\bar{f}}(\co_Z)\colon \R{\bar f}_* \bar{f}^!\co_Z \to\co_Z$ 
then $(\omgs{\bar{f}}, \tin{\!\!\!\bar{f}})$ represents
the functor $\eF \mapsto \Hom_Z(\Rr^n{\bar{f}}_*\eF,\,\co_Z)$ of quasi-coherent sheaves $\eF$
on $X$ (see \cite[(5.1.5)]{fub-abs} for
this argument). Along these lines, if $\omgs{\bar{g}} \set \Hr^{-n}(\bar{g}^!\co_Z)$, and 
$\tin{\!\!\!\bar{g}}\colon \Rr^n\bar{g}_*\omgs{\bar{g}}\to \co_Z$, the map induced
by $\Tr{\bar{g}}(\co_Z)$, then one can make a similar statement about 
$(\omgs{\bar{g}}, \tin{\!\!\!\bar{g}})$.

Let $\overline{X}'=X\times_ZZ'$, $\overline{Y}'=\overline{Y}\times_ZZ'$, and let
$\bar{f}'$, $\bar{g}'$, $\bar{h}'$, $\bar{u}$, $\bar{v}$ be the obvious base changes of
$f$, $g$, $h$, $u$, and $v$, respectively. Let $\tin{\!\!\!h}\colon h_*\omgs{f}\to \omgs{g}$ be
the obvious analogue of $\vin{h}$, namely 
\stepcounter{sth}
\begin{equation*}\label{eq:atr-h}\tag{\thesth}
\tin{\!\!\!h}=\Hr^{-n}\Bigl(h_*f^!\co_Z \xleftarrow{\Iso} h_*h^!g^!\co_Z \xrightarrow{\Tr{h}} g^!\co_Z\Bigr). \end{equation*}
 Similarly define $\tin{\!\!h'}$, $\tin{\!\bar{h}}$, and
$\tin{\!\bar{h}'}$. Since $\bar{g}$ and $\bar{f}$ are proper, 
$\tin{\!\!\!\bar{h}}\colon \bar{h}_*\omgs{\bar{f}}\to \omgs{\bar{g}}$ has
the following alternative description: It is the adjoint to the element
of $\Hom_Y(\Rr^n{\bar{g}}_*\bar{h}_*\omgs{\bar{f}},\,\co_Z)$ given by the composite
\[\Rr^n{\bar{g}}_*\bar{h}_*\omgs{\bar{f}} \xleftarrow{\Iso}\Rr^n{\bar{f}}_*\omgs{\bar{f}}
\xrightarrow{\tin{\!\!\!\bar{f}}} \co_Z.\]

Let $\theta_u^{\bar{f}}\colon \bar{w}^*\omgs{\bar{f}}\to \omgs{\bar{f}'}$ and 
$\theta_u^{\bar{g}}\colon \bar{v}^*\omgs{\bar{g}}\to \omgs{\bar{g}'}$
be the base change isomorphisms defined in \cite[pp.\,738--739, Rmk.\,2.3.2, especially (2.5)]{cm}.
We claim that the following diagram commutes:
\[
\begin{aligned}
{\xymatrix{
\bar{v}^*\bar{h}_*\omgs{\bar{f}} \ar[d]_{\bar{v}^*\tin{\!\!\bar{h}}} \ar@{=}[r] & \bar{h}'_*\bar{w}^*\omgs{\bar{f}} \ar[r]^{\bar{h}'_*\theta_u^{\bar{f}}}
& \bar{h}'_*\omgs{\bar{f}'}\ar[d]^{\tin{\!\!\!\bar{h}'}}\\
\bar{v}^*\omgs{\bar{g}} \ar[rr]_{\theta_u^{\bar{g}}} && \omgs{\bar{g}'}.
}}
\end{aligned}\tag{\text{$*$}}
\]
Suppose $(*)$ commutes. Restricting $(*)$ to $Y$, and using the Verdier isomorphisms
for $f$, $\bar{f}$, $g$, and $\bar{g}$ and \cite[p.\,739, Theorem 2.3.3, especially (c)]{cm} 
(which states that via these isomorphisms
$\theta_u^f$ and $\theta_u^g$ are the identity maps) we get 
$v^*\vin{h}=\vin{h'}$
as we wish. The commutativity of $(*)$ is equivalent to 
\[
\tin{\!\!\!\bar{g}'}(\Rr^n\bar{g}'_*(\tin{\!\!\!\bar{h}'}\smcirc\bar{h}'_*\theta_u^{\bar{f}}))
= \tin{\!\!\!\bar{g}'}(\Rr^n\bar{g}'_*(\theta_u^{\bar{g}}\smcirc\bar{v}^*\tin{\!\!\bar{h}}))
\leqno{(\dag)}
\]
The proof of $(\dag)$ rests on the fact that the following diagram of functors commutes
\[
\begin{aligned}
{\xymatrix{
u^*\Rr^n\bar{f} \ar[r]^{\Iso} \ar[dd]_{\,\rotatebox{-90}{\makebox[-0.1cm]{\Iso}}}
& u^*\Rr^n\bar{g}_*\bar{h}_* \ar[dr] &\\
& &  \Rr^n\bar{g}'_*\bar{v}^*\bar{h}_* \ar@{=}[dl] \\
\Rr^n\bar{f}'_*\bar{w}^* \ar[r]^{\Iso} & \Rr^n\bar{g}'_*\bar{h}'_*\bar{w}^* \ar@{=}[ru] &
}}
\end{aligned}\tag{\ddag}
\]

In greater detail, consider the following diagram:
\[
{\xymatrix{
u^*\Rr^n\bar{f}_*\omgs{\bar{f}} \ar@{=}[r]  \ar[dddd]_{u^*\tin{\!\!\!f}}
& u^*\Rr^n\bar{f}_*\omgs{\bar{f}}  \ar[d]_{\,\rotatebox{-90}{\makebox[-0.1cm]{\Iso}}}\ar[rrr]^{\Iso}
&&& u^*\Rr^n\bar{g}_*\bar{h}_*\omgs{\bar{f}} \ar@{=}[d] \\
& \Rr^n\bar{f}'_*\bar{w}^*\omgs{\bar{f}} \ar[d]_{\theta_u^{\bar{f}}} \ar[r]^-{\Iso}
& \Rr^n\bar{g}'_*\bar{h}'_*\bar{w}^*\omgs{\bar{f}} \ar@{=}[r]  
 \ar[d]^{\theta_u^{\bar{f}}} \ar@{}[ddr]|\blacksquare
 & \Rr^n\bar{g}'_*\bar{v}^*\bar{h}_*\omgs{\bar{f}} \ar[dd]_{\tin{\!\!\!{\bar{h}}}} 
& u^*\Rr^n{\bar{g}}_*\bar{h}_*\omgs{\bar{f}}  \ar[dd]^{\tin{\!\!\!\bar{h}}}\ar[l]_{\,\,\,\Iso}\\
& \Rr^n\bar{f}'_*\omgs{\bar{f}'} \ar[dd]_{\tin{\!\!\!\bar{f}'}}  \ar[r]^-{\Iso}
  & \Rr^n\bar{g}'_*\bar{h}'_*\omgs{\bar{f}'} \ar[d]^{\tin{\!\!\!\bar{h}'}} 
& & \\
& & \Rr^n\bar{g}'_*\omgs{\bar{g}'} \ar[d]^{\tin{\!\!\!\bar{g}'}} & \Rr^n\bar{g}'_*\bar{v}^*\omgs{\bar{g}} \ar[l]^{\theta_u^{\bar{g}}} 
& u^*\Rr^n\bar{g}_*\omgs{\bar{g}} \ar[l]_{\Iso} \ar[d]^{\tin{\!\!\!\bar{g}}}\\
u^*\co_Z \ar@{=}[r] & \co_{Z'} \ar@{=}[r] & \co_{Z'} \ar@{=}[rr] && u^*\co_Z
}}
\]

The outer border commutes because of our alternate description of $\tin{\!\!\bar{h}}$.
The rectangle on the left commutes because of the definition of $\theta_u^{\bar{f}}$.
The rectangle on the lower right commutes because of the definition of $\theta_u^{\bar{g}}$
The remaining rectangle bordering the bottom edge commutes because of the alternate
description of $\tin{\!\bar{h}'}$. The rectangle on the top right is simply $(\ddag)$ and so
commutes.
All other rectangles, save $\blacksquare$, commute
for functorial reasons. Consider $\blacksquare$. We have two possible routes from its
northeast vertex to $\co_{Z'}$ lying directly below its southwest vertex, namely, south followed by
west followed by south, and west followed by south all the way. We have to show that the two
routes give the same map. This follows from the fact that all the subrectangles  
(except possibly $\blacksquare$) and the outer border commute.
This establishes $(\dag)$ and hence the theorem.
\qed

\medskip

We wish to understand \eqref{map:tr-h-fin} more explicitly. For that we need to work
more locally, with affine schemes, and often in a ``punctual way", i.e., by working
with completions of local rings at points. With this in mind, let us assume that
we are in the situation of diagram \eqref{diag:f=gh}, with a small change in hypothesis,
namely we assume $h$ is separated and quasi-finite, rather than finite. The maps $f$
and $g$ remain smooth of relative dimension $n$.
\[
{\xymatrix{
X \ar[d]_h \ar@/^1.0pc/[rrd]^f& & \\
Y \ar[rr]^g && Z
}}
\]

We are interested in duality for $h$ in terms of $\omega_g$ and $\omega_f$
 ``at a point $x\in X$". To that end we make the following further assumptions.
 \begin{itemize}
 \item $Z=\Spec{\,A}$
 \item $Y=\Spec{\,R}$ and $X=\Spec{\,S}$.
 \end{itemize}
 Let $y\in Y$ assume $h^{-1}(y)$ consists of exactly one point $x$.
 
 Let $R'=\wid{\co_y}$ be the completion of the local ring $\co_{Y,y}$, 
 $S'=\wid{\co_x}$ the completion
 of $\co_{X,x}$, and set $Y'=\Spec{\,R'}$, $X'=\Spec{\,S'}$. Since $h^{-1}(y)=\{x\}$,
 we have a cartesian square
 \[
 {\xymatrix{
 X' \ar[d]_{h'} \ar[r]^v \ar@{}[dr]|\square & X \ar[d]^h \\
 Y' \ar[r]_u & Y
 }}
 \]
with $h'$ \emph{finite}, even though $h$ need not be finite.

To lighten notation, write $\omega_R=\omega_{R/A}$, and $\omega_S=\omega_{S/A}$.
Set $\omega_{R'}=\omega_R\otimes_RR'$, and 
$\omega_{S'}=\omega_S\otimes_SS'= \omega_S\otimes_RR'$.

Since $h$ is flat and Gorenstein of relative dimension $0$, for any quasi-coherent $\co_Y$-module
$\eF$ we have $\Hr^k(h^!\eF)=0$ for $k\neq 0$, and so we identify $h^!\eF$ with $\Hr^0(h^!\eF)$.
Similarly, we identify $h'^!\eG$ with $\Hr^0(h'^!\eG)$ for every quasi-coherent $\co_{Y'}$-module
$\eG$. For an $R$-module $M$, $h^!M$ is defined to be $\Gamma(X,\,h^!\wit{M})$. Similarly,
for an $R'$-module $N$, $(h')^!N$ will denote $\Gamma(X',\,(h')^!\wit{N})$.

Let 
\stepcounter{thm}
\begin{equation*}\label{iso:h^!OR}\tag{\thethm}
\varsigma\colon h^!\omega_R \iso \omega_S
\end{equation*}
denote the isomorphism obtained from $h^!g^!\co_Y\iso f^!\co_Y$ and the Verdier isomorphisms
${\bf v}_g$ and ${\bf v}_f$.
 By (flat) base change, we have
 \stepcounter{thm}
\begin{equation*}\label{iso:h^!OR'}\tag{\thethm}
\varsigma'\colon (h')^!\omega_{R'} \iso \omega_{S'}.
\end{equation*}
In particular we have a trace map (for $h'$ is finite)
\stepcounter{thm}
\begin{equation*}\label{map:trS'}\tag{\thethm}
 \vin{S'}\colon \omega_{S'}\to \omega_{R'}
 \end{equation*}
corresponding to
\[
h'_*\wit{\omega}_{S'} \xleftarrow[{(*)_{R'}}]{\Iso} h'_*(h')^!\wit{\omega}_{R'} \xrightarrow{\Tr{h'}}
\wit{\omega}_{S}
\]

Our interest is in making $\vin{S'}$ explicit.We point out that to define it, it was not necessary
to assume that $x$ is the only point of $X$ lying over $y$. However, by shrinking $X$ around
$x$, we can be in the situation we are in.

Now suppose $h\colon X\to Y$ factors as in the following commutative
diagram
\[
 {\xymatrix{
 X\, \ar@{^(->}[r]^i \ar[dr]_h & P \ar[d]^p \\
 & Y
 }}
 \]
with $P=\Spec{\,E}$, $p\colon P\to Y$ \emph{smooth of relative dimension $d$}, and
$i$ a \textit{closed immersion}.
We have a commutative diagram with each square cartesian
\stepcounter{thm}
\[
\begin{aligned}\label{diag:ijX}
{\xymatrix{
X'\ar[r]^v \ar[d]_j \ar@{}[dr]|\square & X \ar[d]^i \\
P' \ar[r]^w \ar[d]_{p'} \ar@{}[dr]|\square & P \ar[d]^p \\
Y' \ar[r]_u & Y
}}
\end{aligned}\tag{\thethm}
\]
with $h=p\smcirc i$ and $h'=p'\smcirc j$.

Let $E'=E\otimes_RR'$, $P'=\Spec{\,E'}$, $\pi=g\smcirc p$, and consistent with out notations
above, let $\omega_E=\omega_{E/A}$, and $\omega_{E'}=\omega_E\otimes_RR'$.

We remark that $\omega_{R'}$ and $\omega_{S'}$ are the $e$-th graded pieces of the differenital
graded algebras $\wedge^\bullet_{R'}(\Omega^1_{R/A}\otimes_RR')$ and 
$\wedge^\bullet_{S'}(\Omega^1_{S/A}\otimes_SS')=
\wedge^\bullet_{S'}(\Omega^1_{S/A}\otimes_RR')$ respectively. Similarly,
$\omega_{E'}$ is the $(n+e)$-th graded piece of $\wedge^\bullet_{E'}(\Omega^1_{E/A}\otimes_EE')
=\wedge^\bullet_{E'}(\Omega^1_{E/A}\otimes_RR')$

Let
\[\phi\colon \omega_R\otimes_R \omega_{E/R} \iso \omega_E\]
be the isomorphism $\phi=\Gamma(Y,\,\bar{\varphi}_{g,\,p})$, where $\bar{\varphi}_{g,\,p}$
is the map defined in \eqref{iso:phi-gf-bar}. In other words $\phi(\nu\otimes\mu) = \mu\wedge p^*\nu$.
Let 
\[\phi'  \colon \omega_{R'}\otimes_{R'} \omega_{E'/R'} \iso \omega_{E'}\]
be the base change of $\phi$. In greater detail, we have 
$\omega_{E'/R'}=\omega_{E/R}\otimes_RR'$, and therefore
$\omega_{R'}\otimes_{R'} \omega_{E'/R'}= (\omega_R\otimes_R \omega_{E/R})\otimes_RR'$.
Set $\phi'=\phi\otimes{\bf 1}$.

Next, let $I=\ker{E\twoheadrightarrow S}$, $J=\ker{E'\twoheadrightarrow S'}$.
Write $N=(\wedge^d_S I/I^2)^*$ and $N'=(\wedge^d_{S'} J/J^2)^*$. Let
\[b\colon \omega_E\otimes_E N \iso \omega_S\]
be the map given by \eqref{eq:b-X/P}. By base change, as in the definition of $\phi'$, we have
a map $b'\set b\otimes{\bf 1}$:
\[b'\colon \omega_{E'}\otimes_{E'} N'\iso \omega_{S'}.\]
Let $\varrho\colon \omega_R\otimes_R \omega_{E/R}\otimes_EN \to \omega_S$
and $\varrho'\colon\omega_{R'}\otimes_{R'} \omega_{E'/R'}\otimes_{E'}N' \to \omega_{S'}$
be the maps
\stepcounter{thm}
\begin{equation*}\label{def:vrho}\tag{\thethm}
\varrho=b\smcirc (\phi\otimes {\bf 1}_N) \quad \text{and} \quad 
\varrho'=b'\smcirc(\phi'\otimes{\bf 1}_{N'})= \varrho\otimes {\bf 1}_{R'}.
\end{equation*}
Finally, let
\stepcounter{thm}
\begin{equation*}\label{def:pi-h}\tag{\thethm}
\psi\colon \omega_{E/R}\otimes_E N \iso h^!R \quad \text{and} \quad 
\psi'\colon \omega_{E'/R'}\otimes_{E'}N' \iso {h'}^!(R')
\end{equation*}
be the maps defined as in \eqref{iso:p-i!}.

\begin{prop}\label{prop:fin-trace} In the above situation, assume $I=\ker{(E\twoheadrightarrow S)}$ is 
generated by ${\bf u}=(u_1, \dots, u_d)$ and set $f_k=u_k\otimes 1\in E'$,
so that ${\bf f}=(f_1, \dots, f_d)$ generates $J=\ker{(E'\twoheadrightarrow S')}$. Set
$\eN=\eN^d_i(=\wit{N})$ and $\eN'=\eN^d_j (=\wit{N'})$. Let $\omega\in \omega_{S'}$.
\begin{enumerate}
\item The following diagram commutes:
\[
{\xymatrix{
(\omega_{R'}\otimes_{R'}\omega_{E'/R'})\otimes_{S'}N' \ar@{=}[d] \ar[rr]^-{\varrho'} & & \omega_{S'} \\
\omega_{R'}\otimes_{R'}(\omega_{E'/R'}\otimes_{S'}N') 
\ar[d]_{{\bf 1}\otimes \psi'}^{\rotatebox{-90}{\makebox[0.1cm]{\Iso}}} && \\
\omega_{R'}\otimes_{R'} {h'}^!(R') \ar[rr]^{\Iso}_{\chi^{{}_{h'}}} && {h'}^!(\omega_{R'}) 
\ar[uu]_{\varsigma'}^{\rotatebox{-90}{\makebox[0.1cm]{\Iso}}} 
}}
\]
where $\chi^{{}_{h'}}$ is the transitivity map defined in \eqref{iso:chi-h}.
\item If $\omega=s\cdot{(h')}^*(\nu)$, where $\nu\in\omega_{R'}$ and $s\in S'$, then 
\[\vin{S'}(\omega)=\trc{S'/R'}(s)\cdot \nu.\]
\item Let $\eta\in \Omega^n_{E/A}\otimes_RR'$ be any element such that $j^*\eta=\omega$.
Then
\[\vin{S'}(\omega)= \nu\cdot \res{{X', p'}}\begin{bmatrix} x\cdot\mu\\
f_1, \dots, f_d\end{bmatrix}\]
where $x\in E'$, $\mu\in \omega_{E'/R'}$, and $\nu\in\omega_{R'}$ are
related via the formula 
\[\wdd{f_1}{f_d}\wedge\eta = x\cdot\mu\wedge{p'}^{{}_*}\nu.\]
\item Let $\mu$, $\nu$, $\eta$ and $x$ be as in (ii).
Suppose $E=R[T_1,\dots, T_d]_{q({\bf T})}$, where $q({\bf T})\in R[{\bf T}]$. Let
$g_i\in E'\otimes_{R'}S'$ be the elements $g_i=T_i\otimes 1 - 1\otimes \gamma_i$, where the
$\gamma_i\in S'$ is the images of $T_i$, $i=1, \dots, d$, and let
\[\lambda\colon S'\to R'\] 
be the map $\lambda=\lambda({\bf f},\,{\bf g})$ of \Tref{thm:tate2}. 
Then
\[\vin{S'}(\omega) = \lambda(x\vert_{X'})\nu.\]

\end{enumerate}
\end{prop}

\proof We point out that ${\bf u}$ and ${\bf f}$ are necessarily quasi-regular.
We first prove (i). Consider the following diagram.
\[
{\xymatrix{
&& v^*(h^*\omega_g\otimes i^*\omega_p\otimes\eN) 
\ar[rr]^-{v^*(\varrho)}_-{\,\rotatebox{180}{\makebox[-0.1cm]{\Iso}}}
\ar@{=}[d] \ar@/_8.5pc/[ddd]^{v^*({\bf 1}\otimes\psi)}_\simeq 
\ar@{}[drr]|{\phantom{XXXXX}\blacklozenge}
 & & v^*\omega_f \ar@{=}[d] && \\
&& {h'}^{{}_*}\wit{\omega}_{R'}\otimes j^*\omega_{p'}\otimes\eN' 
\ar[rr]_-{\varrho'}^-{\Iso} 
\ar[d]_{{\bf 1}\otimes\psi'}^{\>\rotatebox{90}{\makebox[-0.1cm]{\Iso}}} 
\ar@{}[dll]|{\blacksquare\phantom{X}}
 && \wit{\omega}_{S'}  \ar@{}[drr]|\square  
 && \\
&& {h'}^{{}_*}\wit{\omega}_{R'}\otimes {h'}^!\co_{Y'} \ar@{}[drr]|{\phantom{XXX}\lozenge}
 \ar[rr]_-{\chi^{{{}_{h'}}}}^-\Iso &&
 {h'}^!\wit{\omega}_{R'}  \ar[u]_{\varsigma'}^{\rotatebox{-90}{\makebox[0.1cm]{\Iso}}}
  && \\
&& v^*(h^*\omega_g\otimes h^!\co_Y)
\ar[u]_{\rotatebox{-90}{\makebox[0.1cm]{\Iso}}}^{{\bf 1}\otimes\theta_u^h} 
 \ar[rr]_-{v^*(\chi^{{}_h})}^-\Iso && v^*h^!\omega_g 
\ar[u]^{\,\rotatebox{90}{\makebox[-0.1cm]{\Iso}}}_{\theta_u^h} 
\ar@/_7.5pc/[uuu]^{v^*(\varsigma)}_\simeq 
&&
}}
\]

The rectangle $\blacklozenge$ on the top commutes by definition of $\varrho'$. 
The sub-diagram on the right,
the one labelled $\square$, squeezed between the curved arrow and the vertical column, 
commutes by definition of $\varsigma'$. The rectangle labelled $\lozenge$ at the bottom commutes
by \cite[Prop.\,7.2.10]{fub-abs}.

We now show that the sub-diagram on the left, labelled $\blacksquare$, squeezed between the curved
arrow and the vertical column on the left, commutes. 
First, the composite of isomorphisms, with the middle arrow the base change isomorphism
\[v^*\eN[-d] \xrightarrow{v^*\eta'_i} v^*i^!\co_P \iso j^!\co_{P'} \xrightarrow{{\eta'_j}^{-1}} \eN'[-d]\]
is the identity map on 
$\eN'[-d]$ \cite[Remark 6.2.5]{fub-abs}. Next,
the composite of isomorphisms
\[w^*\omega_p[d] \xrightarrow{w^*{\bf v}_p}w^*p^!\co_Y \iso 
p'^!\co_{Y'} \xrightarrow{{\bf v}_{p'}^{-1}} \omega_{p'}[d]\]
is the identity map on $\omega_{p'}[d]$ \cite[p.\,740, Prop.\,2.3.5\,(b)]{cm}. 
Finally, the transitivity property of base change
\cite[Prop.\,A.1.1\,(ii)]{fub-abs} (see also
\cite[p.\,183, Prop.\,4.6.8]{notes}) tells us that the base change of the composite $p\smcirc i$
with respect to $u$ is compatible with the base change for $p$ and $i$ with respect to $u$ and
$w\colon P'\to P$ respectively. Putting these together, we see that $\blacksquare$ also commutes.

The outer border commutes by \Pref{prop:h!}, after using \Tref{thm:res-thm}
to realise $b$ as a concrete representation of the map $a_{X/P}$

It follows that the rectangle in the middle also commutes. This proves (i).

Next note that the following diagram commutes, by definition of the various isomorphisms
involved.
\[
{\xymatrix{
\omega_{R'}\otimes_{R'}h'^!(R') \ar[d]_{{\bf 1}\otimes \Tr{S'/R'}} \ar[rr]^-{\Iso}_-{\chi^{{}_{h'}}}
&&  {h'}^!(R') \ar[d]_{\Tr{S'/R'}} \ar[rr]^{\Iso}_{\varsigma'} && \omega_{S'} \ar[d]^{\vin{S'}} \\
\omega_{R'} \ar@{=}[rr] &&  \omega_{R'} \ar@{=}[rr] && \omega_{R'}
}}
\]
Let $\ttr{h'}=\ttr{h', p', j}\colon h'_*(j^*\omgs{p'}\otimes\eN') \iso \co_{Y'}$
be the map in \cite[(5.3.2)]{fub-abs}.
Define 
\[\vttr{h'} \colon  h'_*(j^*\omega_{p'}\otimes\eN') \iso \co_{Y'}\]
in the obvious way, namely by substituting $\omgs{p'}$ in the definition of $\ttr{h'}$
by $\omega_{p'}$ via the Verdier isomorphism ${\bf v}_{p'}$. Write $\vttr{S'/R'}$ for
the global sections of $\vttr{h'}$. From part (i) and the above commutative diagram,
we see that
\[\vin{S'}\smcirc \varrho'={\bf 1}\otimes \vttr{S'/R'}. \leqno{(*)}\]
Now suppose $\mu\in \omega_{E'/R'}$ and $\nu\in\omega_{R'}$, 
By \cite[Prop.\,5.4.4]{fub-abs} and $(*)$ we
get
\[\vin{S'}(\rho'(\nu\otimes \mu\otimes{\bf 1/f}))
=\res{{X',p'}}\begin{bmatrix}\mu\\ f_1, \dots, f_d \end{bmatrix}\cdot\nu \leqno{(\dag)} \]
Now if $\omega= s\cdot {h'}^{{}_*}(\nu)$, then by definition of $\varrho'$, if $x\in E'$
is a lift of $s$, we have
\[\varrho'(x\cdot (\nu\otimes\wdd{f_1}{f_d}\otimes {\bf 1/f})= \omega,\]
whence by $(\dag)$
\[\vin{S'}(\omega)= \res{{X',p'}}\begin{bmatrix} x\cdot\wdd{f_1}{f_d}\\ f_1, \dots, f_d \end{bmatrix}
\cdot\nu.\]
The right side is equal to $\trc{S'/R'}(s)\cdot\nu$ by \Tref{thm:R6}. This proves (ii).
Part (iii) is a re-statement of $(\dag)$. Indeed
\[
\begin{aligned}
\vin{S'}(\omega)=\vin{S'}(j^*\eta)&=\vin{S'}(b'(\wdd{f_1}{f_d}\wedge\eta\otimes{\bf 1/f}))\\
&=\vin{S'}(\varrho'(x\cdot(\nu\otimes\mu\otimes{\bf 1/f})))\\
&= \nu\cdot \res{{X', p'}}\begin{bmatrix} x\cdot\mu\\
f_1, \dots, f_d\end{bmatrix}
\end{aligned}
\]
Part\,(iv) follows from (iii) and \Pref{prop:lambda-res}.
\qed 

\begin{rem}\label{rem:q-finite} We have already observed that it was not necessary to
assume $h^{-1}(y)$ consisted of exactly one point, namely $x$, in order to define
$\vin{S'}$. In fact more can be be said.
Suppose $u\colon X\to \overline{X}$ is an open immersion of 
(ordinary) $Y$-schemes
such that the structure map $\bar{h}\colon \overline{X}\to Y$ is quasi-finite, and
${\bar h}^{-1}(y)=\{x_1, \dots, x_m\}$, with $x_1=x$. 
In this case, the fibre dimension of $\bar{f}=g\smcirc\bar{h}$ is $n$. As before,
set $\omgs{{\bar f}}= \Hr^{-n}({\bar f}^!\co_Z)$. 
Now, $\overline{X}'=\overline{X}\times_YY'$
is finite over $Y$, since $Y'$ is the spectrum of a complete local ring. Let 
$\bar{h}'\colon \overline{X}'\to Y'$ be the base change of $\bar{h}$. Now
$\overline{X}'=\Spec{\,\prod_{i=1}^m S_i'}$, where $S_i'$ is the completion
of the local ring $S_i=\co_{\overline{X}, x_i}$. Let $\overline{S}'=\prod_iS_i'$, and
let $X_i'=\Spec{\,S_i'}$, so that
$X_i'$ is open and closed in $\overline{X}'$, and $\overline{X}'=\coprod_i X_i'$.
Let $h_i'\colon X_i'\to Y'$ be the restriction of $\bar{h}'$ to $X_i'$.
Note $X_1'=X'$, $S_1'=S'$, and $h_1'=h'$. If $\omgs{S'_i}=\omgs{\bar{f},\,x_i}\otimes_{S_i}S_i'$,
and $\omgs{\overline{S}'}=\oplus_i\omgs{S'_i}$ (the direct sum thought of as an $\overline{S}'$-
module, then, as in the argument used in \eqref{iso:h^!OR'}, we have, analogous to
$\varsigma'$, isomorphisms
$\Hr^0(h_i'^!\omega_{R'})\iso \omgs{S_i'}$ and 
$\Hr^0(({\bar h}')^!\omega_{R'}) \iso \omgs{\overline{S}'}$,\footnote{Regarding 
$({\bar h}')^!\omega_{R'}$ as a complex of $\overline{S}'$-modules associated to
$({\bar h}')^!\wit{\omega}_{R'}$ etc.} whence 
abstract trace maps 
\[\tin{\!\!\!S_i'}\colon \omgs{S_i'} \lra \omega_{R'} \qquad (i=1,\dots m)\]
and 
\[\tin{\!\!\!\overline{S}'}\colon \omgs{\overline{S}'} \lra \omega_{R'}.\]
Clearly
$\tin{\overline{S}'}= \sum_i\tin{S_i'}$.
The Verdier isomorphism ${\bf v}_f\colon \omega_f\iso \omgs{f}$ base changes
to ${\bf v}_{S'}\colon \omega_{S'}\iso \omgs{S'}$ and it is clear that
$\vin{S'} = \tin{\!\!S'}\smcirc {\bf v}_{S'}$.
Finally, \textit{if $\bar{h}\colon \overline{X}\to Y$ is finite}, say 
$\overline{X}=\Spec{\,\overline{S}}$, then we have a map
$\tin{\!\!{\bar h}}\colon {\bar h}_*\omgs{{\bar f}}\to \omgs{g}$ defined in \eqref{eq:atr-h}. 
Consistent with the above notations, set 
$\omgs{\,\overline{S}}=\Gamma(\overline{X},\,\omgs{\bar{f}})$.
Let 
$\tin{\!\!\overline{S}}\colon \omgs{\overline{S}}\to \omega_R$
be the map $\Gamma(Y,\, {\bf v}_g^{-1}\smcirc \tin{\!\!{\bar h}})$.
Then clearly
\[\tin{\!\!\!\overline{S}}\otimes_RR' = \tin{\!\!\overline{S}'} = \sum_k\tin{\!\!\!S_k'},\]
where $\tin{\!\!\!\overline{S}}$ is the global sections of $\tin{\!\!\!h}$ defined in \eqref{eq:atr-h}.
In particular, if $\bar{f}\colon \overline{X}\to Z$ is smooth, then with $\vin{\overline{S}}=\Gamma(\overline{X},\,\vin{f})$ we have
\stepcounter{sth}
\begin{equation*} \label{eq:trS=sum-trS'}\tag{\thesth}
\vin{\overline{S}}\otimes_RR'=\sum_k\vin{S_k'}.
\end{equation*}
\end{rem}

\setcounter{subsubsection}{\value{thm}} 
\subsubsection{The Kunz-Lipman trace} \label{sss:joe-kunz} \stepcounter{thm} 
Suppose, as we have for most of this section,
 $X=\Spec{\,S}$, $Y=\Spec{\,R}$, and $Z=\Spec{\,A}$, and as before suppose $f\colon X\to Z$ and
 $g\colon Y\to Z$ are smooth of relative dimension $n$, $f=g\smcirc h$, and now assume
 $h$ is \textit{finite}, and not merely separated and quasi-finite. In this case (and in
 more general situations) we have a trace map
 \[\sigma_{S/R} \colon \omega_S \lra \omega_R\]
 or, in sheaf-theoretic terms,
 \[\sigma_h\colon h_*\omega_f \lra \omega_g\]
 due to Lipman and Kunz, defined in Kunz's book \cite[p.\,254, 16.4]{kd}. The idea is attributed
 by Kunz to Lipman (see footnote in \textit{loc.cit.})
 
 The \emph{Kunz-Lipman trace} $\sigma_{S/R}$ can be understood punctually. In greater
 detail, the Tate trace $\lambda({\bf f}, {\bf g})$ of \Tref{thm:tate2} is denoted $\tau^x_f$ in
 \cite[p.\,370, (F.20)]{kd} (and studied in some detail in F.18--F.28 of \textit{ibid}).
 Now suppose $y$ is a point in $Y$. Fix $x\in h^{-1}(y)$ and pick an affine open subscheme
 $U=\Spec{\,S_U}$ of $X$ such that $h^{-1}(y)\cap U =\{x\}$, and a presentation
 \[R[T_1, \dots, T_d]_{q({\bf T})}/(u_1, \dots, u_d) = S_U.\]
 Such a $U$ and presentation always exists. Let $R'$ be as before, the completion of the
 local ring $\co_{Y,y}$, and let $S'$ be the completion of the local ring $\co_{X,x}$. 
 Let $E=R[{\bf T}]_{q({\bf T})}$, and $E'=E\otimes_RR'$. Let $f_1, \dots, f_d$ be the images
 of $u_1, \dots, u_d$ in $E'$. We continue to denote the image of the variables $T_k$ in
 $E'$ as $T_k$. Let $\omega_R$, $\omega_{R'}$, $\omega_{S_U}$, $\omega_{S'}$, $\vin{S'}$ etc.,
  be as before. Let $\gamma_k\in S'$ be the image of $f_k$, and set 
 $g_k= X_k\otimes 1 - 1\otimes \gamma_k \in E'\otimes_{R'}S'$.
 Finally let
 \[\lambda\colon S'\to R'\]
 be the Tate trace $\lambda({\bf f}, {\bf g})$ of \Tref{thm:tate2}. Since $\omega_{S'}$ is
 a direct summand of $\omega_S\otimes_RR'$, the map $\sigma_{S/R}$ restricts to
a map 
\[\sigma_{S'}\colon \omega_{S'}\to R'.\] 

For $\omega\in \omega_{S'}$ and $\eta\in \omega_{E'}$ a pre-image of $\omega$ under
the natural surjective map $\omega_{E'}\to \omega_{S'}$, suppose $x\in E'$, $\nu\in \omega_{R'}$
are such that
\[\wdd{f_1}{f_d}\wedge\eta = x\cdot \wdd{T_1}{T_d}\wedge\nu.\]

 Using properties Tr 3) and Tr 4) of \cite[pp.\,245-246, \S\,16]{kd}, proved in
 [\textit{ibid}, p.\,254, Thm.\,16.1], the definition of the Kunz-Kipman trace in [\textit{ibid}, p.\,254, 16.4]
 gives
 \[\sigma_{S'}=\lambda(\bar{x})\cdot \nu,\]
 where ${\bar x}\in S'$ is the image of $x\in E'$. This means, by the formula in \Pref{prop:fin-trace}\,(iv),
 \[\sigma_{S'}=\vin{S'}.\] 
 Once again by the above mentioned properties Tr 3) and Tr 4) of 
 $\sigma_{S/R}$,  and by \eqref{eq:trS=sum-trS'},
 this gives $\widehat{\sigma_{h,y}}= \widehat{\vin{h,y}}$,
 where $\widehat{(\boldsymbol{-})}$ denote completion of an $\co_{Y,y}$-module with respect
 to the maximal ideal. Since $y$ is arbitrary in $Y$, we have
 \stepcounter{sth}
 \begin{equation*}\label{eq:sigma=tr}\tag{\thesth}
 \sigma_h =\vin{h}.
 \end{equation*}
 Clearly, we don't need $X$, $Y$ and $Z$ to be affine for the argument to go through.

\begin{thm}\label{thm:kunz-tr1} 
Let $f\colon X\to Z$ and $g\colon Y\to Z$ be smooth separated maps of ordinary schemes
of relative dimension $n$, and suppose $f=g\smcirc h$, where $h\colon X\to Y$ is
a finite map. 
\begin{enumerate}
\item Suppose $Z=\Spec{\,A}$, $Y=\Spec{\,R}$ and $X=\Spec{\,S}$.
Let $\omega= s(h^*(\nu))$ where $s\in S$ and $\nu\in\Omega^n_{R/A}$. Then
\[\vin{h}(\omega)=\trc{S/R}(s)\,\nu.\]
\item If $\sigma\colon h_*\omega_f\to \omega_g$ is the Kunz-Lipman trace then
\[\sigma_h =\vin{h}.\]
\end{enumerate}
\end{thm}

\proof Part\,(ii) of \Pref{prop:fin-trace} gives (i). Part\,(ii) is simply \eqref{eq:sigma=tr}.
\qed

\subsection{Regular Differentials again}\label{ss:reg-diff-again} 
The fact that $\vin{h}$ agrees with the Lipman-Kunz trace
$\sigma_h$ allows us to prove
\Tref{thm:reg-ver} in a different way. Suppose $f$, $g$, $h$ are as above, with the caveat that
we no longer assume that $f$ is smooth, but assume $f$ is of finite type, the smooth locus of
$f$, $X^{\rm{sm}}$ contains all the  associated points of $X$. The map $g$ remains smooth,
and $h$ finite. Assume further that:
\begin{enumerate}
\item $X$, $Y$ and $Z$ are excellent have no embedded points;
\item $X=\Spec{\,S}$,
$Y=\Spec{\,R}$, and $Z=\Spec{\,A}$;
\item $R\to S$ is injective.
\end{enumerate}

We use the notations of \Sref{s:reg-diffs}. Thus $f^K=\Hr^{-n}(f^!)$ is as in \Ssref{ss:kleiman},
and $\oreg{X/Z}$ is the sheaf of regular differential $n$-forms discussed
in \Ssref{ss:reg-ver}. Since we are in the affine situation, we work with modules and algebras
over $A$, $R$, and $S$, and choose appropriate notations. To that end, let $k(R)$ and $k(S)$
be the total ring of fractions of $R$ and $S$ respectively. Set 
$\omega_R=\Gamma(X,\,\omega_g)$, $\omega_{k(S)}=\Gamma(X,\,\Omega^n_{k(X)/k(Z)})$,
$\Omega_{k(R)}=\Gamma(X,\,\Omega^n_{k(Y)/k(Z)})$. 

Standard arguments show that there is an scheme theoretically dense open subscheme $U$
of $Y$, such that $h^{-1}(U)$ is in $X^{\rm{sm}}$ and is scheme theoretically dense in $X$
(e.g., $U=Y\smallsetminus h(X\smallsetminus X^{\rm{sm}})$). We have the trace map
$\vin{h_U}\colon (h\vert_U)_*\omega_{f\vert_{h^{-1}U}} \to \omega_{g\vert_U}$, where
$h_U\colon h^{-1}(U)\to U$ is the restriction of $h$. By taking stalks at generic points (we
have no embedded points!) we get a map
\[\vin{k(S)}\colon \omega_{k(S)} \lra \omega_{k(R)}.\]
We point out that $\omega_R\subset \omega_{k(R)}$. The content of the next result is that
$\overline{\omega}_S$ is a ``complementary module" in the sense of Kunz and Waldi
\cite[\S\,4]{kw}. It is equivalent to \Tref{thm:reg-ver}, via \Tref{thm:kunz-tr1}, but we give
a direct proof along the lines of the proof given of a related statement in \cite{ast117}.

\begin{thm}\label{thm:kunz-verd2} Let $\overline{\omega}_S \subset \omega_{k(S)}$ be
the image of the injective map $\omgs{S} \to \omega_{k(S)}$ defined in
\eqref{map:can}. Then
\[\overline{\omega}_S = \{\nu\in \omega_{k(S)}\mid \vin{k(S)}(s\nu) \in \omega_R, \forall s\in S\}.\]
\end{thm}
\proof
The proof is \textit{mutatis mutandis} the proof given in \cite[p.\,34, Lemma\,(2.2)]{ast117}.
We give it here, with the necessary changes, for completeness.
We have a natural isomorphism $\Hom_R(S,\,\omega_R) \iso \omgs{S}$
obtained by applying $\Hr^{-n}$ to $h^\flat \omega_g[n] \iso f^!\co_Z$, whence an
isomorphism
\[ \overline{\omega}_S \iso \Hom_R(S,\,\omega_R).\]
One checks (by using the open set $U=Y\smallsetminus h(X\smallsetminus X^{\rm{sm}})$ as
an intermediary if necessary) that the following diagram commutes
\[
{\xymatrix{
\Hom_R(S,\,\omega_R)\, \ar@{^(->}[r] 
& 
\Hom_{k(R)}(k(S),\,\omega_{k(R)})  \\
\overline{\omega}_S\, \ar@{^(-}[r] \ar[u]^{\,\rotatebox{90}{\makebox[0.1cm]{\Iso}}}
& \omega_{k(S)} \ar[u]_{\>\rotatebox{90}{\makebox[-0.1cm]{\Iso}}}
}}
\] 
where the isomorphism on the right is $\nu\mapsto (x\mapsto \vin{k(S)}(x\nu))$,
 for $\nu\in\omega_{k(S)}$ and $x\in k(S)$. The result follows since the image of
 $\Hom_R(S,\,\omega_R)$ in $\Hom_{k(R)}(k(S),\,\omega_{k(R)})$ consists of $k(R)$-linear maps
 $\psi\colon k(S)\to \omega_{k(R)}$ such that $\psi(s)\in \omega_R$ for every $s\in S$.  In other
 words, such $\psi$ are characterised by the property that
 ${\bf{e}}(s\psi)\in \omega_R$ for every $s\in S$, where 
 \[{\bf{e}}\colon \Hom_{k(R)}(k(S),\,\omega_{k(R)}) \to \omega_{k(R)}\]
 is ``evaluation at $1$". Since ${\bf{e}}$ corresponds to $\vin{k(S)}$ under the upward arrow
 on the right in the above diagram, we are done.
 \qed
 
 \smallskip
 
 The next statement is a re-statement of \Tref{thm:reg-ver}, but the point is that it is also
 a consequence of \Tref{thm:kunz-verd2}.
 \begin{cor}\label{cor:kunz-verd} Let $\oreg{S}$ be the $S$-module whose associated 
 quasi-coherent sheaf is $\oreg{X/Y}$. Then $\oreg{S}=\overline{\omega}_S$.
 \end{cor}
 \proof
 Let $U=Y\smallsetminus h(X\smallsetminus X^{\rm{sm}})$ and $h_U\colon h^{-1}(U)\to U$
 the restriction of $h$. From \Tref{thm:kunz-tr1}\,(ii), 
 $\vin{h_U}=\sigma_{h_U}$. The result follows
 from the characterisation of $\overline{\omega}_S$ as a complementary module
 in \Tref{thm:kunz-verd2} and the definition of regular differentials in \cite[p.58]{hk1}.
 \qed
 
 \setcounter{subsubsection}{\value{thm}} 
\subsubsection{} \label{sss:reg-ver} \stepcounter{thm} We would like draw out the differences
between the approach in \Sref{s:reg-diffs} and that of this subsection. In the former, we 
treat the theory of regular differential forms as a settled theory, and freely use the
results in \cite{hk1}, \cite{hk2}, and \cite{ajm} to arrive at a proof of \Tref{thm:reg-ver} 
using our characterisation of the Verdier isomorphism in terms of standard residues along
sections. In the ``settled theory" mentioned above, $\oreg{X/Z}$ is defined via local 
quasi-normalisations, i.e.~ via quasi-finite maps from open subschemes of $X$
to  $\mathbb{A}^n_Z$, their compactifications by Zariski's Main Theorem and
complementary modules \cite{kw}. The theory of residues and traces used there make
no reference to Verdier's isomorphism, and are developed \textit{ab initio} for the 
purpose at hand. In \Ssref{s:reg-diffs}, we mapped our theory on to all of that.

In contrast, from the results in this subsection, if $\oreg{X/Z}$ is {\emph{defined}} as
the image of the injective map $f^K\co_Y\to \Omega^n_{k(X)/k(Z)}$ as in \eqref{map:can},
then we show that every time one has a finite dominant $Z$-map $h\colon X\to Y$ of
schemes, such that $Y\to Z$ is smooth, then $\oreg{X/Z}$ is necessarily the complementary
module on the right side of \Tref{thm:kunz-verd2} (which can clearly be defined even when
$Y$ is not affine). Using \cite[${\rm{IV}}_3$, (13.3.2)]{ega} and Zariski's main theorem, as
in the first two paragraphs proof of \Pref{prop:kunz-sigma}, we see that locally we can
always arrange matters so that $X$ is covered by affine open subschemes, each of
which is finite over an affine smooth $Z$-scheme (in fact an affine open subscheme of
$\mathbb{P}^n_Z$), and hence $\oreg{X/Z}$ has a local description via complementary
modules. This gives a different proof, than that given in \cite{kw}, that these complementary
modules glue, and do not depend on the choice of the various finite maps of the sort just
discussed. Finally the theory of residues and traces presented here in this paper means that
all the important results in \cite{hk1}, \cite{hk2}, and \cite{ajm} can be recovered. 

Our approach (in this subsection) is closer
in spirit to the approach to these matters in \cite{ast117}, though even here it is necessarily
different, since we use, consistently, Verdier's isomorphism, and we work over an arbitrary
(noetherian) base rather than over a perfect field. 
It should be said that  in \cite{kw} and \cite{kd}, 
the theory is for general differential algebras, and that in \cite{kw}, generic
complete intersection algebras $A\to S$ are considered.

\section{\bf The Residue Symbol}\label{s:res-sym}

\subsection{Definition} Let $f\colon X\to Y$ be a separated smooth map of
relative dimension $r$, $t_1, \dots, t_r \in \Gamma(X,\,\co_X)$ such that if $\I$ is
the quasi-coherent ideal sheaf generated by ${\bf t}=(t_1,\dots, t_r)$ , then
$Z\set {\boldsymbol{\Spec}}\,(\co_X/\I)$ is finite over $Y$. Let $i\colon Z\to X$ be the closed
immersion and $h\colon Z\to Y$ the finite map. In this case it is well-known that
$h$ is flat and ${\bf t}$ is a regular  $\co_{X,z}$-sequence for every $z\in Z$
 (\cite[${\rm{IV}}_3$, Th\'eor\`eme (11.3.8)]{ega} or \cite[p.\,177, Corollary to Thm.\,22.5]{matsumura}).
In particular $h_*\co_Z$ locally free over $\co_Y$.  

In this situation, according to \cite[(5.3.2)]{fub-abs},
we have a map
\[\ttr{h}=\ttr{h,f,i}\colon h_*i^*\omgs{f}\otimes_{\co_Z}\wnor{i} \lra \co_Y,\]
allowing us to define
\stepcounter{thm}
\begin{equation*}\label{map:-vfinflat-h}\tag{\thethm}
\vttr{h}=\vttr{h,f,i}\colon h_*i^*\omega_{f}\otimes_{\co_Z}\wnor{i} \lra \co_Y
\end{equation*}
as the composite
\[  h_*i^*\omega_{f}\otimes_{\co_Z}\wnor{i} 
\xrightarrow[\phantom{X}\text{via $\bar{\bf v}$}\phantom{X}]{\Iso}
h_*i^*\omgs{f}\otimes_{\co_Z}\wnor{i} \xrightarrow{\phantom{XX}\ttr{h}\phantom{X}} \co_Y.\]

If $\bar{t_i}\in \Gamma(Z\, \I/\I^2)$ is the section generated by the image of $t_i$, then
$\bar{t_1}\wedge\dots\wedge{\bar t}_i$ is a generator of the free rank one $\co_Z$-module
$\wedge^r_{\co_Z}\I/\I^2$. As before, let ${\bf{1/t}}\in \Gamma(Z,\,\wnor{i})=
\Hom_Z(\I/\I^2, \co_Z)$ be the dual generator.
For $\omega\in\Gamma(X,\omega_f)$ let 
$\omega/{\bf t}\in \Gamma(Z,\,i^*\omega_f\otimes_{\co_Z}\wnor{i})$ be the image
of $\omega\otimes {\bf 1/t} \in \Gamma(Z,\,i^*\omega_f)\otimes \Gamma(Z,\,\wnor{i})$.
With these notations, we folow \cite{RD} and \cite{conrad} and define the {\emph{residue symbol}} as
\stepcounter{thm}
\begin{equation*}\label{def:Res-symb}\tag{\thethm}
\Res{X/Y}{\omega}{t_1}{t_r} \set \Gamma(Z,\,\vttr{h})(\omega/{\bf t}) \in \Gamma(Z,\,\co_Z).
\end{equation*}


In \cite[III, \S9]{RD} a list of statement about the residue symbol are made without proof. 
The statements (with minor corrections to the statements in \cite{RD}) 
have been proved by Conrad in \cite[A.2, Appendix A]{conrad}. Since our approach to
residues and the residue symbol follows a different route (via Verdier's isomorphism) we
provide independent proofs of these statements in \S\S\ref{ss:proofs} below. Here are the
statements (R1)--(R10), as in \cite{conrad}, with modifications to take care of our
conventions. In the statements, $\omega$, $f\colon X\to Y$, $Z$, $t_1, \dots, t_r$ are as
above, except in (R4).

\subsection*{(R1)} Let $s_i = \sum_ju_{ij}t_j$ where $u_{ij}\in \Gamma(X, \co_X)$, $1\le i,j \le r$,
and suppose the closed subscheme of $X$ cut out by the $s_i$'s is finite over $Y$.
Then
\[\Res{X/Y}{\omega}{t_1}{t_r}
=\Res{X/Y}{\det{(u_{ij})}\,\omega}{s_1}{s_r}\]

\subsection*{(R2)} (\textit{Localisation}) We give the version in \cite[p.\,239]{conrad}. Suppose
$g\colon X'\to X$ is separated and \'etale,  $Z'=g^{-1}(Z)$, and the map $g'\colon Z'\to Z$
is finite, where $g'$ is induced from $g$. We have a commutative diagram
of schemes
\[
{\xymatrix{
Z' \ar[r]^{i'} \ar[d]_{g'} \ar@{}[dr]|\square & X' \ar[d]^g \\
Z \ar[r]^i \ar[dr]_{h} & X \ar[d]^f\\
& Y
}}
\]
where, as indicated in the diagram, the square on the top is cartesian. Assume
that the function on $Z$ given by $z\mapsto {\mathrm{rank}}_{\co_{Z,z}}{g'}_*(\co_{Z'})_z$
extends to a locally constant function $\mathrm{rk}_{Z'/Z}$ in a Zariski open neighbourhood $V$
of $Z$ in $X$. Then, for every $\omega\in \Gamma(X,\,\omega_f)$, we have,
\[\Res{V/Y}{\omega\cdot\mathrm{rk}_{Z'/Z}}{t_1}{t_r} = \Res{X/Y}{\omega'}{t_1'}{t_r'},\]
where $t_i'=g^*(t_i)\in\Gamma(X',\,\co_{X'})$ and 
$\omega'=g^*(\omega)\in \Gamma(X',\,\omega_{fg})$.

\subsection*{(R3)} (\textit{Restriction}) Suppose we have a commutative diagram of schemes
\[
 {\xymatrix{
 X\, \ar@{^(->}[r]^i \ar[dr]_f & P \ar[d]^\pi \\
 & Y
 }}
 \]
with $f$ smooth and separated 
of relative dimension $r$, $\pi$ smooth and separated
of relative dimension $n=d+r$, $i$
a closed immersion, with $X$ cut out by $s_1,\dots, s_d\in \Gamma(P,\,\co_P)$, and
suppose $t'_1, \dots, t'_r\in \Gamma(P,\,\co_P)$ are such that $s_1, \dots, s_d, t'_1, \dots, t'_r$
cut out a scheme $Z$ which is finite over $Y$, and finally suppose $t_j$ is the restriction of
$t'_j$ to $X$ for $j=1, \dots, r$. Then for every $\nu\in \Gamma(P,\,\Omega^r_{P/Y})$, 
\[
\Res{P/Y}{\wdd{s_1}{s_d}\wedge\nu}{s_1}{s_d, \,t'_1,\dots, t'_r}
=\Res{X/Y}{i^*\nu}{t_1}{t_r}.
\]

\subsection*{(R4)} (\textit{Transitivity}) Suppose we have a pair of separated maps
$X\xrightarrow{f} Y \xrightarrow{g} Z$ are a pair of separated smooth maps, $f$ of
relative dimension $e$ and $g$ of relative dimension $d$. Suppose 
$s_1, \dots, s_d \in \Gamma(Y,\,\co_Y)$ cuts out a scheme $W'$ in $Y$ which
is finite over $Z$, and with $s_j'=f^*(s_j)$, suppose we have $t_1, \dots, t_e\in \Gamma(X,\,\co_X)$
such that $s_1',\dots, s_d', t_1, \dots, t_e$ cut out a scheme $W$ in $X$ which is finite over $Z$.
For $\mu\in \Gamma(\co_Y,\omega_f)$ and $\nu\in \Gamma(\co_X,\,\omega_g)$ we have:
\[
\Res{Y/Z}{\Res{X/Y}{\mu}{t_1}{t_e}\nu}{s_1}{s_d} = 
\Res{X/Z}{\mu\wedge f^*\nu}{t_1}{t_e, \,s_1', \dots, s_d'}.
\]

\subsection*{(R5)} (\textit{Base Change}) Formation of the residue symbol commutes with 
base change. 

\subsection*{(R6)} (\textit{Trace Formula}) For any $\varphi\in\Gamma(X,\,\co_X)$
\[
\Res{X/Y}{\varphi\cdot\wdd{t_1}{t_r}}{t_1}{t_r}=\trc{Z/Y}(\varphi\vert_Z).
\]

\subsection*{(R7)} (\textit{Intersection Formula}) For any collection of positive integers 
$k_1, \dots, k_r$ not all equal to $1$,
\[
\Res{X/Y}{\wdd{t_1}{t_r}}{t_1^{k_1}}{t_r^{k_r}} =0.
\]

\subsection*{(R8)} (\textit{Duality}) (See \cite[p.\,240, (R8)]{conrad}.) If $\omega\vert_Z=0$, then
\[\Res{X/Y}{\omega}{t_1}{t_r}=0.\]
Conversely, let $\{Y_j\}$ be an \'etale covering of $Y$ such that $Y_j$ is affine,
 $Z_j= Z\times_YY_j$ decomposes
into a finite disjoint union of $Z_{jk}$'s with each $Z_{jk}$ contained in an open subscheme
$X_{jk}$ of $X_j\set X\times_YY_j$, with $X_{jk}\cap Z_{jm}=\emptyset$ for $m\neq k$. Also
assume that $\Gamma(X_{jk},\,\co_{X_{jk}})\to \Gamma(Z_{jk},\,\co_{Z_{jk}})$ is 
surjective\footnote{Such $Y_j$'s, $Z_{jk}$'s, and $X_{jk}$'s always exist, using direct limit arguments.}.
If
\[\Res{X_{jk}/Y_j}{f\omega}{t_1}{t_r}=0\]
for all $f\in\Gamma(X_{jk},\,\co_{X_{jk}})$, then $\omega\vert_Z=0$.

\subsection*{(R9)} (\textit{Exterior Differentiation}) For $\nu\in \Gamma(X,\,\Omega^{r-1}_{X/Y})$
and positive integers $k_1, \dots, k_r$,
\[
\Res{X/Y}{\mathrm{d}\nu}{t_1^{k_1}}{t_n^{k_r}} =
\sum_{i=1}^r k_i\cdot\Res{X/Y}{\mathrm{d}t_i\wedge \nu}{t_1^{k_1}}{t_i^{k_i+1}, \dots, t_n^{k_r}}.
\]

\subsection*{(R10)} (\textit{Residue Formula}) Let $h\colon X'\to X$ be a finite map, with $X'$
smooth over $Y$ of  relative dimension $r$. Let $t_j'=h^*(t_j)\in \Gamma(X',\,\co_{X'})$. Then 
\[
\Res{X'/Y}{\nu}{t_1'}{t_r'} = \Res{X/Y}{\vin{h}(\nu)}{t_1}{t_r},
\]
for every $\nu\in \Gamma(X',\,\omega_{fh})$, where $\vin{h}\colon h_*\omega_{fh}\to \omega_f$
is the map in \eqref{map:tr-h-fin}\footnote{We point out that $\vin{h}$ has an explicit description (in
terms of the Kunz-Lipman trace) given in \Tref{thm:kunz-tr1}\,(ii).}.
\subsection{Proofs}\label{ss:proofs}
For a quasi-coherent $\co_X$-module $\eF$, let
\stepcounter{thm}
\begin{equation*}\label{map:v-sh-foo}\tag{\thethm}
\psi=\psi(\eF) \colon h_*(i^*\eF\otimes_{\co_Z}\wI{\co_Z}{\I}) \xrightarrow{\phantom{XXXX}} 
\Rr^r_Zf_*\eF
\end{equation*}
be defined by applying $\Hr^0$ to the composite
\stepcounter{thm}
\begin{equation*}\label{map:v-sh-foo2}\tag{\thethm}
h_*i^\btrg\eF[r] \iso \Rfs i_*i^\btrg\eF[r] \xrightarrow[\eta_i]{\Iso} \Rfs i_*i^\flat\eF[r] \lra
\Rfs \R\iG{Z}\eF[r].
\end{equation*}
where $\eta_i\colon i^\btrg\iso i^\flat$ is the isomorphism in
\cite[(C.2.11)]{fub-abs}.
(See \cite[(5.3.3) and (5.3.4)]{fub-abs}.)
According to \cite[Thm.\,5.3.8]{fub-abs}, the following
diagram commutes
\stepcounter{thm}
\[
\begin{aligned}\label{diag:res-symbol}
{\xymatrix{
 h_*(i^*\omega_{f}\otimes_{\co_Z}\wI{\co_Z}{\I}) \ar[d]_{\psi(\omega_{f})} \ar[rr]^-{\vttr{h}}
 & & \co_Z\ar@{=}[d]\\
\Rr^r_Zf_*\omega_{f} \ar[rr]^{\res{Z}} & & \co_Z 
}}
\end{aligned}\tag{\thethm}
\]

A few things are worth pointing out. First,
 $\mathrm{Res}_{X/Y}\Bigl[\begin{smallmatrix}\omega\\ t_1,\dots,t_r \end{smallmatrix}\Bigr]$
is linear in $\omega$, and since $\vttr{h}$ is a map of sheaves, the residue
symbol is local over $Y$. Moreover,
according to \cite[Remark 5.3.9]{fub-abs}, if $U$ is an open 
subscheme of $X$ containing $Z$, then
\stepcounter{thm}
\begin{equation*}\label{Res:U-X}\tag{\thethm}
\Res{X/Y}{\omega}{t_1}{t_r} = \Res{U/Y}{\omega}{t_1}{t_r}.
\end{equation*}
From \eqref{diag:res-symbol} we see easily that if $Z$ is a disjoint union of $Z_1,\dots, Z_m$
and $X_i$ is open in $X$ with $X_i\cap Z= Z_i$, then as in \cite[p.\,239, (A.1.5)]{conrad}, we have,
\stepcounter{thm}
\begin{equation*}\label{Res:sum}\tag{\thethm}
\Res{X/Y}{\omega}{t_1}{t_r}=\sum_{i=1}^m\Res{X_i/Y}{\omega}{t_1}{t_r}.
\end{equation*}
We also note that by \cite[Thm.\,6.3.2]{fub-abs}, and
\cite[p.\,740,\,Thm.\,2.3.5\,(b)]{cm},
the residue symbol \eqref{def:Res-symb} is stable under arbitrary (noetherian) base change.
This proves (R5).

If $Y=\Spec{\,A}$ and there is an open affine subscheme $U=\Spec{\,R}$ of $X$ containing
$Z$, then by \cite[Prop.\,5.4.4]{fub-abs} and
\eqref{Res:U-X}, we see that 
\stepcounter{thm}
\begin{equation*}\label{Res:resZ}\tag{\thethm}
\Res{X/Y}{\omega}{t_1}{t_r}=\res{Z}\begin{bmatrix}\omega\\t_1, \dots, t_r \end{bmatrix}.
\end{equation*}

Since the formation of the residue symbol is compatible with arbitrary noetherian base change,
i.e., since (R5) is true, we can prove a number of things by assuming $Y$ is the spectrum of
an artin local ring, or of a complete local ring. In greater detail, many of the formulas we
have to prove are of the form $\alpha=\beta$ where $\alpha, \beta\in \Gamma(Y,\co_Y)$. 
It is clearly enough to prove that the germs $\alpha_y$ and $\beta_y$ are equal at every $y\in Y$.
So suppose $y\in Y$ and  $A=\co_y$, the local ring at $y$, and $\fm$ is the maximal  ideal of
$A$. To show $\alpha_y=\beta_y$ it is clearly enough to show 
$\alpha_y\otimes_AA/\fm^n=\beta_y\otimes_AA/\fm^n$ are equal for every $n\in\mathbb{N}$, 
and by (R5), to prove this for a given positive integer
$n$, it is enough to assume $Y=\Spec{A/\fm^n}$. Once we are in this situation, using 
\eqref{Res:U-X}, \eqref{Res:sum}, we are in a situation where \eqref{Res:resZ}
applies. Note that we are in the situation where \eqref{Res:resZ} applies even when
pass to the completion of $A$ with respect to $\fm$. Occasionally, by a further faithful
flat base change on $Y$, we may assume $Y$ is a strictly henselian local ring, or even
a strictly henselian artin local ring.

With this in mind, (R1) follows from 
\cite[Thm.\,5.4.5]{fub-abs},
(R3) from \Cref{cor:res-thm}, (R4) from \Tref{thm:res-res}, (R6) from \Tref{thm:R6}. 
For (R10), first note that $\vin{h}$ is compatible with arbitrary noetherian base
change by \Pref{prop:kunz-sigma}. So once again,
the problem is stable under base change, and we may assume we are in 
a situation where \eqref{Res:resZ} applies. This gives us (R.10) via
\Pref{prop:R10}. We have already seen that (R5) is true.

It remains to prove (R2), (R7), (R8) and (R9).

For (R2), we may assume, as in the proof of (R2) in \cite{conrad}, that $Y$
 is the spectrum of a strictly henselian artin local ring. We are immediately reduced, via
 \eqref{Res:sum}, to the case where $Z$ and $Z'$ consist of a single component each, 
 $Z'=Z$, and $g$ is the identity map. In this case the completion of $X'$ along $Z'$ is the same as the
 completion of $X$ along $Z$, whence, since $\res{Z}$ and $\res{Z'}$ are really only
 dependent on the formal schemes, we are done.  
 
 To prove (R7) we assume without loss of generality that $Y=\Spec{\,A}$, where $A$
 is an artin local ring, that $Z$ is supported at one point, say $z_0$. By shrinking $X$ around $Z$ (via
 \eqref{Res:U-X}) if necessary, we may assume that the map 
 $\pi\colon X\to {\mathbb A}^r_A=\Spec{\,A[T_1,\dots,T_r]}$ 
 defined by ${\bf t}$ is a quasi-finite and that $\pi^{-1}(W)=Z$, where 
 $W$ is the closed subscheme of ${\mathbb A}^r_A$ cut out by $T_1,\dots, T_r$. 
 By Zariski's Main Theorem, 
 we have a finite map $\bar{\pi}\colon \overline{X}\to \mathbb{A}^r_A$ which is a compactification
 of $\pi$, in the sense that there exists an open immersion $u\colon X\to \overline{X}$ such
 that $\bar{\pi}\smcirc u =\pi$. Let 
 $P= {\mathbb A}^r_A \smallsetminus \bar{\pi}(\overline{X}\smallsetminus X)$. Then $P$
 is open in ${\mathbb A}^r_A$, $W\subset P$, and $\bar{\pi}^{-1}(P)\subset X$. Replacing
 $X$ by $\bar{\pi}^{-1}(P)$ if necessary, we may assume $\pi\colon X\to P$ is finite. Shrinking $P$
 around $W$, we may assume $P$ and $X$, are affine, say $P=\Spec{\,D}$ and $X=\Spec{\,E}$. 
 The map $\pi$ is
 flat by \cite[p.\,174, Thm.\,22.3\,($3'$)]{matsumura}, since $Z$ is flat over $Y$, and 
 ${\mathrm{Tor}}^D_1(A, E)=0$ (the latter by noting that
 $K^\bullet({\bf T})\otimes_DE=K^\bullet({\bf t})$).
 By (R10) and 
 \Tref{thm:kunz-tr1}\,(i) we have
 \begin{align*}
 \Res{X/Y}{\wdd{t_1}{t_r}}{t_1^{k_1}}{t_r^{k_r}}
 &= \Res{P/Y}{\trc{Z/W}(1)\cdot\wdd{T_1}{T_r}}{T_1^{k_1}}{T_r^{k_r}} \\
 &=\mathrm{rk}_{B/A}\cdot \Res{P/Y}{\wdd{T_1}{T_r}}{T_1^{k_1}}{T_r^{k_r}}
 \end{align*}
 where $B=\co_{Z,z_0}$.
The last expression is zero if $k_1, \dots, k_r$ are not all equal to $1$, since $W\to Y$
is an isomorphism. This proves (R7)

For (R8), one direction is obvious, namely if $\omega\vert_Z=0$ then
$\mathrm{Res}_{X/Y}\Bigl[\begin{smallmatrix}\omega\\t_1,\dots,t_r \end{smallmatrix}\Bigr] =0$,
for in this case $i^*\omega\otimes{\bf 1/t}=0$. For the ``converse", by faithful flat descent
we may assume $j=1$, $Y=Y_j$, i.e., we may assume $Y=\Spec{\,A}$. Moreover,
via \eqref{Res:U-X} and \eqref{Res:sum}, we may
replace $X$ by $X_{ij}$ if necessary, and assume that
$\Gamma(X,\,\co_X)\to \Gamma(Z,\,\co_Z)$ is surjective. Since $h\colon Z\to Y$ is finite,
$Z$ is affine, say $Z=\Spec{\,B}$. Write $\omega_{B/A}$ for 
$\Gamma(Z,\,i^*\omega_f\otimes\wnor{i})$. The map $\vttr{h}$ induces a natural isomorphism
$\omega_{B/A}\iso \Hom_A(B,\,A)$, which for any $\nu\in\Gamma(X,\,\omega_f)$, sends
 $\nu/ {\bf t}\in \omega_{B/A}$ to $\varphi_\nu\in \Hom_A(B, A)$ where
\[
\varphi_\nu(g) =\Res{X/Y}{\wit{g}\cdot \nu}{t_1}{t_r}\qquad (g\in B)
\]
where $\wit{g}\in \Gamma(X,\,\co_X)$ is any pre-image of $g$ under the surjective map
$\Gamma(X,\co_X)\to \Gamma(Z,\,\co_Z)$. It is clear that under our hypotheses,
$\varphi_\omega =0$, whence the section $\omega/{\bf t}=0$. This means $\omega\vert_Z=0$.

It remains to prove (R9)

\subsubsection*{Proof of (R9)}
As before, we reduce to the case where $Y=\Spec{\,A}$, $A$ an artin local ring, $X=\Spec{\,R}$ and
$Z_{\mathrm{red}} =\{z_0\}$, where $z_0$ is a closed point of $X$ lying over
the closed point of $Y$. We may assume $A$ has an algebraically closed
residue field \cite[$0_{\rm{III}}$, 10.3.1]{ega}. Recall from
\cite[\S\S\,C.5]{fub-abs}, especially
[\textit{ibid.}, (C.5.1)], that for an
$R$-module $M$, we have the notion of a {\emph{stable Koszul complex}}
$K^\bullet_\infty({\bf t},\,M)$. We need this notion for arbitrary sheaves $\eF$
of abelian groups on $X$ (which need not even be $\co_X$-modules). To that
end, let $U_i=\{t_i\neq 0\}=\Spec{\,R_{t_i}}$, $i=1, \dots, r$, $\mathfrak{U}=\{U_i\}$,
and for a sheaf of abelian groups $\eF$ on $X$, 
and $\mathrm{C}^\bullet(\mathfrak{U},\,\eF)$ the {\emph{ordered}} 
\v{C}ech complex associated with $\mathfrak{U}$. 
Since $\check{\Hr}^0(\mathfrak{U},\,\eF)=\eF(U)$, the natural restriction
map $\eF(X)\to \eF(U)$ gives us a complex $K^\bullet_\infty({\bf t},\,\eF)$ defined
as
\[0\lra \eF(X) \lra {\mathrm C}^0(\mathfrak{U},\,\eF) 
\lra {\mathrm C}^1(\mathfrak{U},\,\eF) 
\lra \dots \lra {\mathrm C}^{r-1}(\mathfrak{U},\,\eF) \lra 0\]
with $K^0_\infty({\bf t},\,\eF) = \eF(X)$, $K^{i+1}_\infty({\bf t},\,\eF)=\mathrm{C}^i(\mathfrak{U},\,\eF)$
for $i\ge 0$, the first map being the composite
\[\eF(X) \lra \eF(U) \lra {\mathrm{C}}^0(\mathfrak{U},\,\eF)\] 
and the remaining maps the usual
coboundary maps on \v{C}ech cohomology. If $M$ is an $R$-module, clearly
$K^\bullet_\infty({\bf t},\,\wit{M})= K^\bullet_\infty({\bf t},\,M)$.
Note that $K^\bullet_\infty({\bf t},\,\eF)$ is functorial in $\eF$, as $\eF$ varies over sheaves
of abelian groups on $X$. In what follows, following standard conventions, we write 
$U_{i_1\dots i_p}\set U_{i_1}\cap \dots \cap U_{i_p}$ for $1\le i_1 < \dots < i_p \le r$

Since $\Hr^0(K^\bullet_\infty({\bf t},\,\eF))=\Gamma_Z(X,\,\eF)$,
we have a functorial map of complexes
\stepcounter{thm}
\begin{equation*}\label{map:GZ-K}\tag{\thethm}
\Gamma_Z(X,\,\eF)[0] \lra K^\bullet_\infty({\bf t},\,\eF)
\end{equation*}
which is one readily checks is a quasi-isomorphism when $\eF$ is flasque.
If $\eG^\bullet$ is a complex of {\emph{flasque}} sheaves of abelian groups
on $X$, and $\D(\vert X\vert)$ denotes the derived
category of sheaves of abelian groups on $X$, then \eqref{map:GZ-K} gives
us a pair of isomorphisms in $\D(\vert X\vert)$
\stepcounter{thm}
\begin{equation*}\label{map:RGZ}\tag{\thethm}
\R\Gamma_Z(X,\,\eG^\bullet) \xrightarrow{\phantom{X}\Iso\phantom{X}} 
\Gamma_Z(X,\,\eG^\bullet) \xrightarrow[\eqref{map:GZ-K}]{\Iso} 
{\mathrm{Tot}}(\mathrm{C}^\bullet(\mathfrak{U},\,\eG^\bullet)).
\end{equation*}
The first isomorphism is from general principles (since flasque sheaves have no higher
cohomologies with support), and the second is from the fact that \eqref{map:GZ-K} is
a quasi-isomorphism on flasque sheaves. 

Now suppose $\eF$ is a sheaf of abelian groups on $X$ and $\eF\to \eG^\bullet$
a {\emph{flasque resolution}} of $\eF$. Since $\R\Gamma_Z(X\,\eF) \iso \Gamma_Z(X,\,\eG^\bullet)$, 
\eqref{map:RGZ} gives us an isomorphism
\stepcounter{thm}
\begin{equation*}\label{iso:F-T}\tag{\thethm}
\R\Gamma_Z(X\,\eF) \iso {\mathrm{Tot}}(\mathrm{C}^\bullet(\mathfrak{U},\,\eG^\bullet)),
\end{equation*}
where the right side is the total complex of the double complex 
$\mathrm{C}^\bullet(\mathfrak{U},\,\eG^\bullet)$

By examining the ``columns" of the double
complex $\mathrm{C}^\bullet(\mathfrak{U},\,\eG^\bullet)$ one obtains a map
of complexes
\stepcounter{thm}
\begin{equation*}\label{map:F-col}\tag{\thethm}
K^\bullet_\infty({\bf t},\,\eF) \lra {\mathrm{Tot}}(\mathrm{C}^\bullet(\mathfrak{U},\,\eG^\bullet)).
\end{equation*}

We therefore have a map in $\D(\vert X\vert)$, which is functorial in 
$\eF$ varying over sheaves of abelian groups,
\stepcounter{thm}
\begin{equation*}\label{map:KRGZ}\tag{\thethm}
K^\bullet_\infty({\bf t},\,\eF) \lra \R\Gamma_Z(X,\,\eF)
\end{equation*}
given by $\eqref{map:KRGZ}=\eqref{iso:F-T}^{-1}\smcirc\eqref{map:F-col}$. 

If $\eF$ is {\emph{quasi-coherent}} then \eqref{map:F-col} is a quasi-isomorphism, 
since $U_{i_1\dots i_p}$ and $X$ are affine, whence 
$\eF(U_{i_1\dots i_p })[0]\to \eG^\bullet(U_{i_1\dots i_p })$ and $\eF(X)[0]\to \eG^\bullet(X)$ are quasi-isomorphisms. This means, {\emph{\eqref{map:KRGZ} is an isomorphism when $\eF$ is quasi-coherent}}. In fact, in this case, by definition it agrees with 
\cite[(C.5.2)]{fub-abs}.

Let $\mathrm{d}^{r-1}_{X/Y}\colon \Omega^{r-1}_{X/Y} \to \Omega^r_{X/Y}$ be the standard
exterior derivative map. Note that $\mathrm{d}^{r-1}_{X/Y}$ is not $\co_X$-linear. Nevertheless
our discussion above gives us a commutative diagram:
\stepcounter{thm}
\[
\begin{aligned}\label{diag:GZ-K}
{\xymatrix{
K^\bullet_\infty({\bf{t}}, \Omega^{r-1}_{X/Y}) 
\ar[d]^{\rotatebox{90}{\makebox[0.1cm]{\Iso}}}_{\eqref{map:KRGZ}}
\ar[rr]^{\mathrm{d}^{r-1}_{X/Y}} & &  K^\bullet_\infty({\bf{t}}, \Omega^r_{X/Y})
\ar[d]_{\,\rotatebox{-90}{\makebox[-0.1cm]{\Iso}}}^{\eqref{map:KRGZ}}\\
\R\Gamma_{Z}(X,\,\Omega^{r-1}_{X/Y}) \ar[rr]_{\mathrm{d}^{r-1}_{X/Y}} && 
\R\Gamma_{Z}(X,\,\Omega^r_{X/Y})
}}
\end{aligned}\tag{\thethm}
\]
Using the generalised fraction notation in 
\cite[(C.5.3)]{fub-abs} and the fact that 
\eqref{map:KRGZ} is described for quasi-coherent sheaves by 
\cite[(C.5.2)]{fub-abs}, 
the commutativity of \eqref{diag:GZ-K} gives:
\stepcounter{thm}
\[
\begin{aligned}\label{map:H-d-r-1}
\Hr^r_Z(\mathrm{d}^{r-1}_{X/Y})\begin{bmatrix} \eta\\ 
t_1^{k_1}, \dots, t_r^{k_r} \end{bmatrix}  =
& \begin{bmatrix}\mathrm{d}\eta\\ t_1^{k_1}, \dots, t_r^{k_r} \end{bmatrix} \\
& - \sum_{j=1}^r k_j \begin{bmatrix} \mathrm{d}t_j \wedge \eta\\
t_1^{k_1}, \dots, t_j^{k_j+1}, \dots, t_r^{k_r}\end{bmatrix}.
\end{aligned}\tag{\thethm}
\]
Thus (R9) is equivalent to:
\stepcounter{thm}
\begin{equation*}\label{eq:R9-res}\tag{\thethm}
\res{Z}\smcirc \Hr^r_Z(\mathrm{d}^{r-1}_{X/Y})=0.
\end{equation*}
If $I$ is the ideal of $R$ generated by ${\bf t}$, $R^*$ the completion of $I$ in the $I$-adic
topology, $I^*={\bf t}R^*=IR^*$, $\X=\Spf{(R^*, I^*)}$, then \eqref{eq:R9-res} is equivalent to
\stepcounter{thm}
\begin{equation*}\label{eq:R9-vin}\tag{\thethm}
\vin{\X/Y}\smcirc \Hr^r_{I^*}(\mathrm{d}^{r-1}_{\X/Y})=0,
\end{equation*}
where $\mathrm{d}^{r-1}_{\X/Y}\colon \Omega^{r-1}_{\X/Y}\to \Omega^r_{\X/Y}$ is the exterior
differentiation on the exterior algebra of universally finite differential forms on $\X/Y$.

Since the residue field of $A$ is algebraically closed, the formal $Y$-scheme $\X$ is isomorphic
as a $Y$-scheme to $\Spf{\,A[|T_1, \dots, T_r|]}$, where $A[|T_1, \dots, T_r|]$ is given the
${\bf T}$-adic topology. And using the equivalence of \eqref{eq:R9-res} and \eqref{eq:R9-vin}
the other way, we are done if we prove \eqref{eq:R9-res} for $R=A[{\bf T}]$ and $Z$ the
scheme cut out by ${\bf T}$. In this case, $\eta$ is a finite sum of $(r-1)$-forms of the kind
\[\eta_{j, a_1,\dots, a_r} 
=T_1^{a_1}\dots T_r^{a_r}\cdot{\mathrm{d}}T_1\wedge\dots\wedge\wid{{\mathrm{d}}T_j}
\wedge \dots \wedge {\mathrm{d}}T_r,\] 
where $a_i$ are non-negative integers. Since $Z\to Y$ is an isomorphism in this case,
 $\res{Z}$ is the standard residue
which we know explicitly, the right side of \eqref{map:H-d-r-1} 
$\eta=\eta_{j, a_1,\dots, a_r}$ and $t_i=T_i$, is trivially seen to vanish.

\begin{ack} This paper has been a long time in the making. The outlines were clear as the
 the results in \cite{cm} were being established by the second author (and in fact were the motivation for the results in the last few sections of that paper). Joe Lipman prodded us,
with timely stimulating questions, to write, in fits and starts, little bits of the results we had
been claiming privately. He is the one who encouraged us when the writing slowed down because
of our other commitments. He commented on countless earlier versions of this manuscript
and its first part \cite{fub-abs}. For all of this, and much more, we are very grateful to him. 
\end{ack}



\end{document}